\documentclass[a4paper]{article}
\usepackage{amsmath}
\usepackage{stmaryrd}
\usepackage{amssymb}
\usepackage{latexsym}
\usepackage{algorithm}
\usepackage{algorithmic}
\usepackage{fullpage}
\usepackage{subfigure}
\usepackage{graphics, color}
\usepackage{graphicx}
\usepackage{caption}
\usepackage{multirow}
\usepackage{rotating}

\newtheorem{remark}{Remark}

\begin{document}

\title{Non-intrusive reduced order modeling of parametric electromagnetic scattering problems through Gaussian process regression}

\author{Ying Zhao, Liang Li and Kun Li}

\author{Ying Zhao$^{a}$, \mbox{  } Liang Li$^{a}$,\mbox{  } Kun Li$^{b}$\\
{\small{\it a. School of Mathematical Sciences,}}\\
{\small{\it University of Electronic Science and Technology of China, Chengdu, Sichuan, P. R. China}}\\
{\small{\it b. School of Economic Mathematics, }}\\
{\small{\it Southwestern University of Finance and Economics, Chengdu, P.R. China}}}

\maketitle

\begin{abstract}
This paper is concerned with the design of a non-intrusive model order reduction (MOR) for the system of parametric time-domain Maxwell equations. A time- and parameter-independent reduced basis (RB) is constructed by using a two-step proper orthogonal decomposition (POD) technique from a collection of full-order electromagnetic field solutions, which are generated via a discontinuous Galerkin time-domain (DGTD) solver. The mapping between the time/parameter values and the projection coefficients onto the RB space is approximated by a Gaussian process regression (GPR). Based on the data characteristics of electromagnetic field solutions, the singular value decomposition (SVD) is applied to extract the principal components of the training data of each projection coefficient, and the GPR models are trained for time- and parameter-modes respectively, by which the final global regression function can be represented as a linear combination of these time- and parameter-Gaussian processes. The extraction of the RB and the training of GPR surrogate models are both completed in the offline stage. Then the field solution at any new input time/parameter point can be directly recovered in the online stage as a linear combination of the RB with the regression outputs as the coefficients. In virtue of its non-intrusive nature, the proposed POD-GPR framework, which is equation-free, decouples the offline and online stages completely, and hence can predict the  electromagnetic solution fields at unseen parameter locations quickly and effectively. The performance of our method is illustrated by a scattering problem of a multi-layer dielectric cylinder.

\end{abstract}

\section{Introduction}
In the field of science and engineering, many applications involve the parametric mathematical modeling which often requires solving partial differential equations under a range of parameter values \cite{1,2,3}. In particular, when different frequencies, incident directions or materials are considered, the mathematical models describing electromagnetic wave propagation are parametric, which are usually described by the time-domain Maxwell equations\cite{4, 5}. The development of numerical methods and computer technology in the past decades has enabled effective simulations of electromagnetic radiation/scattering problems, for which the discontinuous Galerkin time-domian (DGTD) solver \cite{6, 7, 8, 9, 10} is one of the most popular tools with easy adaption to complex geometry and material composition, parallelism, and localization. However, a large number of degrees of freedom (DOFs) with a detailed numerical discretization is usually required to guarantee the accuracy of an electromagnetic solution field. Therefore, the huge computational burden in CPU time and memory makes the high-fidelity simulations too costly to allow repeated solution for varying parameters, which motivates the research on the solution measures with reduced cost.

Due to the inherent pattern in the full-order solution fields under parameter variations, researches on the reduced order modeling (ROM) has emerged during the past decades \cite{11, 12} to reduce the computational cost of detailed, high-fidelity simulations. The ROM aims to explore the low-dimensional structure of full-order model without significantly compromising the accuracy, so as to reduce the computational cost. A multitude of ROM techniques, such as the Krylov subspace method based on the Pad\'{e} approximation \cite{13}, the balanced truncation method \cite{14}, and the proper orthogonal decomposition (POD) method \cite{15, 16}, have been developed, and the ROM techniques have been employed in various applications \cite{17}, such as those to the mechanical systems \cite{18, 19, 20}, the control systems \cite{21}, the navigation systems \cite{22, 23}, as well as to the positioning and measurement systems \cite{24}.

Among the existing methods, the data-driven POD proposed by Sirovich \cite{25}, also known as Karhuen-Loéve expansion \cite{26}, principal component analysis \cite{27}, or empirical orthogonal function \cite{28}, is the most widely prevalent in today's big data trend. In the offline stage, ROM of the POD type starts with a collection of high-fidelity snapshot vectors at several time/parameter locations, and a set of POD basis is then extracted by a singular value decomposition (SVD) of the snapshot matrix. A reduced space is then spanned by the POD basis and captures the dominant features of the original full-order system. In the online stage, the approximate solutions at new time/parameter points are recovered as a linear combination of the POD basis vectors in the reduced space.

Conventionally, the combination coefficients of the RB are determined online via a Galerkin projection onto the reduced space \cite{29}. For the simulation of electromagnetic problems, there is substantial work along this line of Galerkin projection, such as the POD-DGTD scheme \cite{30, 31}, the FD-POD method \cite{32}, and the HDG-POD approach \cite{33, 34}. However, such a projection-based framework of a intrusive nature can provide limited computational gain especially when complex nonlinear problems with non-affine dependence on the parameters are involved\cite{35, 36}.
 
 In the emerging growth of the research on data-driven modeling, a valuable alternative to address this concern has been proposed in the recent years. Non-intrusive regression-based RB methods \cite{37, 38, 39} have been developed  based on supervised learning \cite{40, 41}. In these regression-based ROM frameworks, the high-fidelity solver is only used offline as a 'blackbox' to generate the snapshots for the RB construction and the training data for the regression. The RB coefficients at new time/parameter locations are acquired online by means of rapid evaluation of the regression model, guaranteeing a complete decoupling between the online evaluation and the offline training. In the case of parametric time-dependent problems, the regression models approximate the maps between the time/parameter values and the projection coefficients onto the RB space, and are trained from a certain amount of high-fidelity data. Among a number of existing regression models, artificial neural networks (ANNs) have been incorporated into the ROM for nonlinear Poisson equations \cite{42}, steady-state incompressible Navier-Stokes equations \cite{42}, transient flow problems \cite{43}, as well as for compressible fluid problems and combustion \cite{38}. Alternatively, the Gaussian process regression (GPR) \cite{44, 45}, which takes advantage of a Gaussian process prior corrupted by noise and predicts for unseen input locations through the posterior conditioning on the observed input-output pairs \cite{46, 47, 48, 49}, has been efficiently utilized for steady nonlinear structural analysis \cite{50} and parametric time-dependent problems\cite{51}. 
 
In this paper,  to effectively solve the electromagnetic scattering problems described by parametric time-domain Maxwell's equations, we couple a high order DGTD solver with the non-intrusive ROM framework based on the POD and GPR, which sidesteps the conventional Galerkin projection for the reduced order solutions. It is noted for the time-domain electromagnetic problems, the full snapshot matrices by the DGTD method are usually too large for applying the POD. To address this issue, we leverage a two-step POD algorithm \cite{38} with the aid of a division thought to achieve the offline basis extraction: a first-step SVD is carried out on small snapshot matrix under each parameter value, and a second-step SVD is executed on the composite matrix assembled by the basis vectors in the first step. 
 For the electromagnetic simulations, moreover, the RB are usually evolving with different dynamics in the directions of time and parameters, resulting in a difficult global GPR over the time-parameter domain. Therefore, as a practice of \cite{51}, we employ the SVD to decompose the time- and parameter-modes of the training data of each projection coefficient and express the global regression model as the combination of several time- and parameter-dependent Gaussian process surrogates.

The remainder of this paper is arranged as follow. We briefly introduce the DGTD formulations for the time-domain Maxwell's equations in Section 2. The methodology of the POD and the idea of the two-step POD are detailed in Section 3. Section 4 discusses the regression-based approach to the ROM and gives a precise procedure of the global GPR, on the basis of which a more reliable regression strategy using the SVD is derived. Numerical results of electromagnetic scattering problems, show the effectiveness and efficiency of the proposed technique in Section 5. Finally,  we draw conclusions in Section 6.
\section{Mathematical modeling for electromagnetic scattering problems}
\label{sec:DGTD}
We consider the following normalized form of the time-domain Maxwell's equations, which has been widely used to govern unsteady electromagnetic radiation/scattering problems:
\begin{equation}\label{1}
\begin{cases}
\mu_r\displaystyle\frac{\partial \textbf{H}}{\partial t}+\text{curl} \textbf{E} = 0,&\text{in $\Omega \times \mathcal{T}$, }\\[2ex]
\varepsilon_r\displaystyle\frac{\partial \textbf{E}}{\partial t}-\text{curl} \textbf{H} = 0,&\text{in $\Omega \times \mathcal{T}$, }
\end{cases}
\end{equation}
where $\Omega$ is the spatial domain, $\textbf{E} = (E_x,E_y,E_z)^T$ and $\textbf{H} = (H_x,H_y,H_z)^T$ denote the electric field and the magnetic field, respectively, $\mathcal{T}=[0,T_f)$ is the time interval, and $\varepsilon_r$ and $\mu_r$ denote the relative electric permittivity and magnetic permeability parameters, respectively. Details on the governing equation \eqref{1} can be found in \cite{8}. The first-order Silver-M\"{u}ller absorbing boundary condition (ABC) is considered in this paper as:
\begin{equation}\label{2}
\mathcal{L}(\textbf{E},\textbf{H})=\mathcal{L}(\textbf{E}^{inc},\textbf{H}^{inc}), \text{ on } \partial\Omega,
\end{equation}
where $\partial\Omega$ is the boundary of $\Omega$, $\mathcal{L}$ is defined as $\mathcal{L}(\textbf{E},\textbf{H})=\textbf{n} \times \textbf{E}+\mathrm{Z} \textbf{n} \times (\textbf{n}\times\textbf{H})$, $\textbf{n}$ denotes the outer unit normal vector along $\partial\Omega$, $\textbf{E}^{inc}$ and $\textbf{H}^{inc}$ are the incident fields, and $\mathrm{Z}=\sqrt{\frac{\mu_r}{\varepsilon_r}}$. The initial conditions are prescribed as given functions, i.e.,  
\begin{equation}\label{3}
\begin{cases}
\textbf{E}(\textbf{x},0)=\textbf{E}_0(\textbf{x}),\\
\textbf{H}(\textbf{x},0)=\textbf{H}_0(\textbf{x}),\textbf{x}\in\Omega.
\end{cases}
\end{equation}

The goal of this work is to solve equation \eqref{1} with varying parameter $\varepsilon_r$, i.e., $\varepsilon_r \subset \mathcal{P}$ is the characterizing parameter of the model with $\mathcal{P}$ representing the parameter domain.  

In our study, a high order discontinuous Galerkin time-domain (DGTD) method is used to discretize the Maxwell's equations \eqref{1} on unstructured meshes, and will be resorted to generate snapshot data for the forthcoming model order reduction.

The fully discrete scheme of the DGTD technique is given by
\begin{equation}\label{4}
\begin{cases}
\begin{aligned}
&\mathbb{M}^{\varepsilon_r}\frac{\underline{\textbf{E}}_h(t_{n+1}) - \underline{\textbf{E}}_h(t_{n})}{\Delta t} = (\mathbb{K}-\mathbb{S}^i)\underline{\textbf{H}}_h(t_{n+\frac{1}{2}}) - \mathbb{S}^h \underline{\widehat{\textbf{H}}}_h(t_{n+\frac{1}{2}})
-\mathbf{B}^h(n\Delta t),\\[2ex]
&\mathbb{M}^{\mu_r}\frac{\underline{\textbf{H}}_h(t_{n+\frac{3}{2}}) - \underline{\textbf{H}}_h(t_{n+\frac{1}{2}})}{\Delta t} = (-\mathbb{K}+\mathbb{S}^i)\underline{\textbf{E}}_h(t_{n+1})
+\mathbb{S}^e \underline{\widehat{\textbf{E}}}_h(t_{n+1})+\mathbf{B}^e((n+\frac{1}{2})\Delta t),
\end{aligned}
\end{cases}
\end{equation}
in which the time interval $\mathcal{T}=[0,T_f)$ is divided into $N_t$ equally spaced subintervals as $0=t_{0}<t_{1}<\cdots<t_{N_t}=T_f$ with $t_{n}=n\Delta t$ $(n={0,1,\cdots,N_t})$, $\Delta t $ denoting the time step size. Here $\mathbb{M}^{\varepsilon_r}$ and $\mathbb{M}^{\mu_r}$ are the mass matrices, $\mathbb{K}$ is the stiffness matrix, $\mathbb{S}^i$ is the surface matrix for the interior faces, and $\mathbb{S}^h$ and $\mathbb{S}^e$ are the boundary face matrices. The specific process of the DGTD discretization and the detailed definition of these matrices can be found in \cite{32}. 

\section{Proper orthogonal decomposition}
\label{sec:POD}
In this section, we introduce the POD approach for model order reduction. Given a parameter sampling $\mathcal{P}_h=\{\theta_1,\theta_2,\cdots,\theta_{N_\theta}\}$ over the parameter domain $\mathcal{P}$, a collection of high-fidelity solutions of \eqref{1} can be obtained by running the DGTD solver under different parameter points in $\mathcal{P}_h$. In this work, we equidistantly extract $N_l$ transient solutions $\underline{\textbf{u}}_h(t_{n_i},\theta_j)$ ($0\leq t_{n_1} < \cdots < t_{n_{N_l}} \leq t_{N_t}$, $n_{N_l}\ll N_t$) from the solution  $\underline{\textbf{u}}_h(t_i,\theta_j)=\underline{\textbf{E}}_h(t_i,\theta_j)~\text{or}~\underline{\textbf{H}}_h(t_{i+\frac{1}{2}},\theta_j)$ for each parameter value in $\mathcal{P}_h$. We formulate  $N_h\times N_l$ snapshot matrices with respect to each parameter sampling point, i.e., 
$$
\mathbb{S}_{\textbf{u},j}=\left(
\begin{smallmatrix}
\underline{\textbf{u}}_{h,1}(t_{n_1},\theta_j)&\underline{\textbf{u}}_{h,1}(t_{n_2},\theta_j)&\cdots &\underline{\textbf{u}}_{h,1}(t_{n_{N_l}},\theta_j)\\
\underline{\textbf{u}}_{h,2}(t_{n_1},\theta_j)&\underline{\textbf{u}}_{h,2}(t_{n_2},\theta_j)&\cdots &\underline{\textbf{u}}_{h,2}(t_{n_{N_l}},\theta_j)\\
\vdots&\vdots&\ddots&\vdots\\
\underline{\textbf{u}}_{h,N_h}(t_{n_1},\theta_j)&\underline{\textbf{u}}_{h,N_h}(t_{n_2},\theta_j)&\cdots &\underline{\textbf{u}}_{h,N_h}(t_{n_{N_l}},\theta_j)\\
\end{smallmatrix}
\right),\mathbf{u} = \mathbf{E}~\text{or}~\mathbf{H}, j=1,2,\cdots,N_\theta,
$$
$N_h$ being the number of DOFs, and the global snapshot matrices that assemble all of the $\mathbb{S}_{\textbf{u},j}$'s is
$$
\mathbb{S}_{\textbf{u}}=[\mathbb{S}_{\textbf{u},1}|\cdots|\mathbb{S}_{\textbf{u},N_\theta}]=\left(
\begin{smallmatrix}
\underline{\textbf{u}}_{h,1}(t_{n_1},\theta_1)&\cdots &\underline{\textbf{u}}_{h,1}(t_{n_{N_l}},\theta_1)&\cdots &\underline{\textbf{u}}_{h,1}(t_{n_1},\theta_{N_\theta})&\cdots &\underline{\textbf{u}}_{h,1}(t_{n_{N_l}},\theta_{N_\theta})\\
\underline{\textbf{u}}_{h,2}(t_{n_1},\theta_1)&\cdots &\underline{\textbf{u}}_{h,2}(n_{N_l},\theta_1)&\cdots &\underline{\textbf{u}}_{h,2}(t_{n_1},\theta_{N_\theta})&\cdots &\underline{\textbf{u}}_{h,2}(t_{n_{N_l}},\theta_{N_\theta})\\
\vdots&\ddots&\vdots&\cdots&\vdots&\ddots&\vdots\\
\underline{\textbf{u}}_{h,N_h}(t_{n_1},\theta_1)&\cdots &\underline{\textbf{u}}_{h,N_h}(t_{n_{N_l}},\theta_1)&\cdots &\underline{\textbf{u}}_{h,N_h}(t_{n_1},\theta_{N_\theta})&\cdots &\underline{\textbf{u}}_{h,N_h}(t_{n_{N_l}},\theta_{N_\theta})\\
\end{smallmatrix}
\right),
$$
which is an $N_h\times N_s$ matrix with $N_s=N_\theta\cdot N_l$.

We then perform a low-rank approximation to $\mathbb{S}_{\textbf{u}}$ and construct a low dimensional vector space $\mathcal{V}_{\textbf{u},rb}$ with reduced dimension $d_{\textbf{u}}\ll min\{N_h,N_s\}$, with which we aim to effectively capture the feature of the solution manifold over the parameter variation. Spanned by a group of time- and parameter-independent RB functions, the reduced space is given as
$$
\mathcal{V}_{\textbf{u},rb}=span\{\psi_{\textbf{u},1},\psi_{\textbf{u},2},\cdots,\psi_{\textbf{u},d_{\textbf{u}}}\}, \mathbf{u} = \mathbf{E}~\text{or}~\mathbf{H}.
$$
Thus, the reduced-order solution $\underline{\textbf{u}}_h^d(t,\theta)$ serves as an approximation to the high-fidelity solution $\underline{\textbf{u}}_h(t,\theta)$ and can be represented as
\begin{equation}\label{5}
\underline{\textbf{u}}_h^d(t,\theta)=\sum_{i=1}^{d_{\textbf{u}}} \alpha_{\textbf{u},i}(t,\theta)\psi_{\textbf{u},i}, \mathbf{u} = \mathbf{E}~\text{or}~\mathbf{H},
\end{equation}
 where $\alpha_{\textbf{u}}(t,\theta)=[\alpha_{\textbf{u},1}(t,\theta),\alpha_{\textbf{u},2}(t,\theta),\cdots,\alpha_{\textbf{u},d_{\textbf{u}}}(t,\theta)]^T\in\mathbb{R}^{d_{\textbf{u}}}$ collects the combination coefficients.
 
 We perform SVD to $\mathbb{S}_{\textbf{u}}$ and let
 $$
 \mathbb{W}_\textbf{u}^T\mathbb{S}_\textbf{u}\mathbb{Z}_\textbf{u}=\left(
\begin{matrix}
\Sigma^{\textbf{u}}_{r_{\textbf{u}}\times r_{\textbf{u}}}&0_{r_{\textbf{u}}\times (N_s-r_{\textbf{u}})}\\
0_{(N_h-r_{\textbf{u}})\times r_{\textbf{u}} }&0_{(N_h-r_{\textbf{u}})\times (N_s-r_{\textbf{u}})}
\end{matrix}
\right)=\mathbb{D}_{\textbf{u}}, \mathbf{u} = \mathbf{E}~\text{or}~\mathbf{H},
 $$
 where $\mathbb{W}_\textbf{u}= (\textbf{w}_{\textbf{u},1},\textbf{w}_{\textbf{u},2},\cdots,\textbf{w}_{\textbf{u},N_h})$ and $\mathbb{Z}_\textbf{u}=(\textbf{z}_{\textbf{u},1},\textbf{z}_{\textbf{u},2},\cdots,\textbf{z}_{\textbf{u},N_s})$ are $N_h\times N_h$ and $N_s\times N_s$ unitary matrices, respectively, $\sum^{\textbf{u}}_{r_{\textbf{u}}\times r_{\textbf{u}}}=diag(\sigma_{\textbf{u},1},\sigma_{\textbf{u},2},\cdots,\sigma_{\textbf{u},r_{\textbf{u}}})$ with $\sigma_{\textbf{u},1}\geq \sigma_{\textbf{u},2}\geq  \cdots \geq  \sigma_{\textbf{u},r_{\textbf{u}}}\geq0$ being the singular values of $\mathbb{S}_\textbf{u}$, and hence $r_{\textbf{u}}$ is the rank of $\mathbb{S}_\textbf{u}$. According to the Schmidt-Eckart-Young theorem \cite{52, 53}, the POD basis of dimension $d_{\textbf{u}}$ $(d_{\textbf{u}}<r_{\textbf{u}})$ is the set $\{\psi_{\textbf{u},i}\}_{i=1}^{d_{\textbf{u}}}$ with $\psi_{\textbf{u},i}=\textbf{w}_{\textbf{u},i}$, which can minimize the projection error of the snapshots among all $d_{\textbf{u}}$-dimensional orthogonal bases in $\mathbb{R}^{N_h}$. The error bound can be evaluated using the singular values
\begin{equation}\label{6}
 \sum_{i=1}^{N_s}\|\mathbb{S}_{\textbf{u}}(:,i)-\sum_{j=1}^{d_{\textbf{u}}}(\mathbb{S}_{\textbf{u}}(:,i),\psi_{\textbf{u},j})\psi_{\textbf{u},j}\|^2_{\mathbb{R}^{N_h}}=\
 \sum_{i=1}^{N_s}\|\mathbb{S}_{\textbf{u}}(:,i)-\Psi_\textbf{u}\Psi_\textbf{u}^T\mathbb{S}_{\textbf{u}}(:,i)\|^2_{\mathbb{R}^{N_h}}=\
 \sum_{j=d_{\textbf{u}}+1}^{r_{\textbf{u}}}\sigma_{\textbf{u},j}^2,
\end{equation}
where $\Psi_\textbf{u}=[\psi_{\textbf{u},1},\cdots,\psi_{\textbf{u},d_{\textbf{u}}}]$. One can determine the dimension $d_{\textbf{u}}$ to be the smallest integer such that $d_{\textbf{u}}=\arg\max{\{ \mathcal{E}}(d_{\textbf{u}}):{\mathcal{E}}(d_{\textbf{u}})\geq1-\epsilon\}$ with $\mathcal{E}(d_{\textbf{u}})=\sum_{i=1}^{d_{\textbf{u}}}\sigma_{\textbf{u},i}^2/\sum_{i=1}^{r_{\textbf{u}}}\sigma_{\textbf{u},i}^2$, and $\epsilon$ being the relative error tolerance controlling the accuracy of POD.

Since $N_s$ is large and the SVD of such a large-scale snapshot matrix is expensive, the POD algorithm described above can not be directly applied to the parametric electromagnetic problems with a large number of time steps. To overcome this difficulty, we adopt a two-step POD strategy which can effectively save the computational cost. The process of such a two-step POD is given as follows (shown in Algorithm 1):
\begin{enumerate}
\item POD for the small snapshot matrix for each single parameter point. For each $\mathbb{S}_{\textbf{u},j}, j=1,2,\cdots,N_\theta, \mathbf{u} = \mathbf{E}~\text{or}~\mathbf{H}$, we run the POD process with a relative error tolerance $\epsilon_{\textbf{u},t}$ and obtain a reduced basis $\{\gamma^j_{\textbf{u},1},\cdots,\gamma^j_{\textbf{u},d_{\textbf{u}}^j}\}$ with $d_{\textbf{u}}^j$ vectors, followed by the assembly of a matrix $\mathbb{T}_{\textbf{u},j}=[\gamma^j_{\textbf{u},1},\cdots,\gamma^j_{\textbf{u},d_{\textbf{u}}^j}].$
\item POD for the composite matrix. Put together the matrices obtained in the first step $\mathbb{T}_{\textbf{u},j},j=1,2,\cdots,N_\theta$, to construct a composite matrix $\mathbb{T}_{\textbf{u}}=[\mathbb{T}_{\textbf{u},1}|\cdots|\mathbb{T}_{\textbf{u},N_\theta}]$, and perform the POD on $\mathbb{T}_{\textbf{u}}$ with relative error tolerance $\epsilon_{\textbf{u},\theta }$. Then we get a reduced basis $\{\psi_{\textbf{u},1},\cdots,\psi_{\textbf{u},d_{\textbf{u}}}\}$ comprised of $d_{\textbf{u}}$ basis vectors collected in a matrix $\Psi_\textbf{u}=[\psi_{\textbf{u},1},\cdots,\psi_{\textbf{u},d_{\textbf{u}}}]$.
\end{enumerate}

\renewcommand{\algorithmicrequire}{\textbf{Input:}}
\renewcommand{\algorithmicensure}{\textbf{Output:}}
\begin{algorithm}
\caption{Two-step Proper Orthogonal Decomposition}
\label{alg:A}
\begin{algorithmic}[1]
\REQUIRE ~~\\
Snapshot matrices $\mathbb{S}_{\textbf{u},j}$, $j=1,2,\cdots,N_\theta$, $\mathbf{E}~\text{or}~\mathbf{H}$;\\
Projection error tolerance $\epsilon_{\textbf{u},t}$ and $\epsilon_{\textbf{u},\theta }$;
\ENSURE ~~\\
Reduced number of DOFs $d_{\textbf{u}}$;
Matrix $\Psi_\textbf{u}$ collecting the reduced basis;
\FOR{$j=1,2,\cdots,N_\theta$}
\STATE Perform the SVD: $[\mathbb{W}^j_{\textbf{u}},\mathbb{D}^j_{\textbf{u}},\mathbb{Z}^j_{\textbf{u}}]=\emph{svd}(\mathbb{S}_{\textbf{u},j})$;
\STATE Determine the dimension $d_{\textbf{u}}^j=\textmd{argmin}{\{\mathcal{E}}(d_{\textbf{u}}^j):{\mathcal{E}}(d_{\textbf{u}}^j)\geq1-\epsilon_{\textbf{u},t}\}$;
\STATE Set the basis as: $\gamma^j_{\textbf{u},i}=\textbf{w}^j_{\textbf{u},i}, i=1,\cdots,d_{\textbf{u}}^j$, get $\mathbb{T}_{\textbf{u},j}=[\gamma^j_{\textbf{u},1},\cdots,\gamma^j_{\textbf{u},d_{\textbf{u}}^j}]$;
\ENDFOR
\STATE Assemble $\mathbb{T}_{\textbf{u}}=[\mathbb{T}_{\textbf{u},1}|\cdots|\mathbb{T}_{\textbf{u},N_\theta}]$;
\STATE Perform the SVD: $[\mathbb{W}_{\textbf{u}},\mathbb{D}_{\textbf{u}},\mathbb{Z}_{\textbf{u}}]=\emph{svd}(\mathbb{T}_{\textbf{u}})$;
\STATE Determine the dimension $d_{\textbf{u}}=\textmd{argmin}{\{\mathcal{E}}(d_{\textbf{u}})\geq1-\epsilon_{\textbf{u},\theta }\}$;
\STATE Set the basis as $\psi_{\textbf{u},i}=\textbf{w}_{\textbf{u},i}, i=1,\cdots,d_{\textbf{u}}$, and $\Psi_\textbf{u}=[\psi_{\textbf{u},1},\cdots,\psi_{\textbf{u},d_{\textbf{u}}}]$;
\end{algorithmic}
\end{algorithm}

According to the algebraic projection theory, the reduced-order approximation of the field takes the form
\begin{equation}\label{7}
\underline{\textbf{u}}_h(t,\theta)\approx\underline{\textbf{u}}_h^d(t,\theta)\
=\Psi_{\textbf{u}}\Psi_{\textbf{u}}^T\underline{\textbf{u}}_h(t,\theta)=\Psi_{\textbf{u}}\alpha_{\textbf{u}}(t,\theta), \textbf{u}=\mathbf{E}~\text{or}~\mathbf{H},
\end{equation}
in which $\alpha_{\textbf{u}}(t,\theta)=\Psi_{\textbf{u}}^T\underline{\textbf{u}}_h(t,\theta)$ collects combination coefficients of the RB.

 \begin{remark} According to \eqref{6}, the error bound in the first and second steps of the two-step POD algorithm are written as
\begin{equation}\label{8}
\begin{cases}
\begin{aligned}
&\sum_{i=1}^{N_l}\|\mathbb{S}_{\textbf{u},j}(:,i)-\mathbb{T}_{\textbf{u},j}\mathbb{T}_{\textbf{u},j}^T\mathbb{S}_{\textbf{u},j}(:,i)\|^2_{\mathbb{R}^{N_h}}=\sum_{i=d_{\textbf{u}}^j+1}^{r_{\textbf{u}}^j}(\sigma_{\textbf{u},i}^j)^2\leqslant \epsilon_{\textbf{u},t}\sum_{i=1}^{r_{\textbf{u}}^j}(\sigma_{\textbf{u},i}^j)^2, 1\leqslant j \leqslant N_\theta,\\
&\sum_{j=1}^{N_\theta} \sum_{i=1}^{d_{\textbf{u}}^j}\|\mathbb{T}_{\textbf{u},j}(:,i)-\Psi_\textbf{u}\Psi_\textbf{u}^T\mathbb{T}_{\textbf{u},j}(:,i)\|^2_{\mathbb{R}^{N_h}}=\sum_{i=d_{\textbf{u}}+1}^{r_{\textbf{u}}}(\sigma_{\textbf{u},i})^2\leqslant \epsilon_{\textbf{u},\theta}\sum_{i=1}^{r_{\textbf{u}}}(\sigma_{\textbf{u},i})^2, 
\end{aligned}
\end{cases}
\end{equation}
respectively, where $r_{\textbf{u}}$ and $r_{\textbf{u}}^j$ ($j=1,2,\cdots,N_\theta, \textbf{u}=\mathbf{E}~\text{or}~\mathbf{H}$) are the rank of $\mathbb{S}_\textbf{u}$ and $\mathbb{S}_{\textbf{u},j}$, and $\{\sigma_{\textbf{u},i}\}_{i=1}^{r_{\textbf{u}}}$ and $\{\sigma_{\textbf{u},i}^j\}_{i=1}^{r_{\textbf{u}}^j}$ are the corresponding  singular values. The two-step POD projection error can thus be bounded as
\begin{equation}\label{9}
\begin{split}
\sum_{j=1}^{N_\theta} \sum_{i=1}^{N_l}\|\mathbb{S}_{\textbf{u},j}(:,i)-\Psi_\textbf{u}\Psi_\textbf{u}^T\mathbb{S}_{\textbf{u},j}(:,i)\|&_{\mathbb{R}^{N_h}}\leqslant \sum_{j=1}^{N_\theta}\sum_{i=1}^{N_l}\|\mathbb{S}_{\textbf{u},j}(:,i)-\mathbb{T}_{\textbf{u},j}\mathbb{T}_{\textbf{u},j}^T\mathbb{S}_{\textbf{u},j}(:,i)\|_{\mathbb{R}^{N_h}}\\
&+\sum_{j=1}^{N_\theta}\sum_{i=1}^{N_l}\|\mathbb{T}_{\textbf{u},j}\mathbb{T}_{\textbf{u},j}^T\mathbb{S}_{\textbf{u},j}(:,i)-\Psi_\textbf{u}\Psi_\textbf{u}^T\mathbb{T}_{\textbf{u},j}\mathbb{T}_{\textbf{u},j}^T\mathbb{S}_{\textbf{u},j}(:,i)\|_{\mathbb{R}^{N_h}}\\
&+\sum_{j=1}^{N_\theta}\sum_{i=1}^{N_l}\|\Psi_\textbf{u}\Psi_\textbf{u}^T\mathbb{T}_{\textbf{u},j}\mathbb{T}_{\textbf{u},j}^T\mathbb{S}_{\textbf{u},j}(:,i)-\Psi_\textbf{u}\Psi_\textbf{u}^T\mathbb{S}_{\textbf{u},j}(:,i)\|_{\mathbb{R}^{N_h}}\\
&\leqslant \sqrt{\epsilon_{\textbf{u},t}}\mathcal{L}_1+\sqrt{\epsilon_{\textbf{u},\theta}}\mathcal{L}_2, \textbf{u}=\mathbf{E}~\text{or}~\mathbf{H},
\end{split}
\end{equation}
in which we set $N_\textbf{u}=\sum_{j=1}^{N_\theta}d_{\textbf{u}}^j$ and 
$$
\begin{cases}
\mathcal{L}_1=(1+\|\Psi_\textbf{u}\Psi_\textbf{u}^T\|_F)\sum_{j=1}^{N_\theta}(N_l\sum_{i=1}^{r_{\textbf{u}}^j}(\sigma_{\textbf{u},i}^j)^2)^\frac{1}{2},\\
\mathcal{L}_2=\max_{1\leqslant j \leqslant N_\theta}\sum_{i=1}^{N_l}\|\mathbb{S}_{\textbf{u},j}(:,i)\|_{\mathbb{R}^{N_h}}\max_{1\leqslant j \leqslant N_\theta, 1\leqslant i\leqslant d_{\textbf{u}}^j}\|\mathbb{T}_{\textbf{u},j}(:,i)\|_{\mathbb{R}^{N_h}}(N_\textbf{u}\sum_{i=1}^{r_{\textbf{u}}}(\sigma_{\textbf{u},i})^2)^\frac{1}{2}.
\end{cases}
$$
Therefore, the accuracy of the two-step POD can be controlled by the tolerance $\epsilon_{\textbf{u},t}$ and $\epsilon_{\textbf{u},\theta}$.
\end{remark}
\section{Gaussian process regression}
\label{sec:GPR}
Based on the RB functions obtained by POD, one only needs to calculate the combination coefficients $\alpha_{\textbf{u}}(t^{\ast},\theta^{\ast})$ ($\textbf{u}=\mathbf{E}~\text{or}~\mathbf{H}$) to obtain an approximate solution for any new time/parameter point $(t^{\ast},\theta^{\ast})$. Galerkin projection method is a usually used to determine the combination coefficients of the RB \cite{4,luo2016podfdtd}. However, the Galerkin projection scheme is complicated and relatively expensive because one needs to go back to the continuous function spaces and compute the inverse of some dense matrices, which limits the application of POD method. Therefore, in the proposed non-intrusive reduced-order framework, a regression model, \emph{i.e.} the Gaussian process regression (GPR) model, is used to calculate the approximate solution at new time/parameter points.  The simulation is decomposed into offline/online stages. The offline stage includes the extraction of the RB, the establishment and training of the GPR. In the online stage one only need to calculate the output of the GPR and compute the linear combination of the RB, which are cheap.

\subsection{Regression-based approach to reduced-order solutions}
To approximate the projection coefficients $\alpha_{\textbf{u}}(t,\theta)$ ($\textbf{u}=\mathbf{E}~\text{or}~\mathbf{H}$) for any desired time-parameter location $(t,\theta)\in\mathcal{T}\times\mathcal{P}$, we resort to the technique of nonlinear regression $\widehat{\alpha}_{\textbf{u}}$:
\begin{equation}\label{10}
(t,\theta)\longmapsto\alpha_{\textbf{u}}(t,\theta)=\Psi_{\textbf{u}}^T\underline{\textbf{u}}_h(t,\theta)\approx\widehat{\alpha}_{\textbf{u}}(t,\theta), \textbf{u}=\mathbf{E}~\text{or}~\mathbf{H},
\end{equation}
and the regression models $\widehat{\alpha}_{\textbf{u}}(\cdot,\cdot)$ are constructed from a set of training data
\begin{equation}\label{11}
\begin{cases}
D_\textbf{u} = \{\{(t,\theta),\Psi_{\textbf{u}}^T\underline{\textbf{u}}_h(t,\theta)\}:t\in\mathcal{T}_{tr},\theta\in\mathcal{P}_{tr}\}, \textbf{u}=\mathbf{E}~\text{or}~\mathbf{H},\\
\mathcal{T}_{tr}=\{t^{n_i}:i=1,2,\cdots,N_t^{tr}\}\subset \mathcal{T},\\
\mathcal{P}_{tr}=\{\theta^j:j=1,2,\cdots,N_\theta^{tr}\}\subset \mathcal{P}.
\end{cases}
\end{equation}

These models are then used during the online stage to recover the output $\widehat{\alpha}_{\textbf{u}}(t^{\ast},\theta^{\ast})$ ($\textbf{u}=\mathbf{E}~\text{or}~\mathbf{H}$) for any new input $(t^{\ast},\theta^{\ast})\in\mathcal{T}\times\mathcal{P}$. The corresponding reduced-order solution is written as
\begin{equation}\label{12}
\underline{\textbf{u}}_{h,reg}^d(t,\theta)=\Psi_{\textbf{u}}\widehat{\alpha}_{\textbf{u}}(t^{\ast},\theta^{\ast})=\
\sum_{i=1}^{d_{\textbf{u}}}\widehat{\alpha}_{\textbf{u},i}(t^{\ast},\theta^{\ast})\psi_{\textbf{u},i}, \textbf{u}=\mathbf{E}~\text{or}~\mathbf{H}.
\end{equation}

In this paper, we take advantage of the GPR to realize this regression-based approach.
\subsection{Gaussian process regression}
Regression is a supervised machine learning method, which is used to predict some continuous quantities by using a certain amount of observation data. Define $D_{\textbf{u},l}=\{(\textbf{x}^{(i-1)N_\theta^{tr}+j},y_{\textbf{u},l}^{(i-1)N_\theta^{tr}+j}):i=1,2,\cdots,N_t^{tr}, j=1,2,\cdots,N_\theta^{tr}\}=(\textbf{X},\textbf{y}_{\textbf{u},l})$ ($\textbf{u}=\mathbf{E}~\text{or}~\mathbf{H}$), in which $\textbf{x}^{(i-1)N_\theta^{tr}+j}=(t^{n_i},\theta^j)^T \in \mathcal{D}$ represents the input time-parameter vector, $\mathcal{D}$ is the domain of time-parameter inputs, $y_{\textbf{u},l}^{(i-1)N_\theta^{tr}+j}=\psi_{\textbf{u},l}^T\underline{\textbf{u}}_h(t^{n_i},\theta^j) \in \mathbb{R}$ is the corresponding $l$th projection coefficient scalar with $l=1,2,\cdots,d_\textbf{u}$, $\textbf{X}=[\textbf{x}^1,\textbf{x}^2,\cdots,\textbf{x}^{N_t^{tr}N_\theta^{tr}}]$ and $\textbf{y}_{\textbf{u},l}=[y^1_{\textbf{u},l},y^2_{\textbf{u},l},\cdots,y^{N_t^{tr}N_\theta^{tr}}_{\textbf{u},l}]$ collect the input and output matrices respectively. The task of the regression is to study the mapping between $\textbf{X}$ and $\textbf{y}_{\textbf{u},l}$, thereby predicting the most likely output value $y_{\textbf{u},l}^\ast=\psi_{\textbf{u},l}^T\underline{\textbf{u}}_h(t^\ast,\theta^\ast)$ for a new test point $\textbf{x}^\ast=(t^\ast,\theta^\ast)^T$.

A Gaussian process (GP) assumes that the random variables at any finite set of input locations have a joint Gaussian distribution, and its properties are completely determined by the mean function and the covariance function defined as
\begin{equation}\label{13}
\begin{cases}
m(\textbf{x})=\textit{\textbf{E}}[f(\textbf{x})],\\
\kappa(\textbf{x},\textbf{x}^{'})=\textit{\textbf{E}}[(f(\textbf{x})-m(\textbf{x}))(f(\textbf{x}^{'})-m(\textbf{x}^{'}))],
\end{cases}
\end{equation}
where $(\textbf{x},\textbf{x}^{'})\in\mathcal{D}\times\mathcal{D}$ and the function $f$ here represents a GP. Hence a GP is defined as
\begin{equation}\label{14}
f(\textbf{x})\sim GP(m(\textbf{x}),\kappa(\textbf{x},\textbf{x}^{'})).
\end{equation}

The essence of Gaussian process regression (GPR) is to infer the relationship between the input variable $\textbf{x}$ and the output $y_{\textbf{u},l}$, that is, to determine the conditional distribution of the target output after the input variable is given. In GPR, it is assumed that the prior regression function is a GP $f$ corrupted by an independent Gaussian noise
\begin{equation}\label{15}
y_{\textbf{u},l} = f(\textbf{x})+\epsilon , \epsilon\sim\mathcal{N}(0,\sigma^2_y), \textbf{u}=\mathbf{E}~\text{or}~\mathbf{H}.
\end{equation}

Based on a finite number of training data, we can get the prior joint Gaussian distribution of the observed values as
\begin{equation}\label{16}
\textbf{y}_{\textbf{u},l}|\textbf{X} \sim \mathcal{N}(m(\textbf{X}),\textbf{K}_{\textbf{u},l}), \textbf{K}_{\textbf{u},l}=\textrm{cov}[\textbf{y}_{\textbf{u},l}|\textbf{X}]=\kappa(\textbf{X},\textbf{X})+\sigma^2_{y}\textbf{I}_n, \textbf{u}=\mathbf{E}~\text{or}~\mathbf{H},
\end{equation}
where $\textbf{I}_n$ is the $n\times n$ identity matrix and $n=N_t^{tr}\cdot N_\theta^{tr}$.

Given a new test input denoted by $\textbf{x}^{\ast}$, predictions of the corresponding noise-free outputs $y^{\ast}_{\textbf{u},l}$ is desired. Under the Bayesian principle, the joint density of the observed outputs $\textbf{y}_{\textbf{u},l}$ and the noise-free test output $y^{\ast}_{\textbf{u},l}$ can be written as
\begin{equation}\label{17}
\left[
\begin{matrix}
\textbf{y}_{\textbf{u},l}\\
y^{\ast}_{\textbf{u},l}
\end{matrix}
\right]\sim \mathcal{N}\left(\begin{matrix}
\left[
\begin{matrix}
m(\textbf{X})\\
m(\textbf{x}^{\ast})
\end{matrix}
\right],\left[
\begin{matrix}
\textbf{K}_{\textbf{u},l}&\textbf{K}_{\textbf{u},l}^{\ast}\\
(\textbf{K}_{\textbf{u},l}^{\ast})^T&K_{\textbf{u},l}^{\ast\ast}
\end{matrix}
\right]
\end{matrix}
\right), \textbf{u}=\mathbf{E}~\text{or}~\mathbf{H},
\end{equation}
where $\textbf{K}_{\textbf{u},l}^{\ast}=\kappa(\textbf{X},\textbf{x}^{\ast})$ and $K_{\textbf{u},l}^{\ast\ast}=\kappa(\textbf{x}^{\ast},\textbf{x}^{\ast})$,
and the posterior predictive distribution for $y^{\ast}_{\textbf{u},l}$ can be obtained following the standard rules for conditional Gaussian as
\begin{equation}\label{18}
y^{\ast}_{\textbf{u},l}|\textbf{x}^{\ast},\textbf{X},\textbf{y}_{\textbf{u},l} \sim \mathcal{N}(m(\textbf{x}^{\ast})+(\textbf{K}^{\ast}_{\textbf{u},l})^T\textbf{K}_{\textbf{u},l}^{-1}(\textbf{y}_{\textbf{u},l}-m(\textbf{X})),K^{\ast \ast}_{\textbf{u},l}-(\textbf{K}^{\ast}_{\textbf{u},l})^T\textbf{K}_{\textbf{u},l}^{-1}\textbf{K}^{\ast}_{\textbf{u},l}).
\end{equation}
It is easy to verify that the corresponding posterior process can be represented as
\begin{equation}\label{19}
\begin{split}
&y^{\ast}_{\textbf{u},l}|D_{\textbf{u},l} \sim GP(m^{\ast},C^{\ast}), \textbf{u}=\mathbf{E}~\text{or}~\mathbf{H},\\
&m^{\ast}(\textbf{x})=m(\textbf{x})+\kappa(\textbf{x},\textbf{X})\textbf{K}_{\textbf{u},l}^{-1}(\textbf{y}_{\textbf{u},l}-m(\textbf{X})),\\
&C^{\ast}(\textbf{x},\textbf{x}^{'})=\kappa(\textbf{x},\textbf{x}^{'})-\kappa(\textbf{x},\textbf{X})\textbf{K}_{\textbf{u},l}^{-1}\kappa(\textbf{X},\textbf{x}^{'}).
\end{split}
\end{equation}

The selection of mean function and covariance function plays a key role in the final prediction. For the mean function, we select the constant function $m(\textbf{x}):=\sum_{i=1}^{n}\beta_i$ with $\beta_i$ ($i=1,2,\cdots,n$) being {\color{red}constant variables}. For the covariance function, a frequently used one is the automatic relevance determination (ARD) squared exponential (SE) kernel:
\begin{equation}\label{20}
\kappa(\textbf{x},\textbf{x}^{'})=\sigma^2_f\exp(-\frac{1}{2}\sum^d_{m=1}\frac{(x_m-x_m^{'})^2}{\ell^2_m}),
\end{equation}
which includes an individual correlated  {\color{red}lengthscale} $\ell_m$ for each input, with $d$ being the dimension of $\mathcal{D}$ and $\sigma_f$ being the signal variance.

Therefore, the superparameter set of GPR is $\boldsymbol{\mu}=\{\beta_1,\cdots,\beta_n,\ell_1,,\cdots,\ell_d,\sigma_f,\sigma_y\}$, which makes significant difference on the predictive performance. Based on the Bayesian maximum likelihood theory, we can estimate the optimal hyperparameters $\boldsymbol{\mu}_{opt}$ via solving the following problem
\begin{equation}\label{21}
\begin{split}
\boldsymbol{\mu}_{opt}&=\arg\max_{\boldsymbol{\mu}}(\log p(\textbf{y}_{\textbf{u},l}|\textbf{X},\boldsymbol{\mu}))\\
&=\arg\max_{\boldsymbol{\mu}}\{-\frac{1}{2}(\textbf{y}_{\textbf{u},l}-m(\textbf{X}))^T \textbf{K}_{\textbf{u},l}^{-1}(\boldsymbol{\mu})(\textbf{y}_{\textbf{u},l}-m(\textbf{X}))-\frac{1}{2}\log|\textbf{K}_{\textbf{u},l}(\boldsymbol{\mu})|-\frac{n}{2}\log(2\pi)\},
\end{split}
\end{equation}
where $p(\textbf{y}_{\textbf{u},l}|\textbf{X},\boldsymbol{\mu})$ is the conditional density function of $\textbf{y}_{\textbf{u},l}$ given $\textbf{X}$.
The procedure of the GPR is described in Algorithm 2.
\begin{algorithm}
\renewcommand{\algorithmicrequire}{\textbf{Input:}}
\renewcommand{\algorithmicensure}{\textbf{Output:}}
\caption{GPR}
\begin{algorithmic}[1]
\REQUIRE ~~\\
A training set of $n$ observations $D_{\textbf{u},l}, \textbf{u}=\mathbf{E}~\text{or}~\mathbf{H}$;\\
A chosen mean function $m(\cdot)$;\\
A kernel function $\kappa(\cdot,\cdot)$;\\
Test input $\textbf{x}^{\ast}$;
\ENSURE ~~\\
Test outputs result $\textbf{y}_{\textbf{u},l}^{\ast}|\textbf{x}^{\ast},\textbf{X},\textbf{y}_{\textbf{u},l}, \textbf{u}=\mathbf{E}~\text{or}~\mathbf{H}$;
\STATE Compute the optimal hyperparameters $\boldsymbol{\mu}_{opt}$ by maximizing the likelihood.
\STATE Compute $\textbf{K}_{\textbf{u},l}=\textbf{K}_{\textbf{u},l}(\boldsymbol{\mu}_{opt})$, $m(\cdot)=m(\cdot)(\boldsymbol{\mu}_{opt})$,\ $\textbf{K}_{\textbf{u},l}^{\ast\ast}=\textbf{K}_{\textbf{u},l}^{\ast\ast}(\boldsymbol{\mu}_{opt})$ and $\textbf{K}_{\textbf{u},l}^{\ast}=\textbf{K}_{\textbf{u},l}^{\ast}(\boldsymbol{\mu}_{opt})$;
\STATE Obtain the conditioning mean value $m^{\ast}(\textbf{x}^{\ast})=m(\textbf{x}^{\ast})+\kappa(\textbf{x}^{\ast},\textbf{X})\textbf{K}_{\textbf{u},l}^{-1}(\textbf{y}_{\textbf{u},l}-m(\textbf{X}))$
and covariance value $C^{\ast}(\textbf{x}^{\ast},\textbf{x}^{\ast})=\kappa(\textbf{x}^{\ast},\textbf{x}^{\ast})-\kappa(\textbf{x}^{\ast},\textbf{X})\textbf{K}_{\textbf{u},l}^{-1}\kappa(\textbf{X},\textbf{x}^{\ast})$;
\STATE Define $\textbf{y}_{\textbf{u},l}^{\ast}|\textbf{x}^{\ast},\textbf{X},\textbf{y}_{\textbf{u},l}\sim \mathcal{N}(m^{\ast}(\textbf{x}^{\ast}),C^{\ast}(\textbf{x}^{\ast},\textbf{x}^{\ast}))$.
\end{algorithmic}
\end{algorithm}
\subsection{Regression under singular value decomposition}
In general, the projection coefficients $\alpha_{\textbf{E}}(t,\theta)$ and $\alpha_{\textbf{H}}(t,\theta)$ vary more drastically with time than with parameter, which usually leads to a difficult global GPR. However, the single $1$D regression of time and $1$D (or multi-dimensional) regression of parameter are both easy to implement. Therefore, before constructing GPR models for $\alpha_{\textbf{E}}(t,\theta)$ and $\alpha_{\textbf{H}}(t,\theta)$, we firstly apply SVD to decompose the training data into independent time- and parameter-modes and to extract the corresponding principal components. Then, the GPR models for these decomposed time- and parameter-modes are constructed respectively, by which the global GPR $\widehat{\alpha}_{\textbf{E}}(t,\theta)$ and $\widehat{\alpha}_{\textbf{H}}(t,\theta)$ will be represented as linear combinations of several products of two Gaussian processes, one of time and the other of parameter.

Detailed procedure is presented here only for the case of the electric field. The same procedure can also be applied to the magnetic field.

For the $l$th coefficient $\alpha_{\textbf{E},l}=\psi_{\textbf{E},l}^T\underline{\textbf{E}}_h(t,\theta)$, $l=1,\cdots,d_{\textbf{E}}$, the training data can be written in a matrix as
\begin{equation}\label{22}
\textbf{P}_{\textbf{E},l}=[\alpha_{\textbf{E},l}(t^{n_i},\theta^{j})]_{ij}, 1\leq i\leq N_t^{tr}, 1\leq j\leq N_\theta^{tr}.
\end{equation}
We resort to the SVD to decompose $\textbf{P}_{\textbf{E},l}$ into several time- and parameter-modes
\begin{equation}\label{23}
\textbf{P}_{\textbf{E},l}\approx\widetilde{\textbf{P}}_{\textbf{E},l}=\sum_{k=1}^{Q_{\textbf{E},l}}\zeta_{\textbf{E},k}^l\boldsymbol{\xi}_{\textbf{E},k}^l\
(\boldsymbol{\phi}_{\textbf{E},k}^l)^T, 1\leq l \leq d_\textbf{E},
\end{equation}
where $\boldsymbol{\xi}_{\textbf{E},k}^l$ and $\boldsymbol{\phi}_{\textbf{E},k}^l$ are the $k$th discrete time- and parameter-modes for the $l$th projection coefficient, respectively, $\zeta_{\textbf{E},k}^l$ is the $k$th singular value, and $Q_{\textbf{E},l}$ is the truncation rank corresponding to the error tolerance $\delta_{\textbf{E},l}$, i.e., $Q_{\textbf{E},l}=\arg\max{\{ \mathcal{G}}(Q_{\textbf{E},l}):{\mathcal{G}}(Q_{\textbf{E},l})\geq1-\delta_{\textbf{E},l}\}$ with $\mathcal{G}(Q_{\textbf{E},l})=\sum_{k=1}^{Q_{\textbf{E},l}}(\zeta_{\textbf{E},k}^l)^2/\sum_{k=1}^{R_{{\textbf{E},l}}}(\zeta_{\textbf{E},k}^l)^2$ and $R_{{\textbf{E},l}}$ being the rank of $\textbf{P}_{\textbf{E},l}$. 

With the discrete modes database, GPR models can be trained to approximate the continuous modes as
\begin{equation}\label{24}
\begin{split}
&t \longmapsto \widehat{\xi}_{\textbf{E},k}^l(t),\text{trained from} \{(t^{n_i},(\boldsymbol{\xi}_{\textbf{E},k}^l)_i),i=1,2,\cdots N_t^{tr}\},\\
&\theta \longmapsto \widehat{\phi}_{\textbf{E},k}^l(\theta),\text{trained from} \{(\theta^{j},(\boldsymbol{\phi}_{\textbf{E},k}^l)_j),j=1,2,\cdots N_\theta^{tr}\}.
\end{split}
\end{equation}
Hence, we have
\begin{equation}\label{25}
(\textbf{P}_{\textbf{E},l})_{ij}=\alpha_{\textbf{E},l}(t^{n_i},\theta^{j})\approx\sum_{k=1}^{Q_l}\zeta_{\textbf{E},k}^l\
\widehat{\xi}_{\textbf{E},k}^l(t^{t^{n_i}})\widehat{\phi}_{\textbf{E},k}^l(\theta^{j}),
\end{equation}
with $1\leq i\leq N_t^{tr}$, $1\leq j\leq N_\theta^{tr}.$

The continuous regression function $\widehat{\alpha}_{\textbf{E},l}(t,\theta)$ for the $l$th projection coefficient $\alpha_{\textbf{E},l}(t,\theta)$ can be recovered as
\begin{equation}\label{26}
\alpha_{\textbf{E},l}(t,\theta)\approx\widehat{\alpha}_{\textbf{E},l}(t,\theta)=\sum_{k=1}^{Q_l}\zeta_{\textbf{E},k}^l\
\widehat{\xi}_{\textbf{E},k}^l(t)\widehat{\phi}_{\textbf{E},k}^l(\theta),(t,\theta)\in\mathcal{T}\times\mathcal{P}.
\end{equation}

The model order reduction process based on POD-GPR proposed in this paper is shown in Algorithm 3.

\begin{algorithm}
\caption{POD-GPR Reduced order method for electromagnetic problem}
\label{alg:A}
\begin{algorithmic}[1]
\STATE \textbf{Offline stage:}
\STATE Generate full-order solutions $\underline{\textbf{u}}_h(t_{n_i},\theta_j)$, $i=1,\cdots,N_l$, $j=1,\cdots,N_\theta$ and construct snapshot matrices $\mathbb{S}_{\textbf{u},j}$,  $\textbf{u}=\mathbf{E}~\text{or}~\mathbf{H}$;
\STATE Extract the reduced basis $\Psi_\textbf{u}=[\psi_{\textbf{u},1},\cdots,\psi_{\textbf{u},d_{\textbf{u}}}]$ ($\textbf{u}=\mathbf{E}~\text{or}~\mathbf{H}$) through Algorithm 1;
\STATE Generate full-order solutions $\underline{\textbf{u}}_h(t^{n_i},\theta^j)$, $i=1,2,\cdots,N_t^{tr}, j=1,2,\cdots,N_\theta^{tr}$ and form the training matrices $\textbf{P}_{\textbf{u},l}=[\psi_{\textbf{u},l}^T\underline{\textbf{E}}_h(t^{n_i},\theta^j)]_{ij}$, $l=1,2,\cdots,d_{\textbf{u}}$, $\textbf{u}=\mathbf{E}~\text{or}~\mathbf{H}$;
\STATE Perform SVD with $\delta_{\textbf{u},l}$ on $\textbf{P}_{\textbf{u},l}$, $\textbf{P}_{\textbf{u},l}=\sum_{k=1}^{Q_{\textbf{u},l}}\zeta_{\textbf{u},k}^l\boldsymbol{\xi}_{\textbf{u},k}^l\
(\boldsymbol{\phi}_{\textbf{u},k}^l)^T$, $l=1,2,\cdots,d_{\textbf{u}}$, $\textbf{u}=\mathbf{E}~\text{or}~\mathbf{H}$;
\STATE Construct GPR models $\widehat{\xi}_{\textbf{u},k}^l(t)$ and $\widehat{\phi}_{\textbf{u},k}^l(\theta)$, $k=1,2,\cdots,Q_{\textbf{u},l}$, $l=1,2,\cdots,d_\textbf{u}$, $\textbf{u}=\mathbf{E}~\text{or}~\mathbf{H}$;
\STATE Recover every projection coefficient as $\widehat{\alpha}_{\textbf{u},l}(t,\theta)=\sum_{k=1}^{Q_{\textbf{u},l}}\zeta_{\textbf{u},k}^l\
\widehat{\xi}_{\textbf{u},k}^l(t)\widehat{\phi}_{\textbf{u},k}^l(\theta)$, $\textbf{u}=\mathbf{E}~\text{or}~\mathbf{H}$;
\STATE \textbf{Online stage:}
\STATE Recover output $\widehat{\alpha}_{\textbf{u}}(t^\ast,\theta^\ast)$ for a new parameter value $(t^\ast,\theta^\ast)$;
\STATE Evaluate the reduced-order solution $\underline{\textbf{u}}_{h,reg}^d(t^\ast,\theta^\ast)=\sum_{l=1}^{d_{\textbf{u}}}\widehat{\alpha}_{\textbf{u},l}(t^{\ast},\theta^{\ast})\psi_{\textbf{u},l}=
\Psi_{\textbf{u}}\widehat{\alpha}_{\textbf{u}}(t^\ast,\theta^\ast)$.
\end{algorithmic}
\end{algorithm}

\begin{remark} After assuming that the full-order solutions in step 2 of Algorithm 3 are used both as snapshots and as training data, i.e., $\{\underline{\textbf{u}}_h(t_{n_i},\theta_j): i=1,\cdots,N_l, j=1,\cdots,N_\theta\}=\{\underline{\textbf{u}}_h(t^{n_i},\theta^j): i=1,2,\cdots,N_t^{tr}, j=1,2,\cdots,N_\theta^{tr}\}$, 
the total error between the projection solutions $\underline{\textbf{u}}_h^d(t_{n_i},\theta_j)=\Psi_\textbf{u}\Psi_\textbf{u}^T\underline{\textbf{u}}_h(t_{n_i},\theta_j)$ and the POD-GPR reduced-order solutions $\underline{\textbf{u}}_{h,reg}^d(t_{n_i},\theta_j)$ ($i=1,\cdots,N_l$, $j=1,\cdots,N_\theta$) can be estimated as
\begin{equation}\label{27}
\begin{split}
\sum_{j=1}^{N_\theta} \sum_{i=1}^{N_l}\|\underline{\textbf{u}}_h^d(t_{n_i},\theta_j)-\underline{\textbf{u}}_{h,reg}^d(t_{n_i},\theta_j)\|^2_{\mathbb{R}^{N_h}}&=\sum_{j=1}^{N_\theta} \sum_{i=1}^{N_l}\|\Psi_\textbf{u}\Psi_\textbf{u}^T\underline{\textbf{u}}_h(t_{n_i},\theta_j)-\underline{\textbf{u}}_{h,reg}^d(t_{n_i},\theta_j)\|^2_{\mathbb{R}^{N_h}}\\
&=\sum_{j=1}^{N_\theta} \sum_{i=1}^{N_l}\|\sum_{l=1}^{d_\textbf{u}}(\textbf{P}_{\textbf{u},l})_{ij}\psi_{\textbf{u},l}-\sum_{l=1}^{d_\textbf{u}}(\widetilde{\textbf{P}}_{\textbf{u},l})_{ij}\psi_{\textbf{u},l}\|^2_{\mathbb{R}^{N_h}}\\
&=\sum_{l=1}^{d_\textbf{u}}\|\textbf{P}_{\textbf{u},l}-\widetilde{\textbf{P}}_{\textbf{u},l}\|^2_F=\sum_{l=1}^{d_\textbf{u}}\frac{ \|\textbf{P}_{\textbf{u},l}-\widetilde{\textbf{P}}_{\textbf{u},l}\|^2_F}{\|\textbf{P}_{\textbf{u},l}\|^2_F} \|\textbf{P}_{\textbf{u},l}\|^2_F\\
&\leqslant \sum_{l=1}^{d_\textbf{u}}(\delta_{\textbf{u},l}\sum_{k=1}^{R_{{\textbf{E},l}}}(\zeta_{\textbf{u},k}^l)^2), \textbf{u}=\mathbf{E}~\text{or}~\mathbf{H},
\end{split}
\end{equation}
by which the total recovery error of the POD-GPR can be expressed as
\begin{equation}\label{28}
\begin{split}
&   \sum_{j=1}^{N_\theta} \sum_{i=1}^{N_l}\|\underline{\textbf{u}}_h(t_{n_i},\theta_j)-\underline{\textbf{u}}_{h,reg}^d(t_{n_i},\theta_j)\|_{\mathbb{R}^{N_h}}\\
&\leqslant\sum_{j=1}^{N_\theta} \sum_{i=1}^{N_l}(\|\textbf{u}_h(t_{n_i},\theta_j)-\Psi_\textbf{u}\Psi_\textbf{u}^T\underline{\textbf{u}}_h(t_{n_i},\theta_j)\|_{\mathbb{R}^{N_h}}+\|\Psi_\textbf{u}\Psi_\textbf{u}^T\underline{\textbf{u}}_h(t_{n_i},\theta_j)-\underline{\textbf{u}}_{h,reg}^d(t_{n_i},\theta_j)\|_{\mathbb{R}^{N_h}})\\
&\leqslant \sqrt{\epsilon_{\textbf{u},t}}\mathcal{L}_1+\sqrt{\epsilon_{\textbf{u},\theta}}\mathcal{L}_2+({\color{red}N_s}\sum_{l=1}^{d_\textbf{u}}(\delta_{\textbf{u},l}\sum_{k=1}^{R_{{\textbf{E},l}}}(\zeta_{\textbf{u},k}^l)^2))^\frac{1}{2},  \textbf{u}=\mathbf{E}~\text{or}~\mathbf{H},
\end{split}
\end{equation}
which contains the truncation errors arising from both the POD and the GPR models. This shows that the value of $\epsilon_{\textbf{u},t}, \epsilon_{\textbf{u},\theta}$ and $\delta_{\textbf{u},l}$ ($l=1,2,\cdots d_\textbf{u}, \textbf{u}=\mathbf{E}~\text{or}~\mathbf{H}$) have a great influence on the accuracy of the whole POD-GPR algorithm and can act as a benchmark for the error control.

\end{remark} 
\section{Numerical results}
\label{sec:tests}
In this section, numerical results for two electromagnetic scattering problems are displayed to validate the effectiveness and the accuracy of the proposed method. We consider the solution of the 2-D time-domain Maxwell's equations in the case of transverse magnetic (TM) waves
 \begin{equation}
\begin{cases}
\mu_r\displaystyle\frac{\partial H_x}{\partial t}+\displaystyle\frac{\partial E_z}{\partial y} = 0,\\[2ex]
\mu_r\displaystyle\frac{\partial H_y}{\partial t}-\displaystyle\frac{\partial E_z}{\partial x} = 0,\\[2ex]
\varepsilon_r\displaystyle\frac{\partial E_z}{\partial t}-\displaystyle\frac{\partial H_y}{\partial x}+\displaystyle\frac{\partial H_x}{\partial y} = 0.
\end{cases}
\end{equation}
The excitation in all considered scattering scenarios is an incident plane wave defined as
 \begin{equation}
\begin{cases}
H_x^{inc}(x,y,t)=0,\\[2ex]
H_y^{inc}(x,y,t)=-\cos(\omega t-k x),\\[2ex]
E_z^{inc}(x,y,t)=\cos(\omega t-k x),
\end{cases}
\end{equation}
where $\omega=2\pi f$ is the angular frequency with the wave frequency $f=30 GHz$, and $k=\frac{\omega}{c}$ is the wave number, $c$ is the wave speed in vacuum.

The  relative $L^2$ error between the POD-GPR reduced solution and the DGTD high-fidelity solution is utilized as the metric to evaluate the accuracy of the results
 \begin{equation}
\varepsilon_{POD-GPR}(t,\theta)=\frac{\|\underline{\textbf{u}}_{h}(t,\theta)-\underline{\textbf{u}}_{h,reg}(t,\theta)\|_{L^2}}{\|\underline{\textbf{u}}_{h}(t,\theta)\|_{L^2}}=\frac{\|\underline{\textbf{u}}_{h}(t,\theta)-\Psi_{\textbf{u}}\widehat{\alpha}_{\textbf{u}}(t,\theta)\|_{L^2}}{\|\underline{\textbf{u}}_{h}(t,\theta)\|_{L^2}}, \textbf{u}=\mathbf{E}~\text{or}~\mathbf{H},
\end{equation}
which will be compared with the relative projection $L^2$ error committed by POD
\begin{equation}
\varepsilon_{Projection}(t,\theta)=\frac{\|\underline{\textbf{u}}_{h}(t,\theta)-\Psi_{\textbf{u}} \Psi_{\textbf{u}} ^T \textbf{u}_{h}(t,\theta)\|_{L^2}}{\|\underline{\textbf{u}}_{h}(t,\theta)\|_{L^2}}, \textbf{u}=\mathbf{E}~\text{or}~\mathbf{H}.
\end{equation}

Simulations are run on a Macbook equipped with an Intel Core i5 1.8 GHz CPU and 8 GB memory, and GPR models are constructed by the MATLAB function $\verb"fitrgp"$.
\subsection{Scattering of a plane wave by a dielectric cylinder}
We first investigate the electromagnetic scattering of a plane wave by a dielectric cylinder.
The computation domain is artificially truncated by the square $\Omega = [-2.6m,2.6m]\times[-2.6m,2.6m]$, on which the first order Silver-M\"{u}ller ABC boundary condition is imposed. The cylinder is located at the origin and its radius is $0.6m$. Our interest relative permittivity of the cylinder is $\varepsilon_r\in[1,5]$ (i.e., $\mathcal{P}=[1,5]$) and we set $\mu_r=1$ (i.e., nonmagnetic material). The medium exterior to the dielectric cylinder is assumed to be vacuum, i.e. $\varepsilon_{r,1}=1$ and $\mu_{r,1}=1$.

The full-order simulations are performed on an unstructured triangular mesh with 1733 nodes and 3380 elements, in which 780 elements are located inside the cylinder, and the minimal and maximum mesh size is $3.92\times10^{-2}m$ and $3.359\times10^{-1}m$, leading to the number of DOFs of the full-order model $N_h=20280$. The total simulation time corresponds to 50 periods of the incident wave oscillation.

During the offline stage, we apply DGTD solver to obtain full-order solutions at $N_\theta=81$ equidistant parameter sampling points (i.e., $\theta\in\mathcal{P}_h=\{1,1.05,1.10,\cdots,4.95,5\}$), among which we consider $N_l=218$ equally collected time points in the last oscillation period (i.e., $t\in\mathcal{T}_{h}=\{49.0040,49.0084,49.0128,\cdots,49.9678\}$). All these transient vectors will be used both as snapshots and training data. With the two-step POD criteria of $\epsilon_{\textbf{E},t}=\epsilon_{\textbf{H},t}=1\times e^{-3}$ and $\epsilon_{\textbf{E},\theta}=\epsilon_{\textbf{E},\theta }=1\times e^{-4}$, $d_{E_z}=18$ POD bases are extracted for $E_z$, $d_{H_y}=16$ for $H_y$ and $d_{H_x}=159$ for $H_x$. As for the training of GPR models, the SVD truncation tolerance $\delta_{\textbf{u},l}$ ($l=1,2,\cdots,d_\textbf{u}$, $\textbf{u}=\mathbf{E}~\text{or}~\mathbf{H}$) are set in groups as is shown in Table \ref{tab:1}, which depends on the fact that lower-order projection coefficients contains more dominant information about the origin model. 

 \begin{table}[]
\caption{Scattering of a plane wave by a dielectric cylinder: The SVD truncation tolerance.}
\label{tab:1} 
\centering
\begin{tabular}{cc}
\hline
The projection coefficient item& $\delta_{\textbf{E},l}, \delta_{\textbf{H},l}$\\
\hline
$l\leq5$ & $1\times 10^{-4}$\\
$5<l\leq10$ & $5\times 10^{-4}$\\
$10<l\leq20$ & $1\times 10^{-3}$\\
$20<l\leq30$ & $2\times 10^{-3}$\\
$30<l\leq40$ & $3\times 10^{-3}$\\
$40<l\leq55$ & $4\times 10^{-3}$\\
$55<l$ & $5\times 10^{-3}$\\
\hline
\end{tabular}
\end{table}

\begin{figure}[htbp]
\centering
\includegraphics[height=0.395\textheight, width=0.495\textheight]{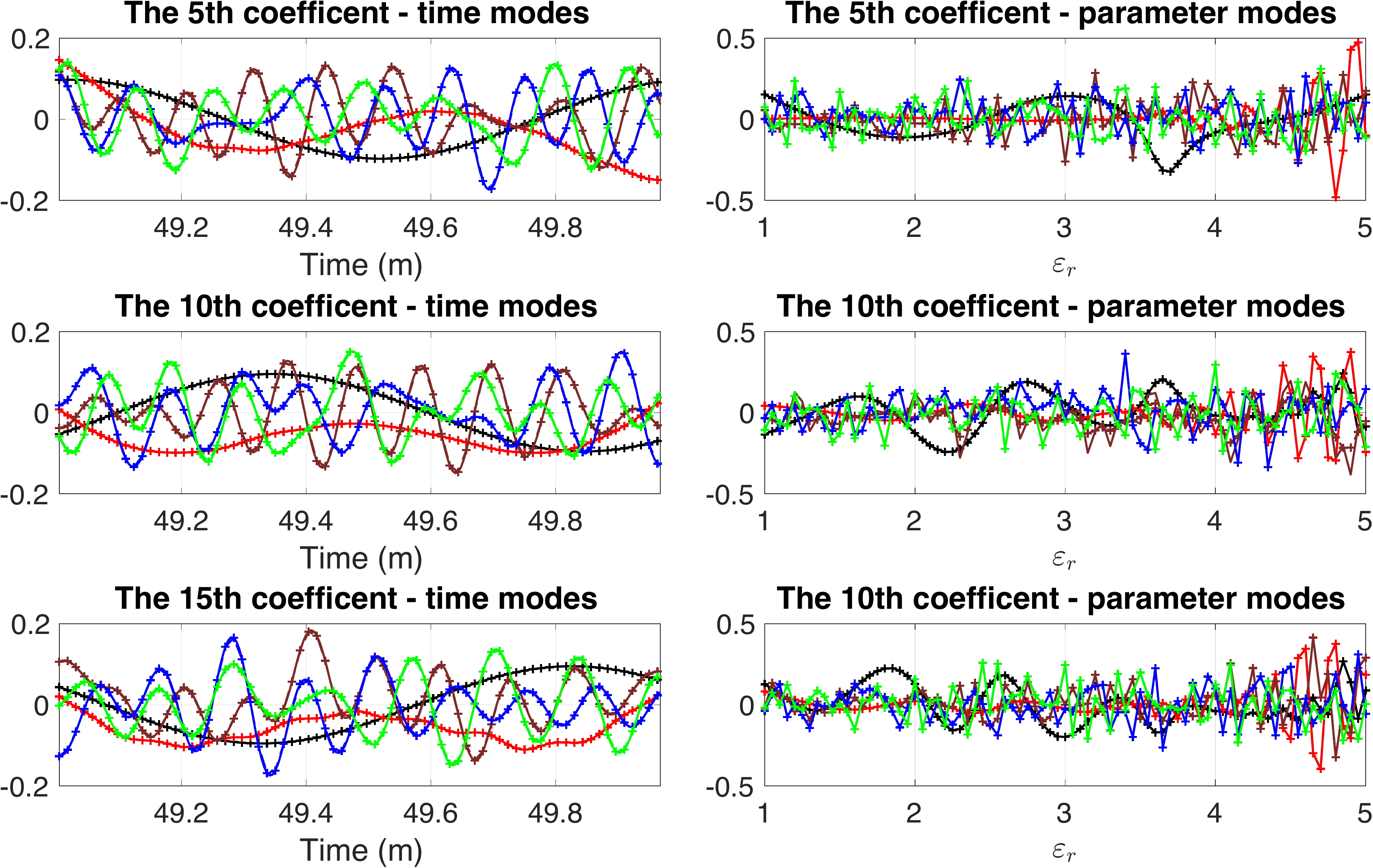}
\caption{Scattering of a plane wave by a dielectric cylinder: Time- and parameter-modes for the 5th, 10th and 15th projection coefficients for $E_z$: the 2nd modes-black, the 4th modes-red, the 6th modes-brown, the 8th modes-blue, the 10th modes-green. (For interpretation of the colors in the figure(s), the reader is referred to the web version of this article.) }
\label{fig:2} 
\end{figure}

\begin{figure}[htbp]
\centering
\includegraphics[height=0.395\textheight, width=0.495\textheight]{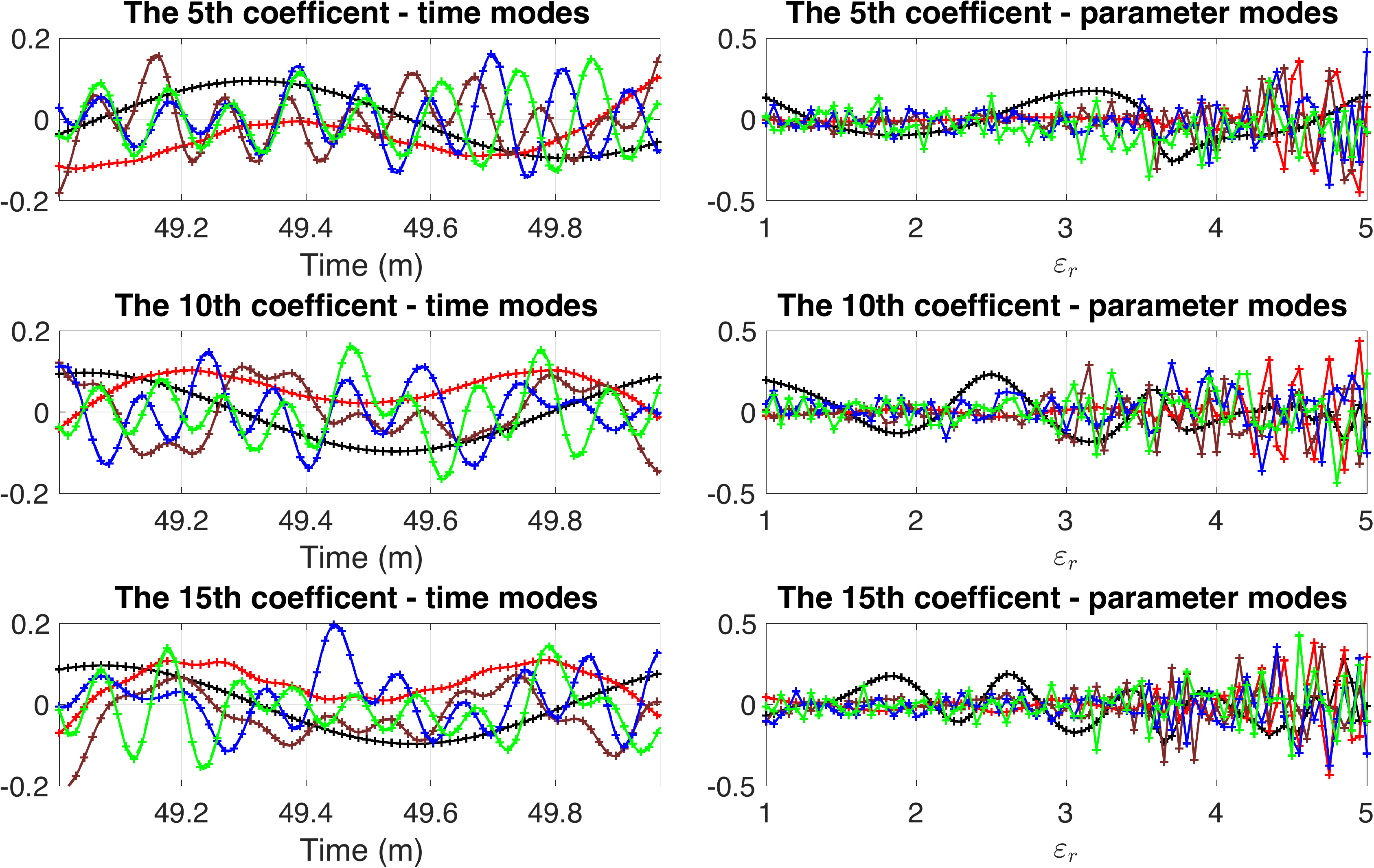}
\caption{Scattering of a plane wave by a dielectric cylinder: Time- and parameter-modes for the 5th, 10th and 15th projection coefficients for $H_y$: the 2nd modes-black, the 4th modes-red, the 6th modes-brown, the 8th modes-blue, the 10th modes-green.}
\label{fig:3} 
\end{figure}

\begin{figure}[htbp]
\centering
\includegraphics[height=0.395\textheight, width=0.495\textheight]{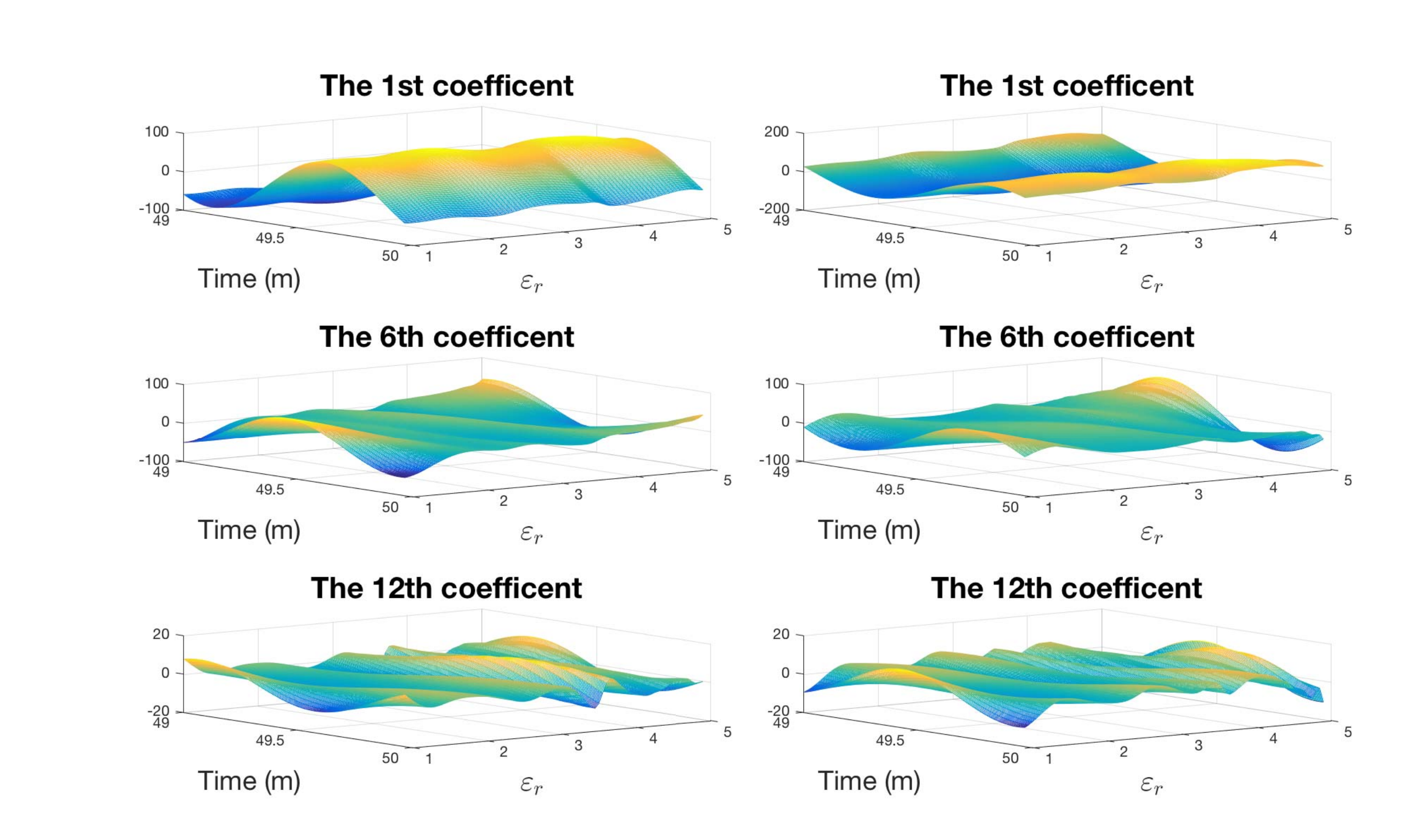}
\caption{Scattering of a plane wave by a dielectric cylinder: Global fitting results for the 1st, 6th and 12th projection coefficients for $E_z$ (left) and $H_y$ (right).}
\label{fig:4} 
\end{figure}

Based on SVD procedures on training matrices, the time- and parameter-modes are approximated via GPR models, some of which are shown in Fig.\ref{fig:2} for $E_z$ and  Fig.\ref{fig:3} for $H_x$. The dotted line represents the training point and the solid line represents the regression function. As we can see, for each projection coefficient, whether it is time or parameter regression, the lower order mode is smoother. Although for parameter regression, the shock intensity of high order modes is relatively large, but all of them have little influence on the performance of the whole model, which is attributed to the dominant role played by low order modes. This also reflects the robustness of the regression model. Therefore, the global GPR of all coefficients can be obtained and some of which is displayed in Fig.\ref{fig:4}.


To verify the performance of POD-GPR model built here, the reduced-order electromagnetic fields are recovered under some non-trained arbitrarily chosen parameters: $\varepsilon_r=1.215$, $\varepsilon_r=2.215$, $\varepsilon_r=3.215$ and $\varepsilon_r=4.215$, which are calculated from the POD basis, with their coefficients obtained as direct outputs from regression models. These arbitrary test instances are then compared with the corresponding DGTD high-fidelity solutions.

\begin{figure}
\centering
\subfigure[]
{
\includegraphics[width=2.2in]{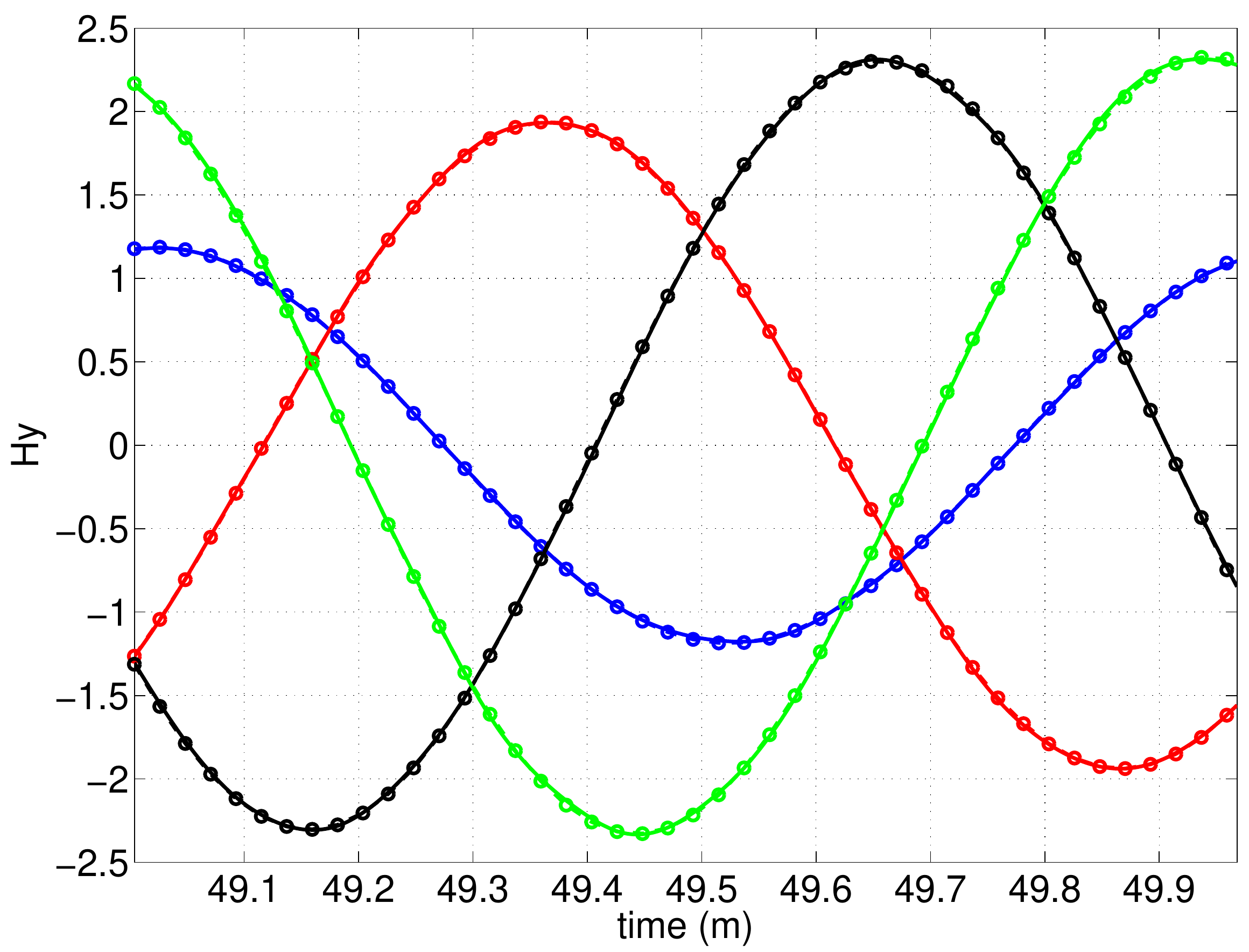}
}
\hspace{0.15cm}
\subfigure[]
{
\includegraphics[width=2.8in]{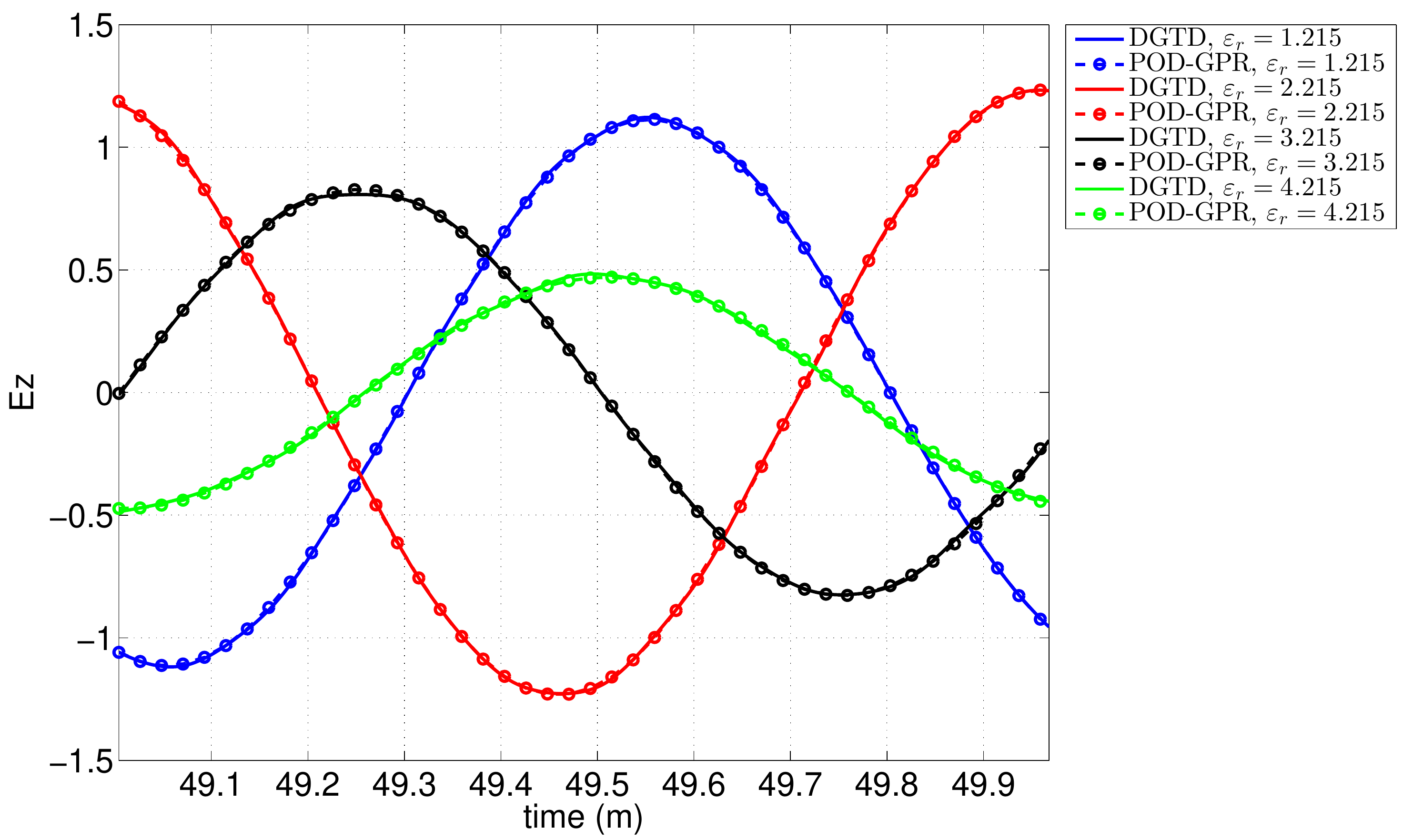}
}
\caption{Scattering of a plane wave by a dielectric cylinder: Comparison of the time evolution of the field (a) $H_y$ and (b) $E_z$ at a given point.}
\label{fig:5}
\end{figure}

\begin{figure}
\centering
\subfigure
{
\includegraphics[width=2.2in]{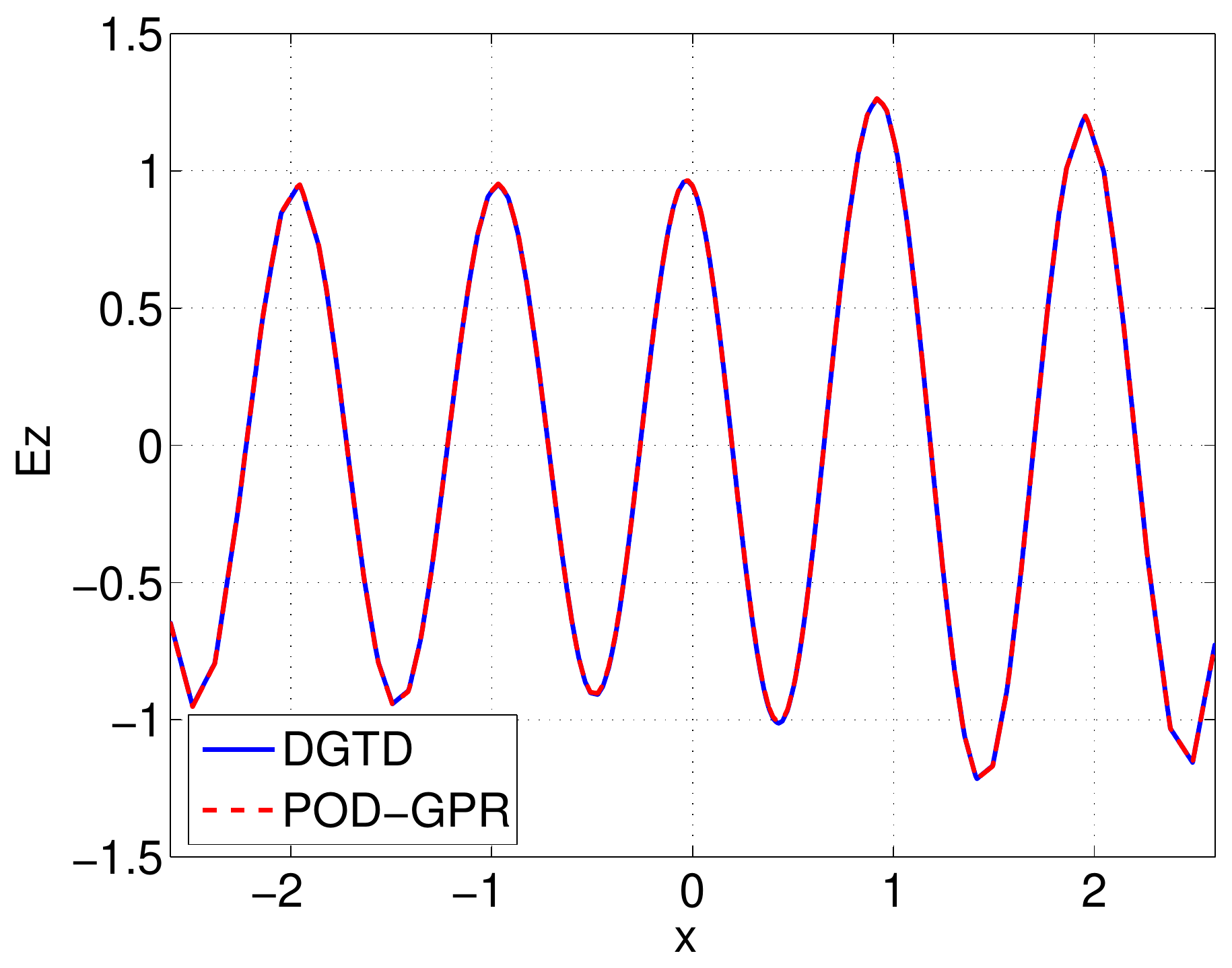}
}
\hspace{0.15cm}
\subfigure
{
\includegraphics[width=2.17in]{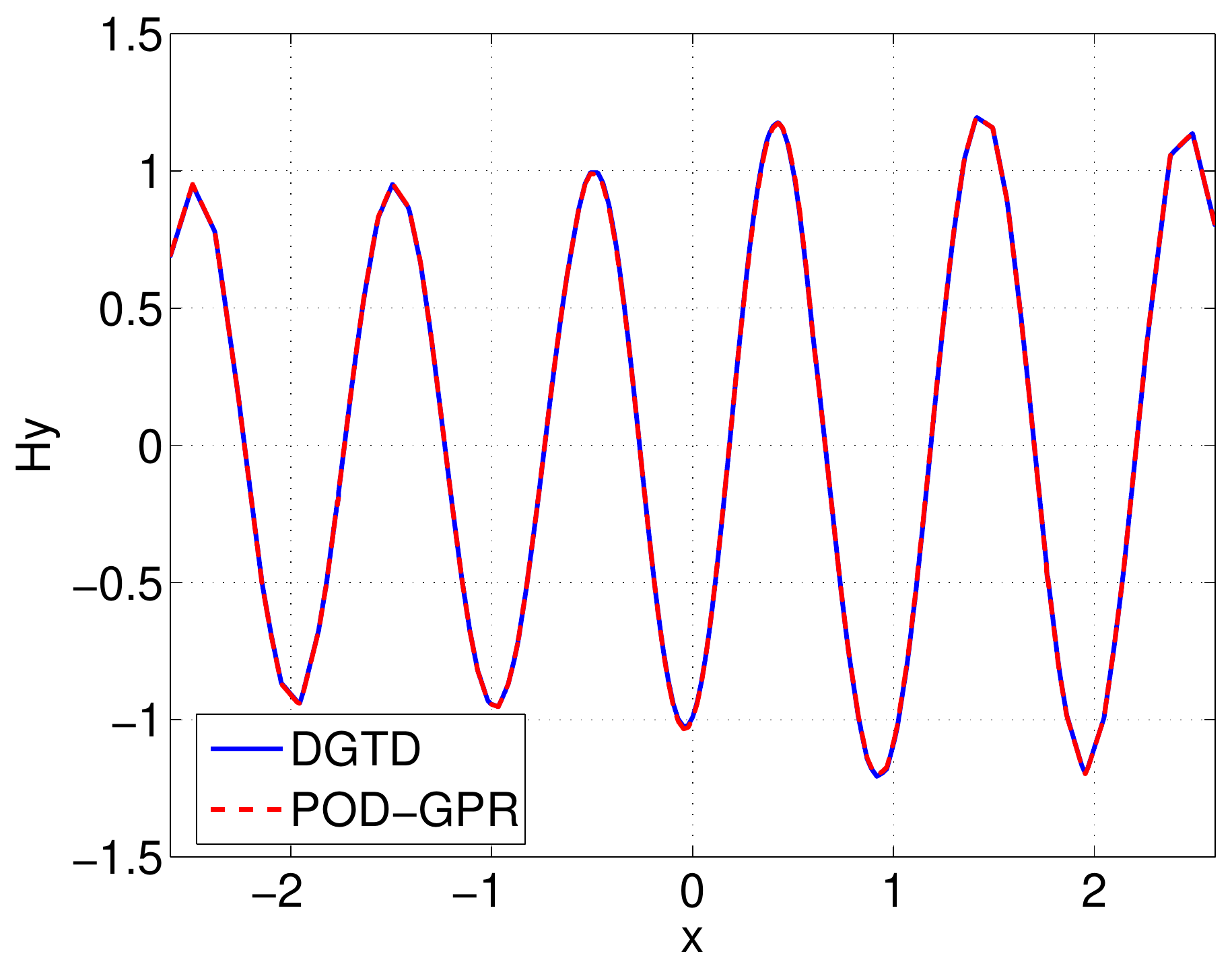}
}\\
\subfigure
{
\includegraphics[width=2.2in]{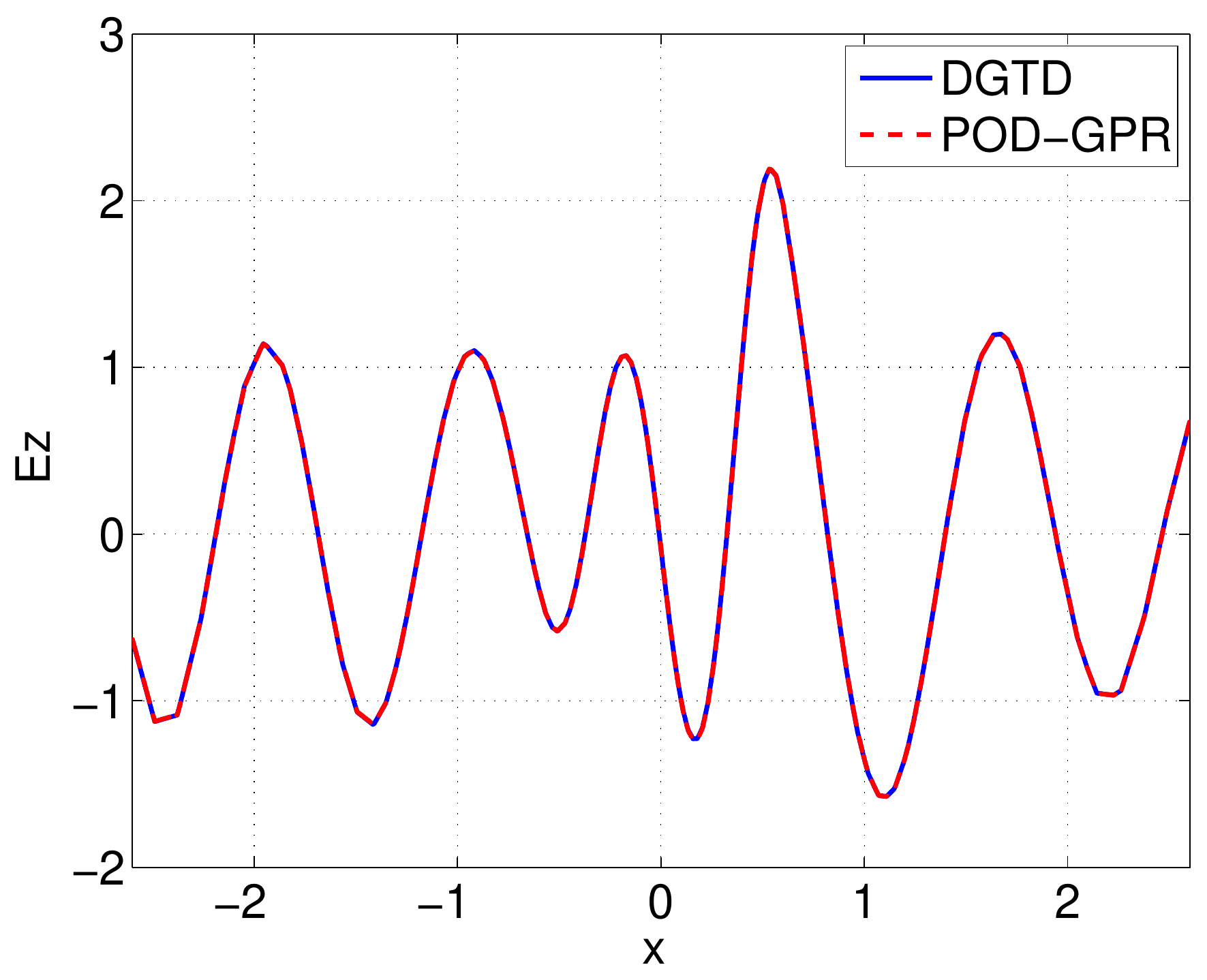}
}
\hspace{0.15cm}
\subfigure
{
\includegraphics[width=2.24in]{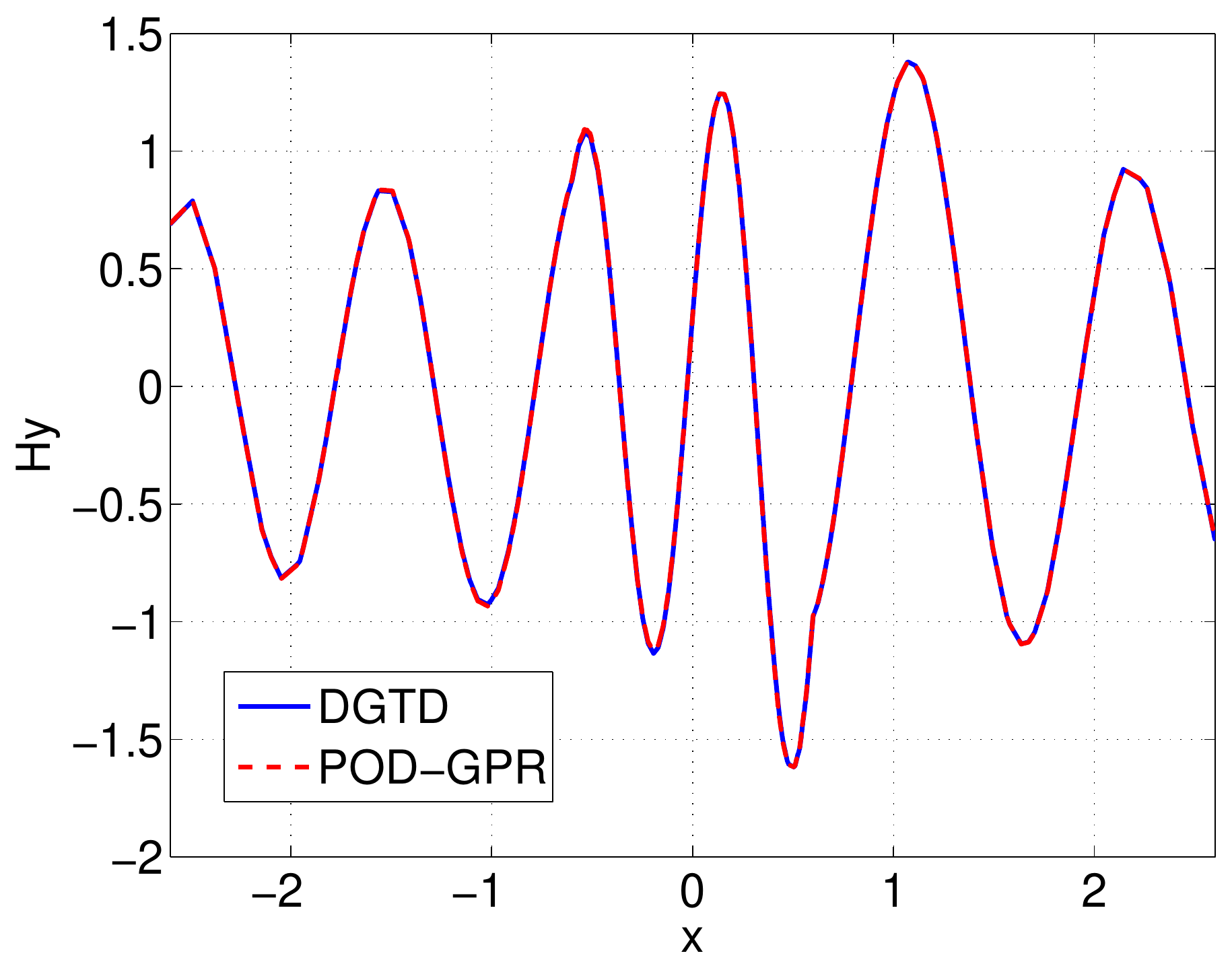}
}\\
\subfigure
{
\includegraphics[width=2.2in]{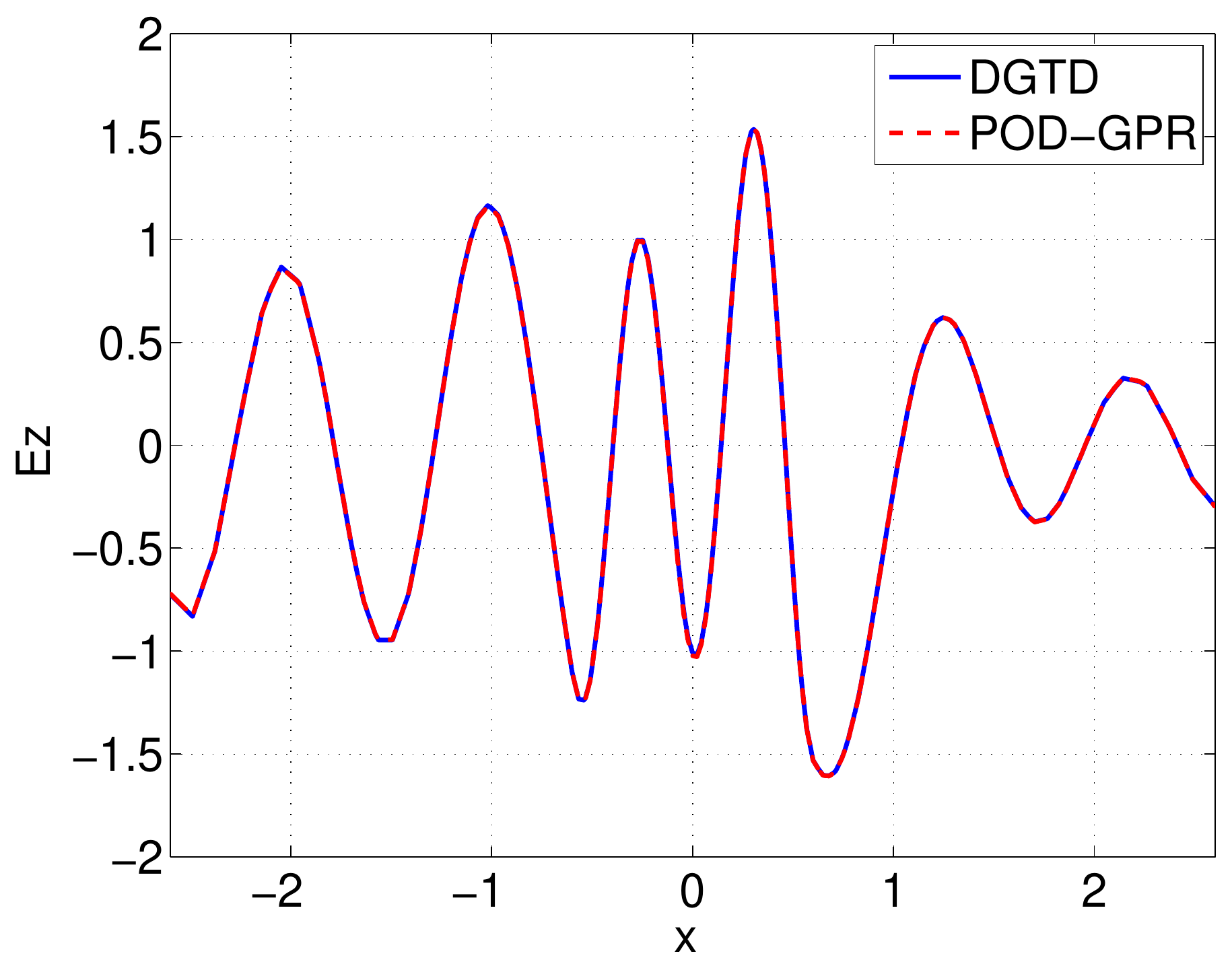}
}
\hspace{0.15cm}
\subfigure
{
\includegraphics[width=2.13in]{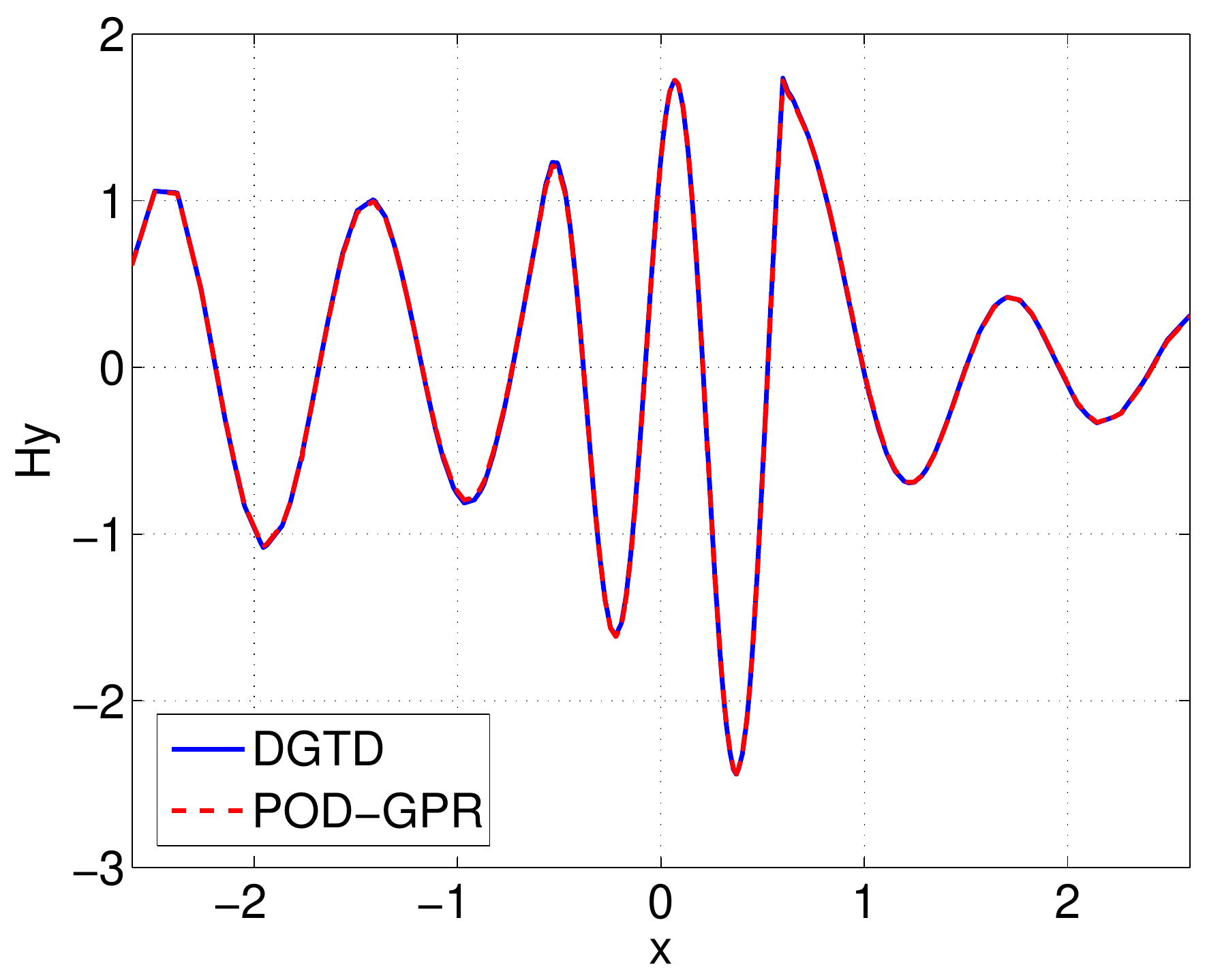}
}\\
\subfigure
{
\includegraphics[width=2.2in]{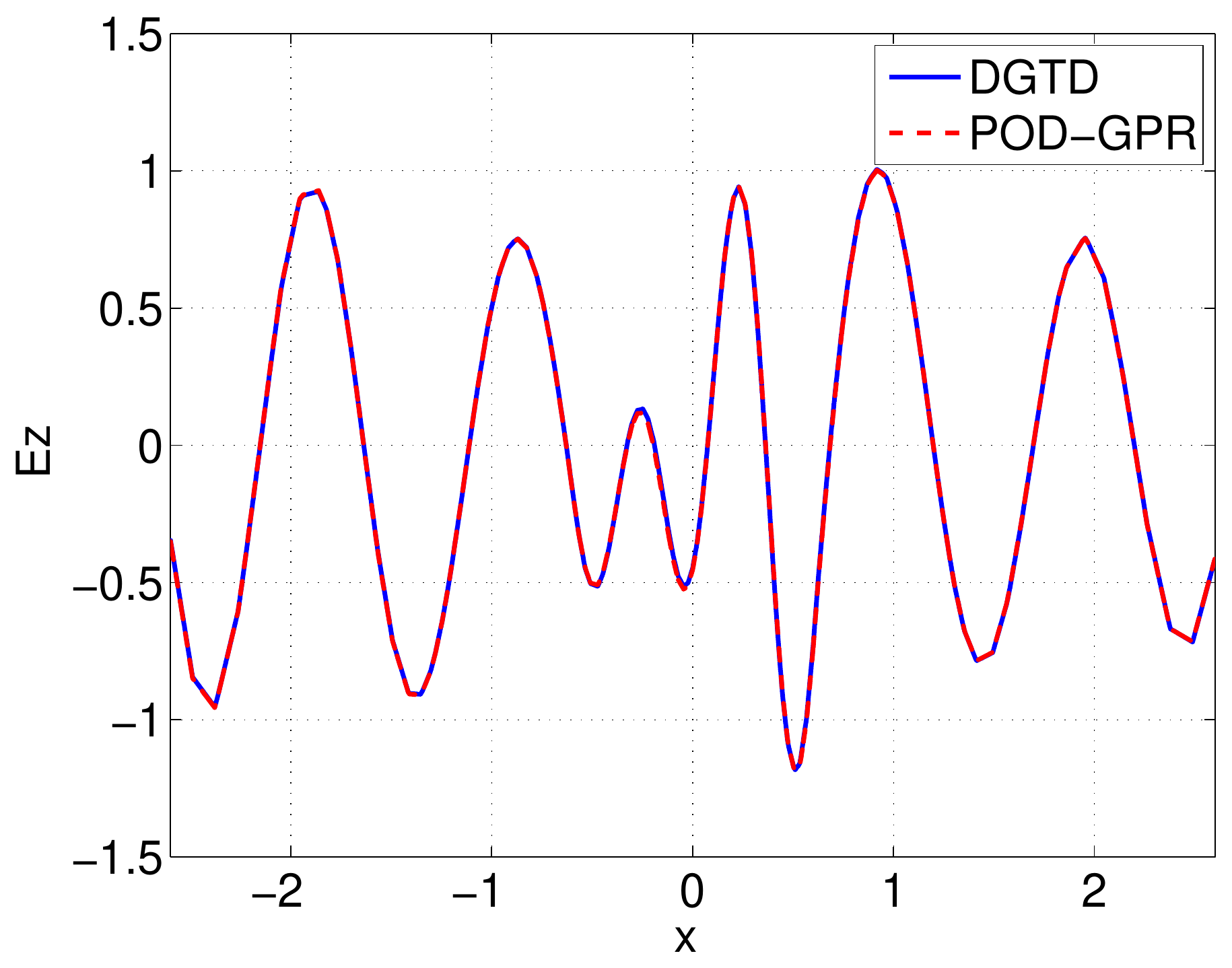}
}
\hspace{0.15cm}
\subfigure
{
\includegraphics[width=2.13in]{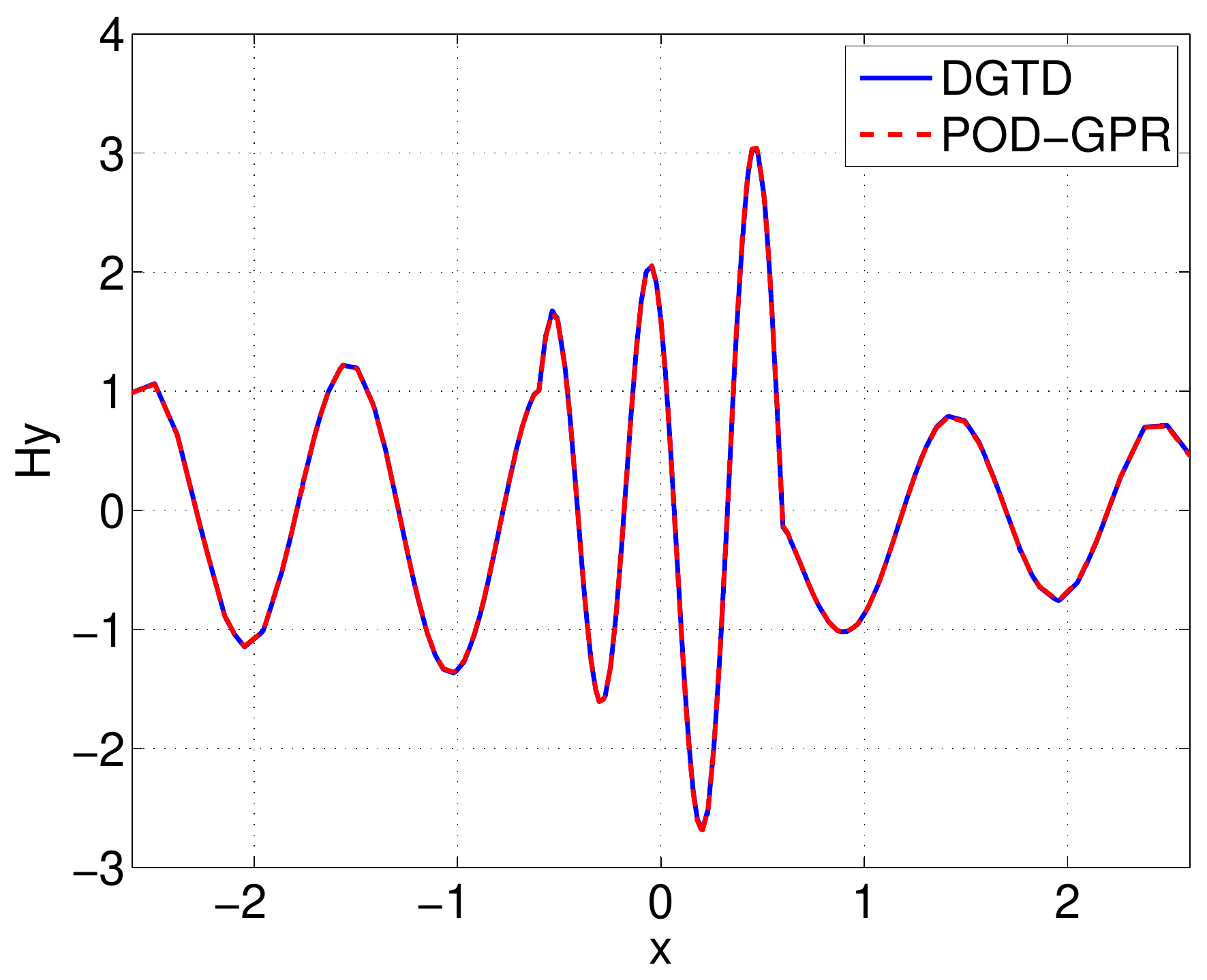}
}\\
\caption{Scattering of a plane wave by a dielectric cylinder: Comparison of the 1D x-wise distribution along $y=0$ of the real part of $E_z$ (left) and $H_y$ (right) of four test points: $\varepsilon_r=1.215$ (1st row), $\varepsilon_r=2.215$ (2nd row), $\varepsilon_r=3.215$ (3rd row) and $\varepsilon_r=4.215$ (4th row).}
\label{fig:6} 
\end{figure}

\begin{figure}
\centering
\subfigure
{
\includegraphics[width=2.2in]{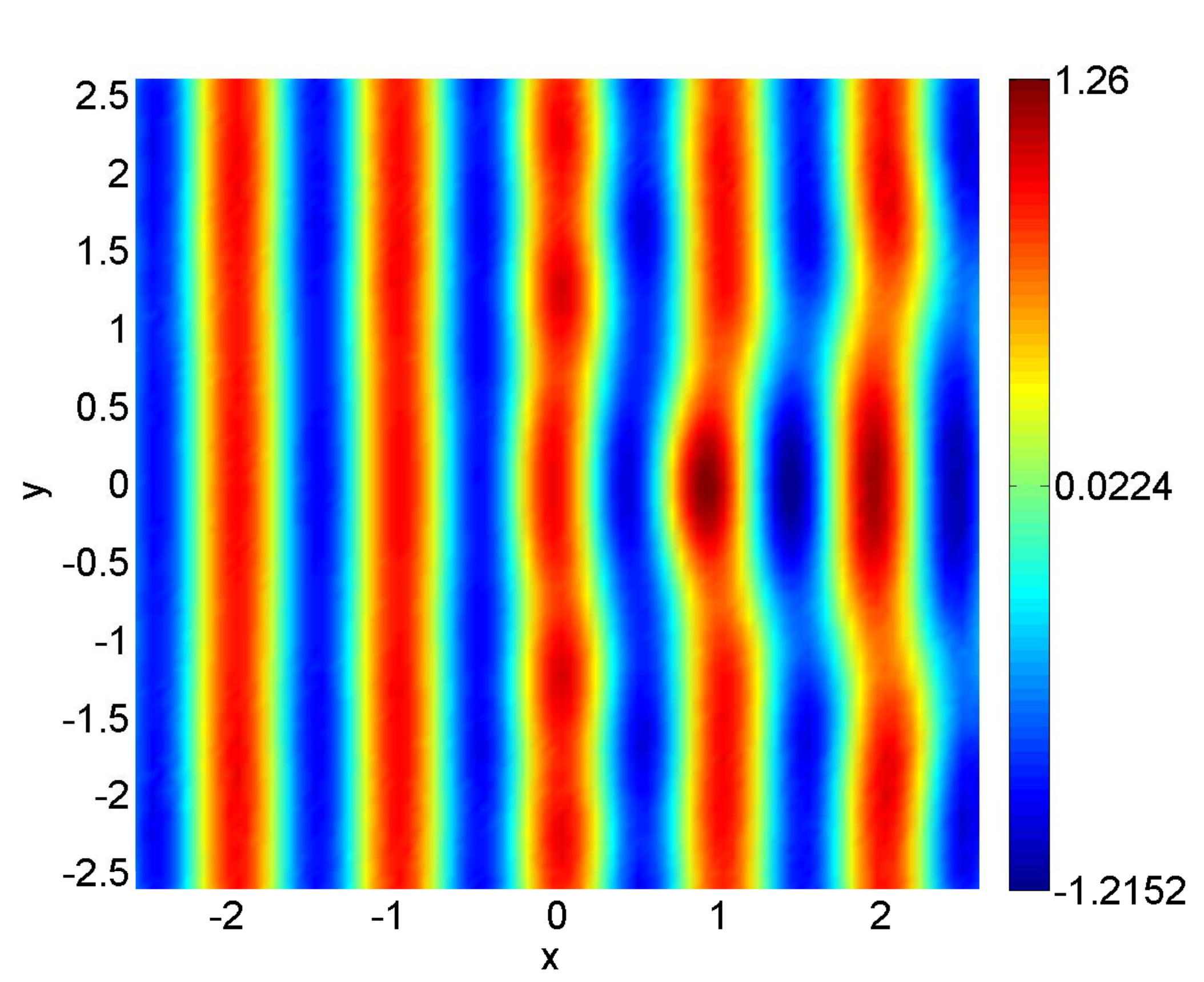}
}
\hspace{0.15cm}
\subfigure
{
\includegraphics[width=2.2in]{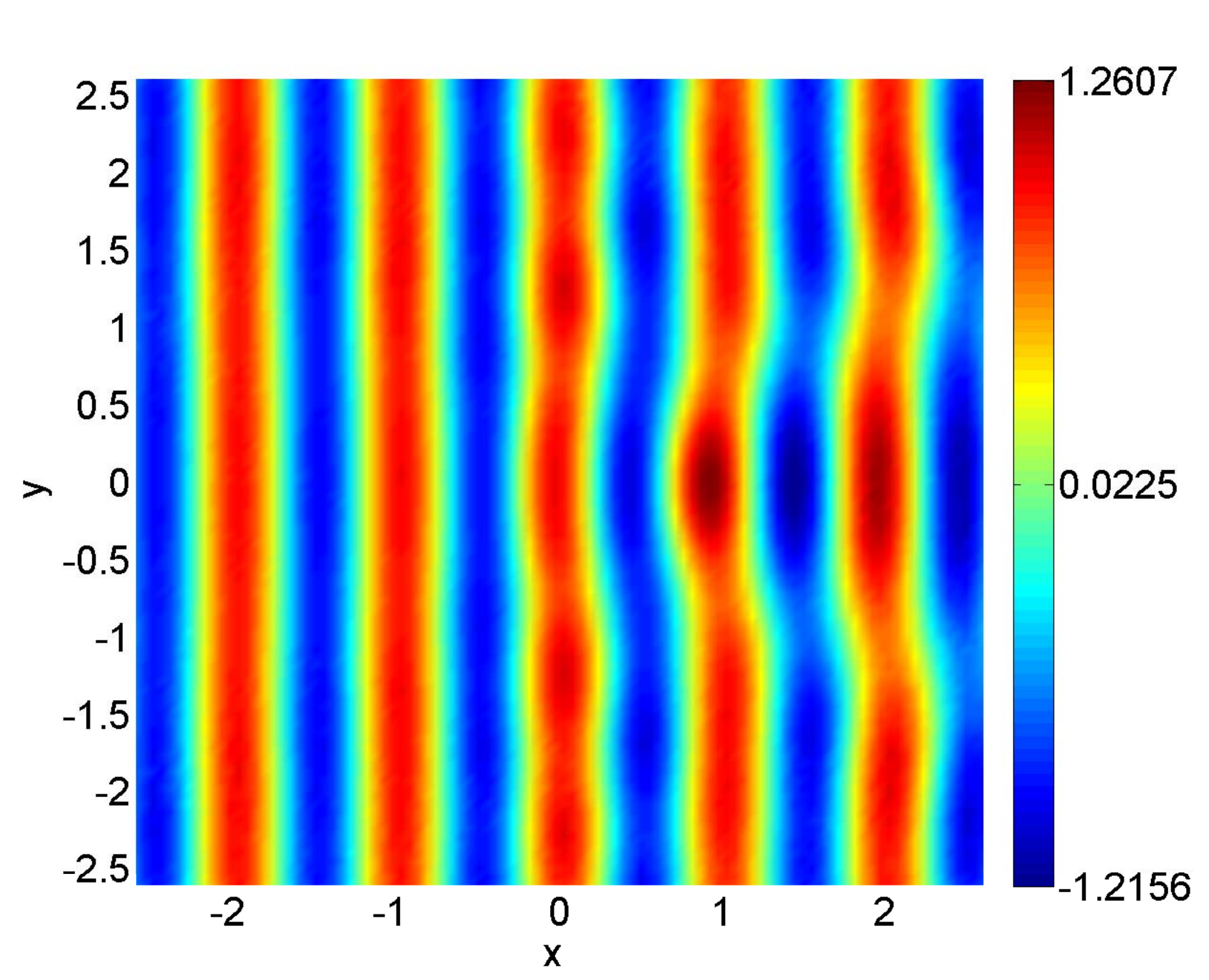}
}\\
\subfigure
{
\includegraphics[width=2.2in]{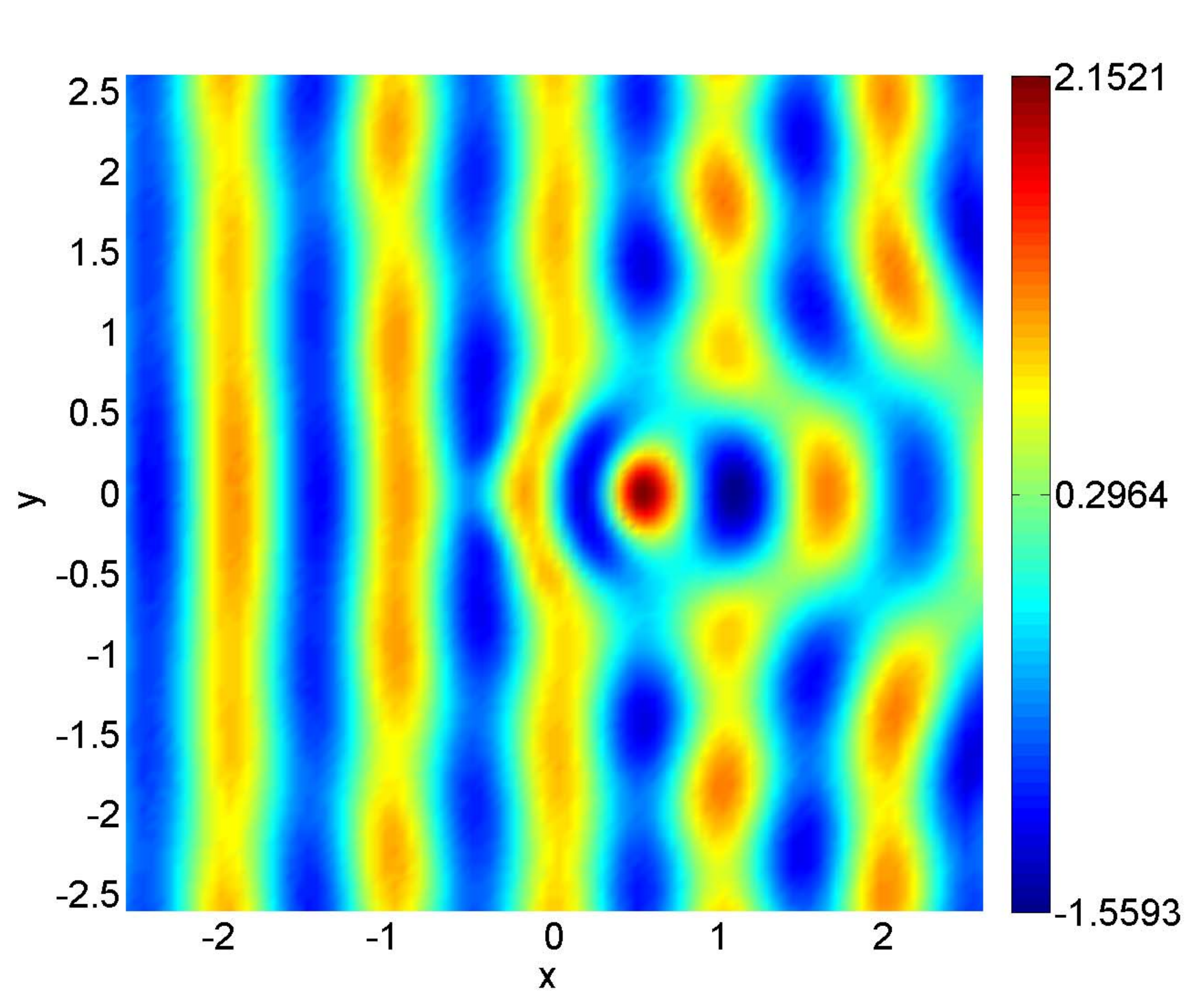}
}
\hspace{0.15cm}
\subfigure
{
\includegraphics[width=2.2in]{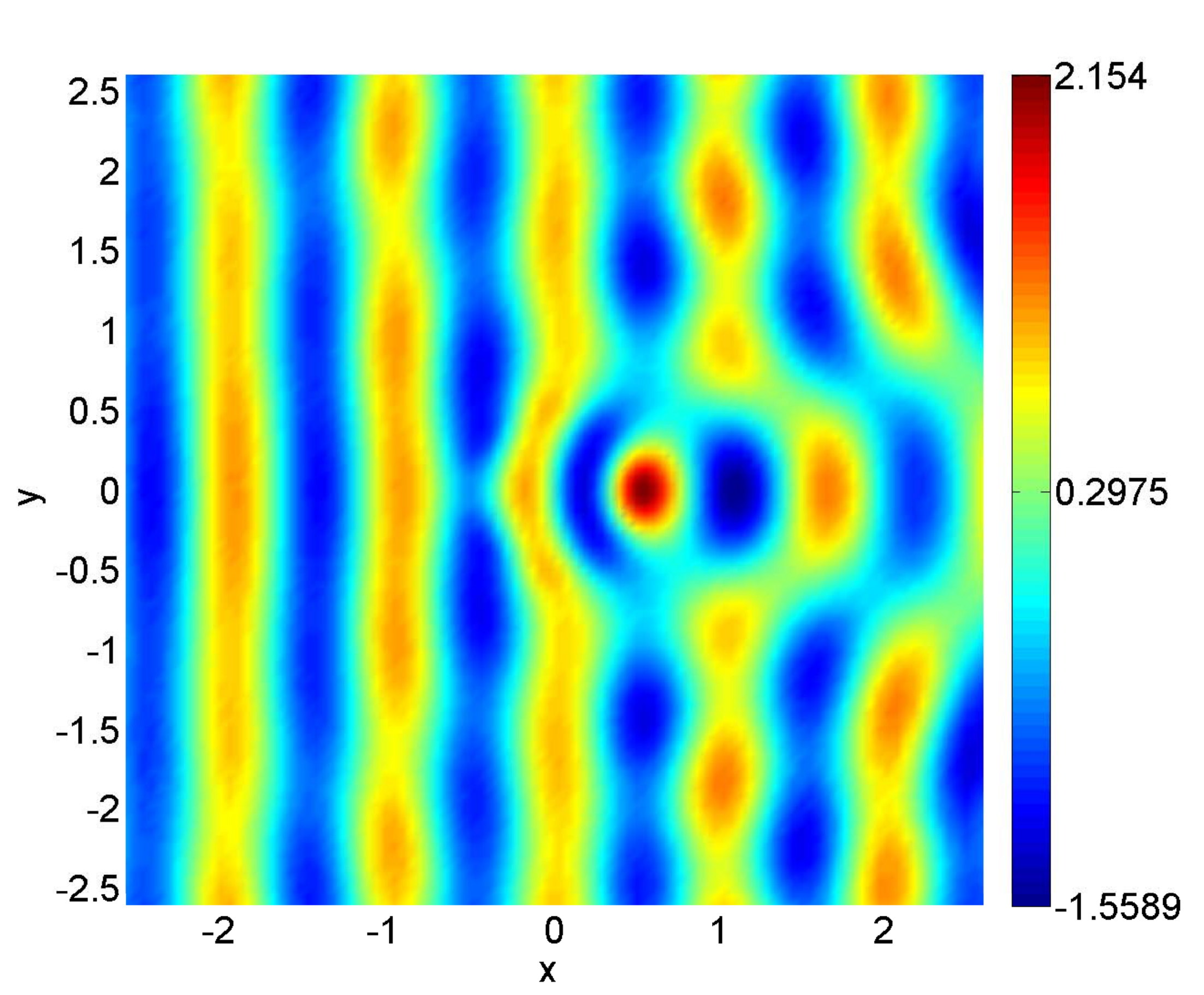}
}\\
\subfigure
{
\includegraphics[width=2.2in]{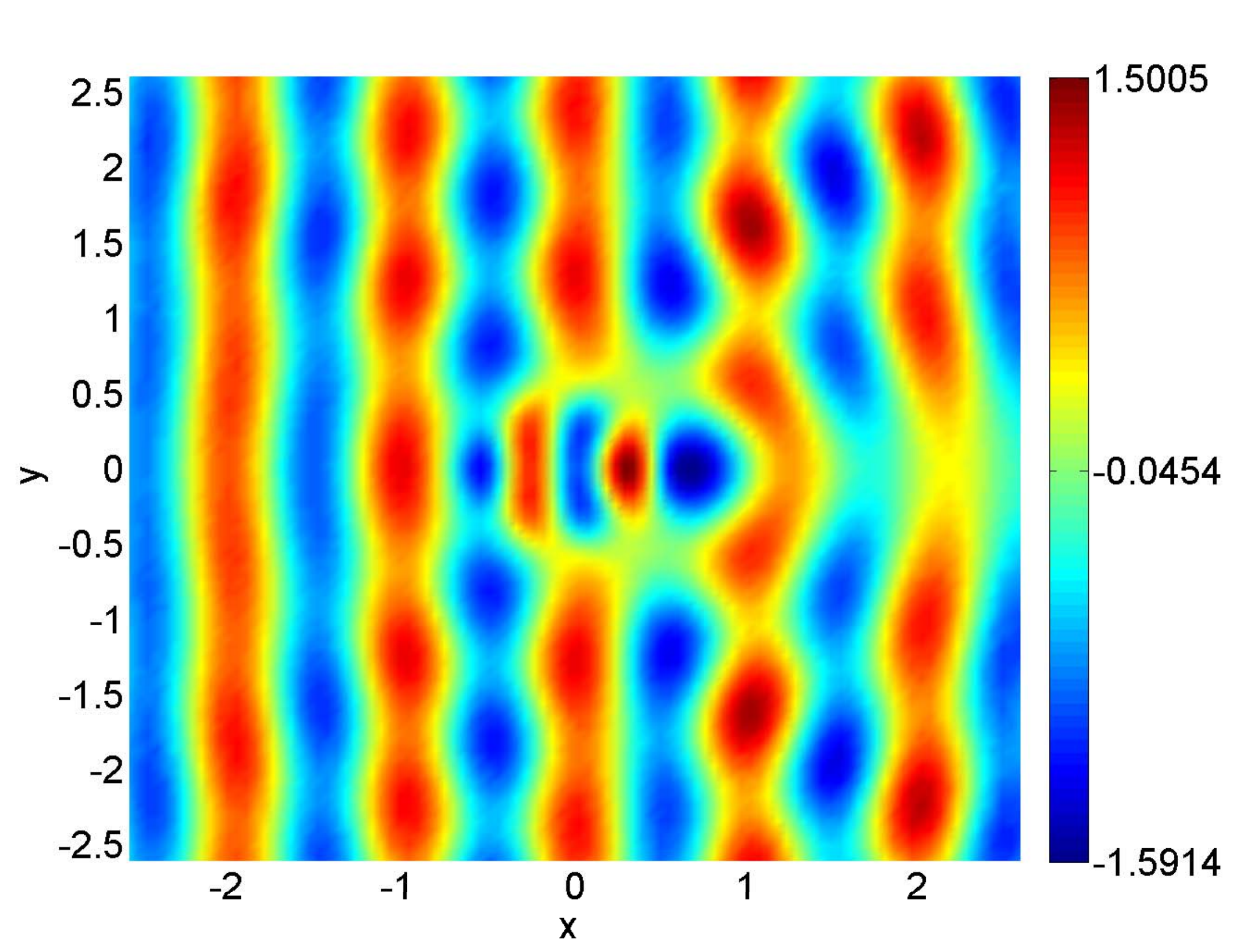}
}
\hspace{0.15cm}
\subfigure
{
\includegraphics[width=2.2in]{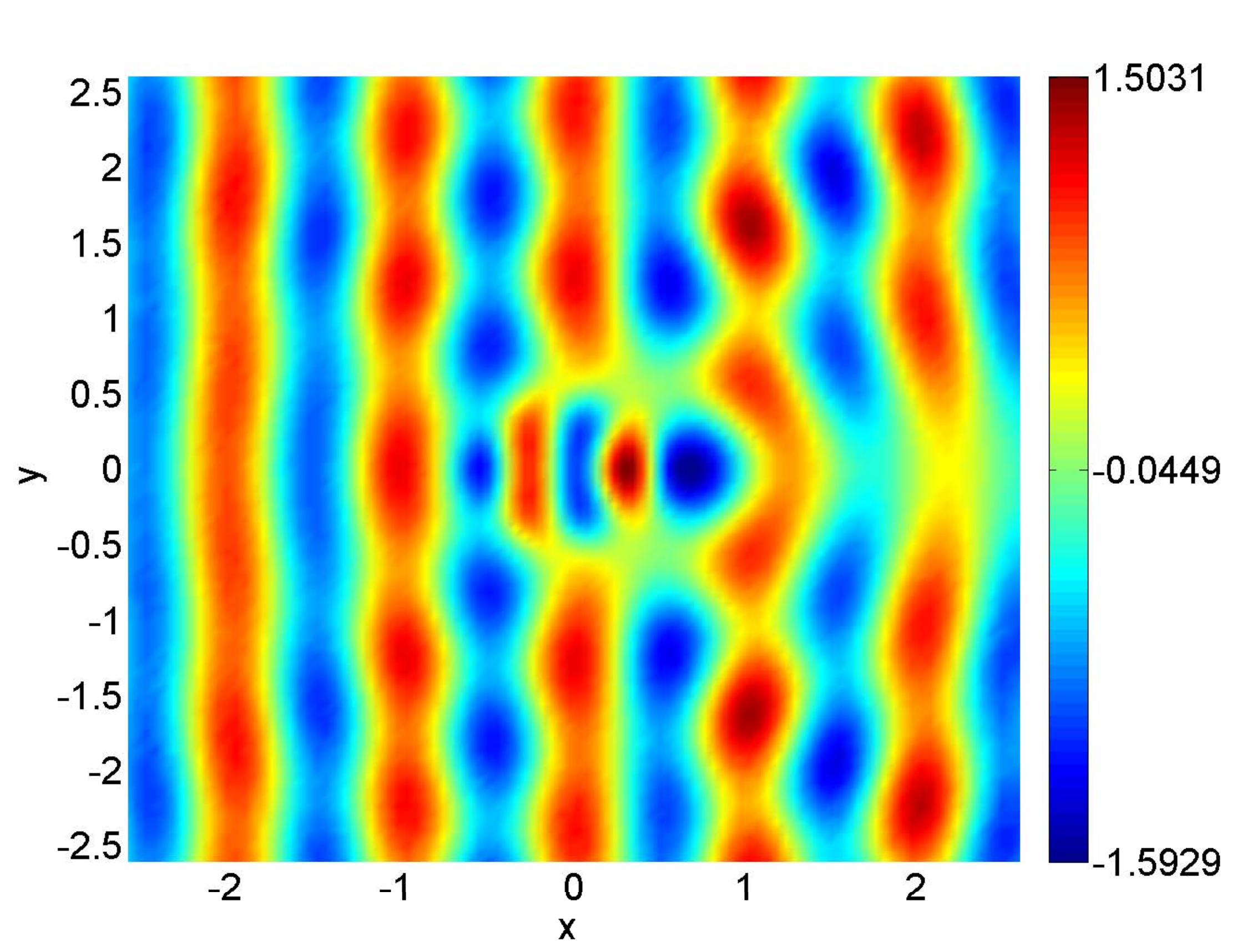}
}\\
\subfigure
{
\includegraphics[width=2.2in]{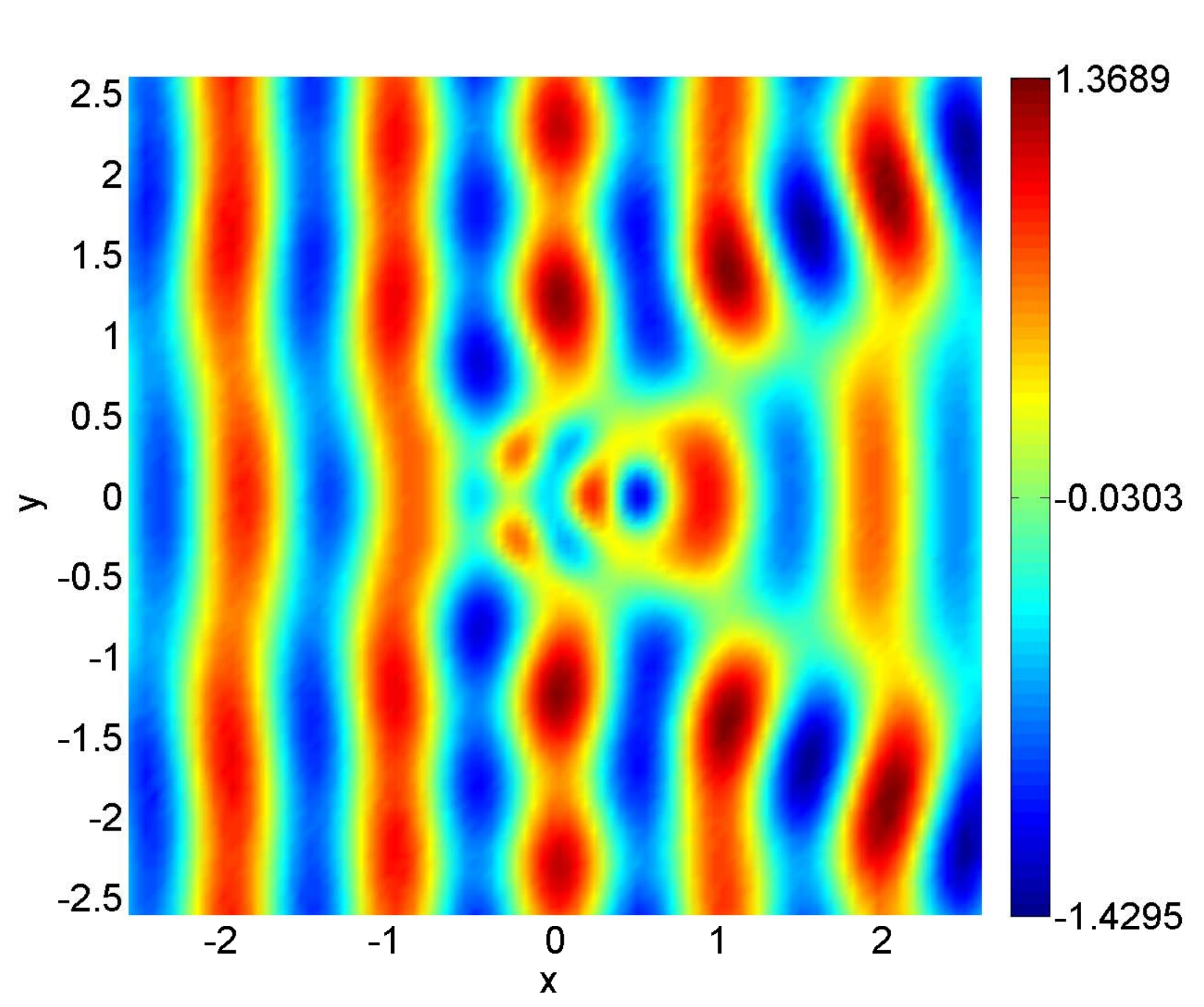}
}
\hspace{0.15cm}
\subfigure
{
\includegraphics[width=2.2in]{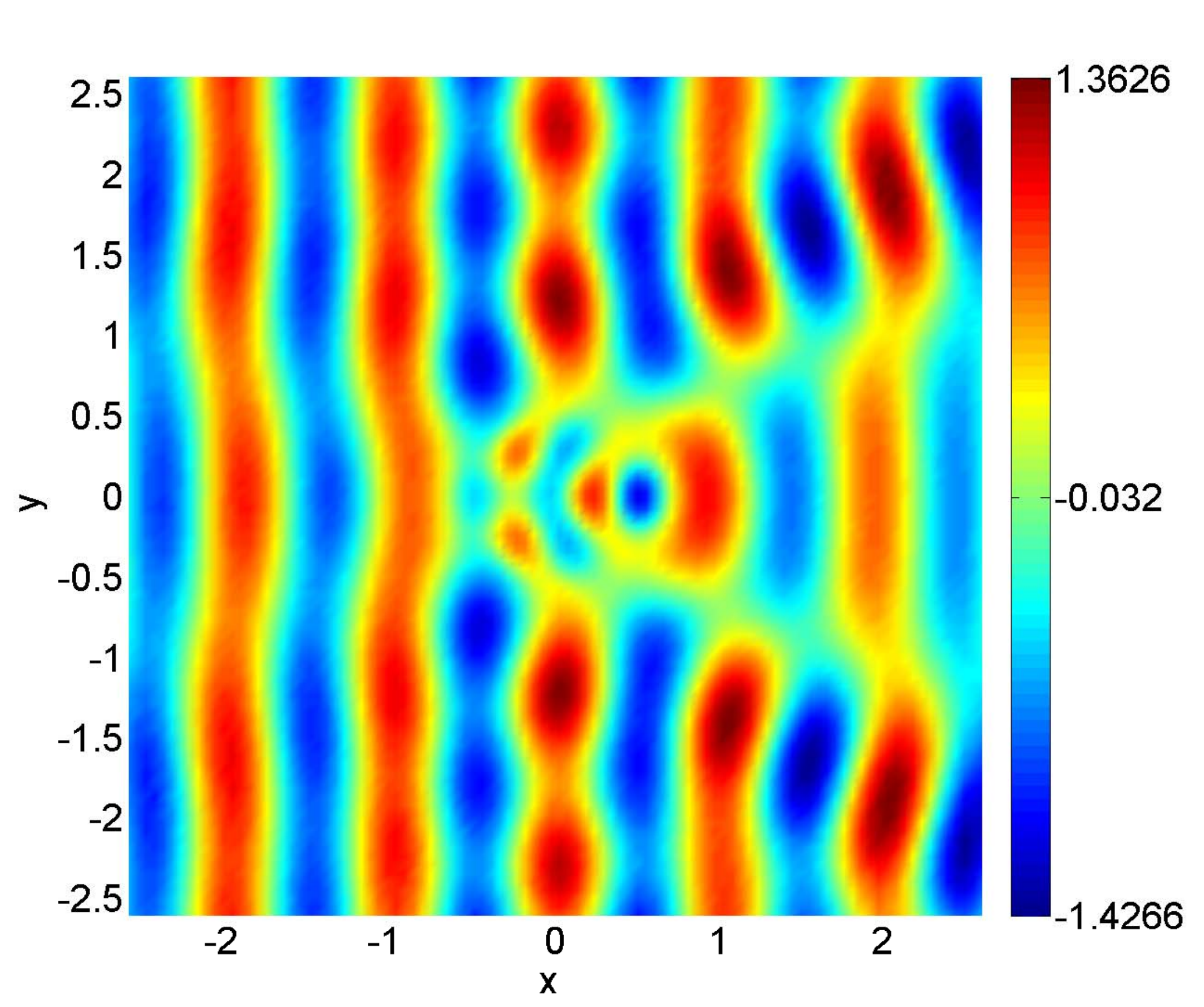}
}\\
\caption{Scattering of a plane wave by a dielectric cylinder: Comparison of the 2D distribution of the real part of $E_z$ between DGTD(left) and POD-GPR (right) of four test points: $\varepsilon_r=1.215$ (1st row), $\varepsilon_r=2.215$ (2nd row), $\varepsilon_r=3.215$ (3rd row) and $\varepsilon_r=4.215$ (4th row).}
\label{fig:7} 
\end{figure}

\begin{figure}
\centering
\subfigure
{
\includegraphics[width=2.2in]{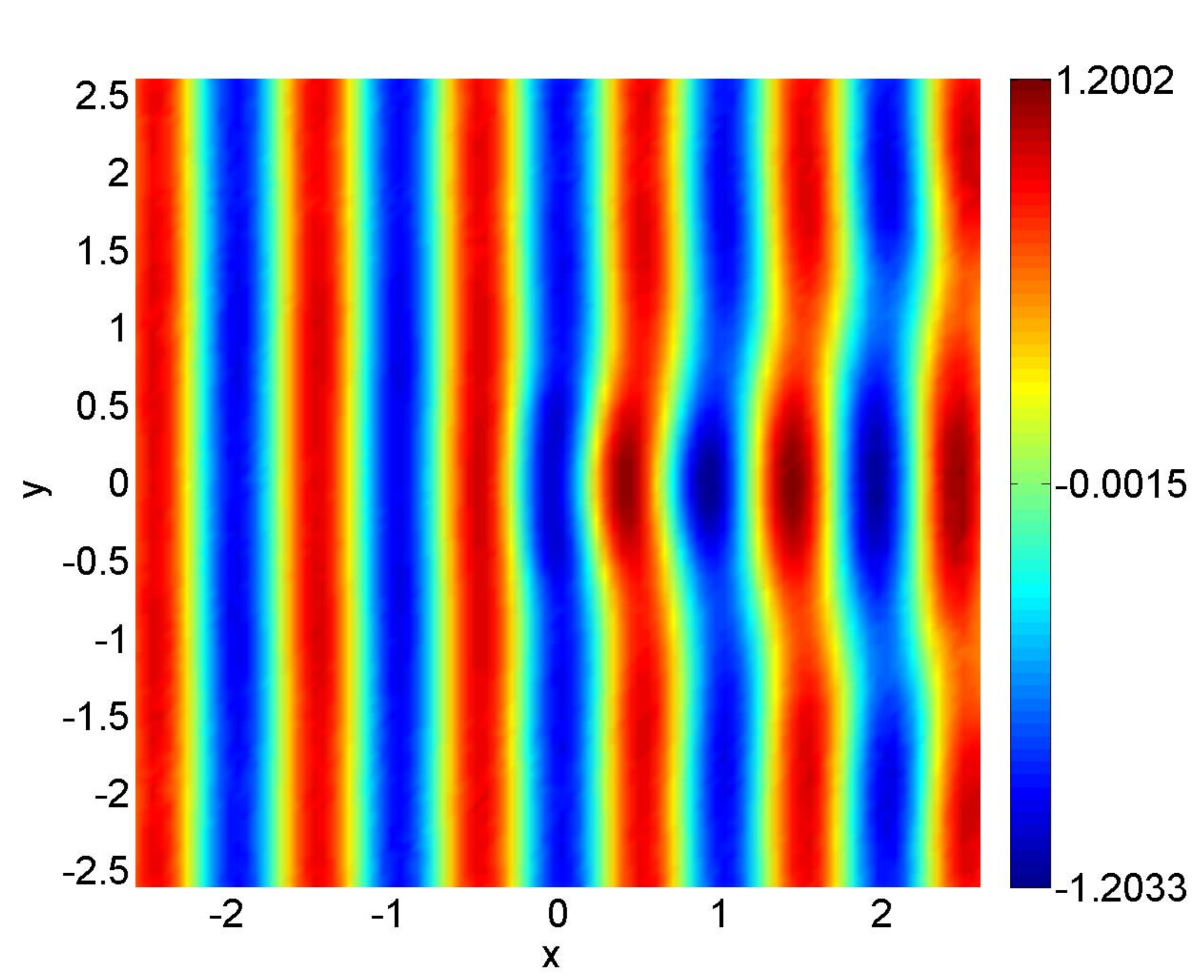}
}
\hspace{0.15cm}
\subfigure
{
\includegraphics[width=2.2in]{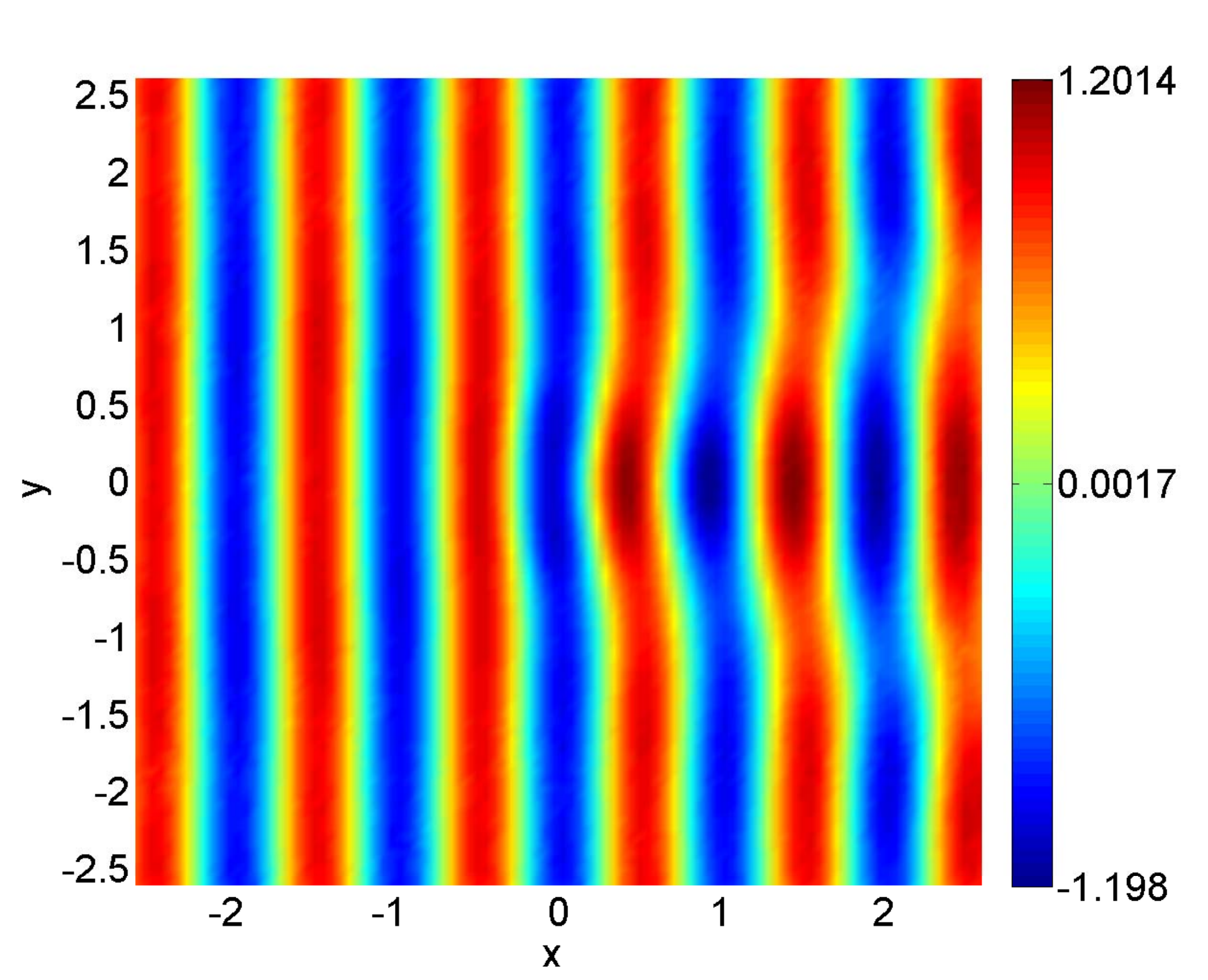}
}\\
\subfigure
{
\includegraphics[width=2.2in]{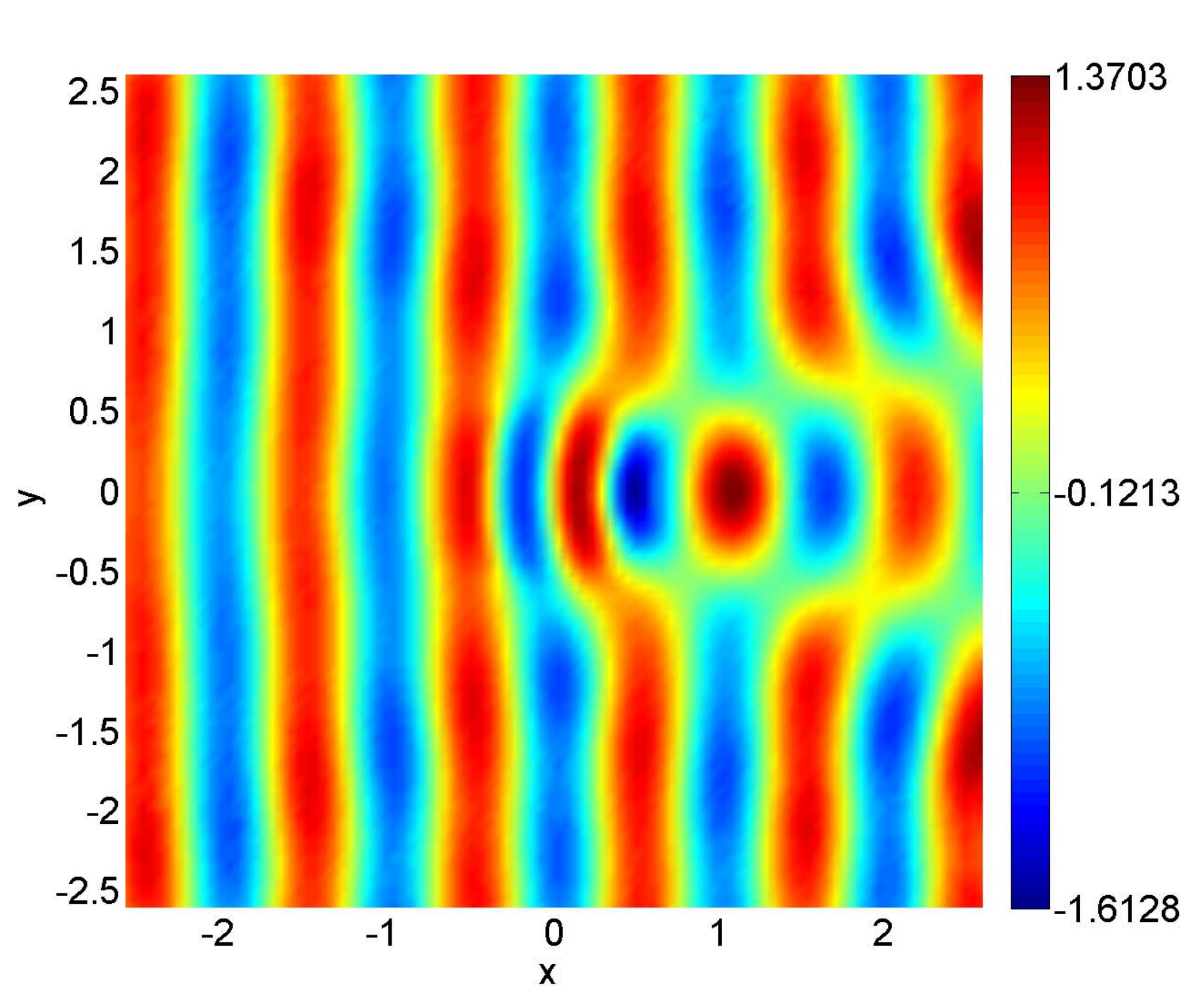}
}
\hspace{0.15cm}
\subfigure
{
\includegraphics[width=2.2in]{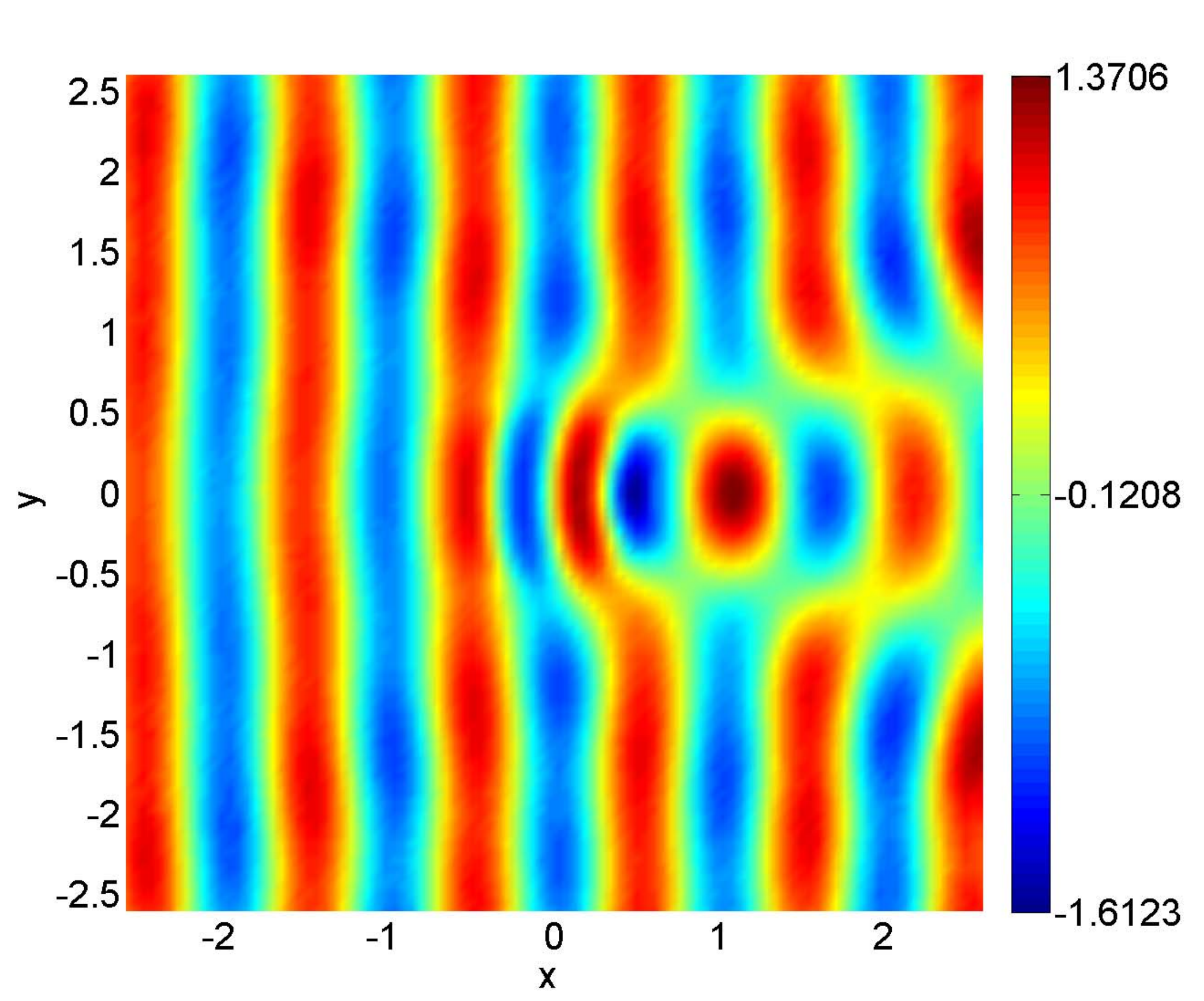}
}\\
\subfigure
{
\includegraphics[width=2.2in]{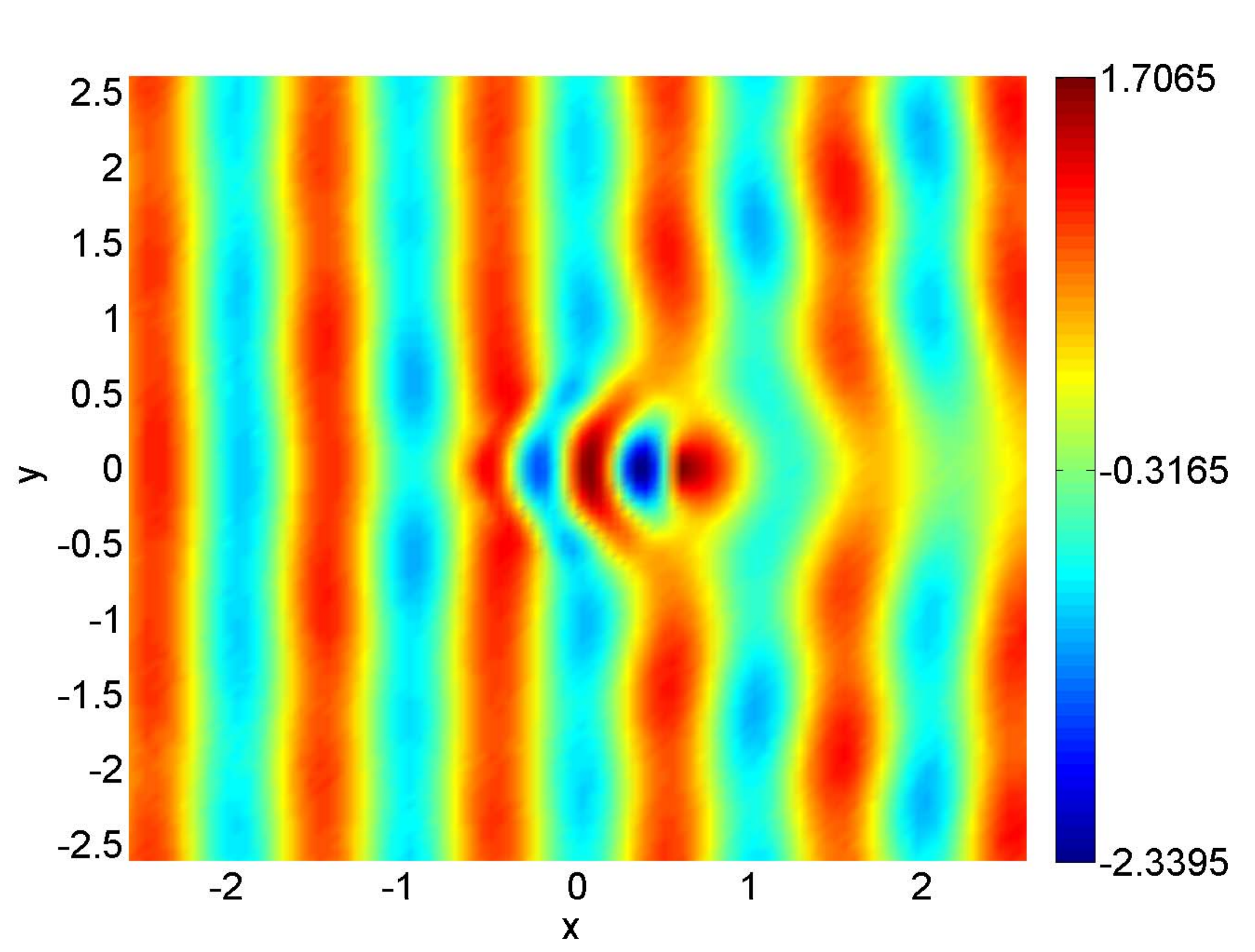}
}
\hspace{0.15cm}
\subfigure
{
\includegraphics[width=2.2in]{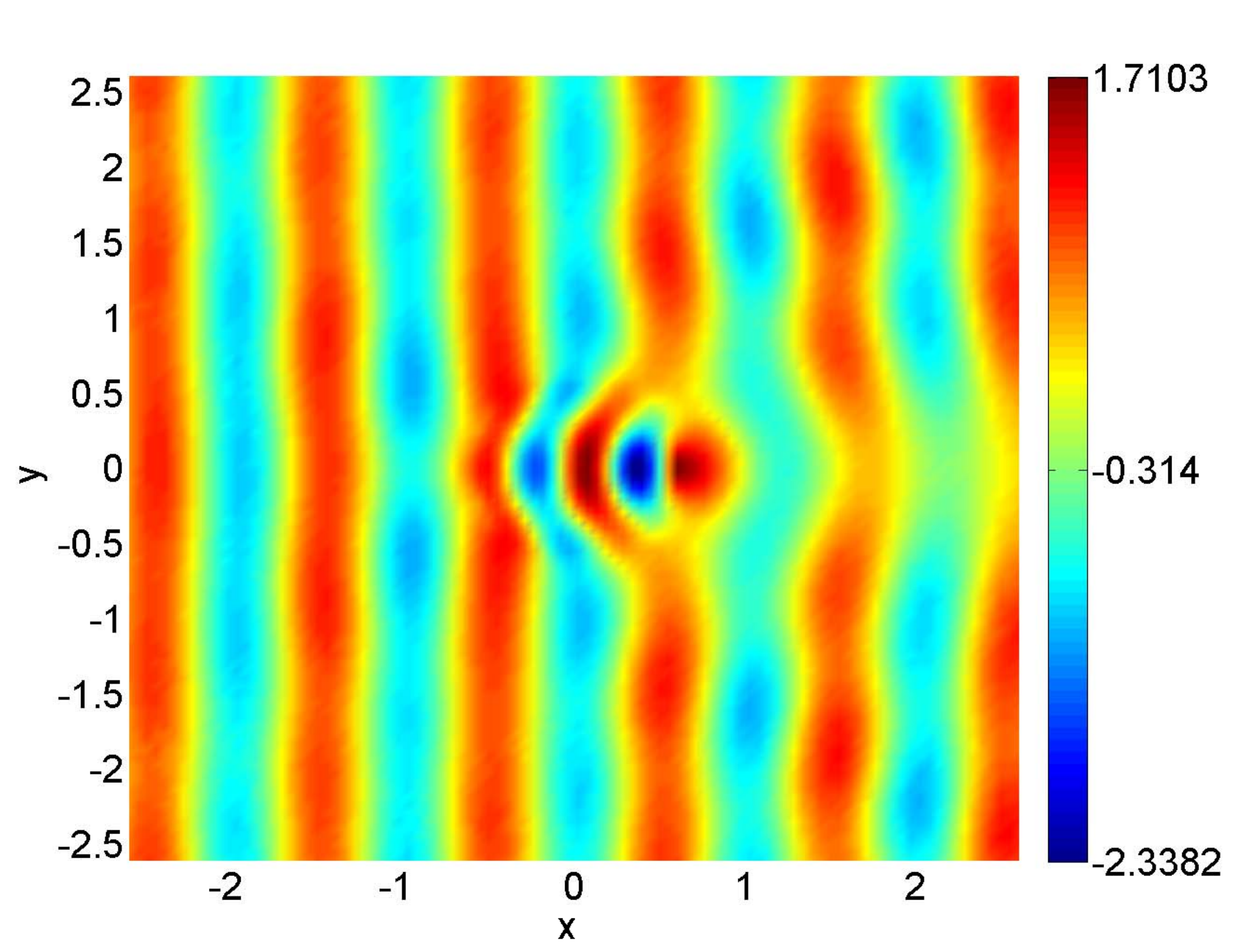}
}\\
\subfigure
{
\includegraphics[width=2.2in]{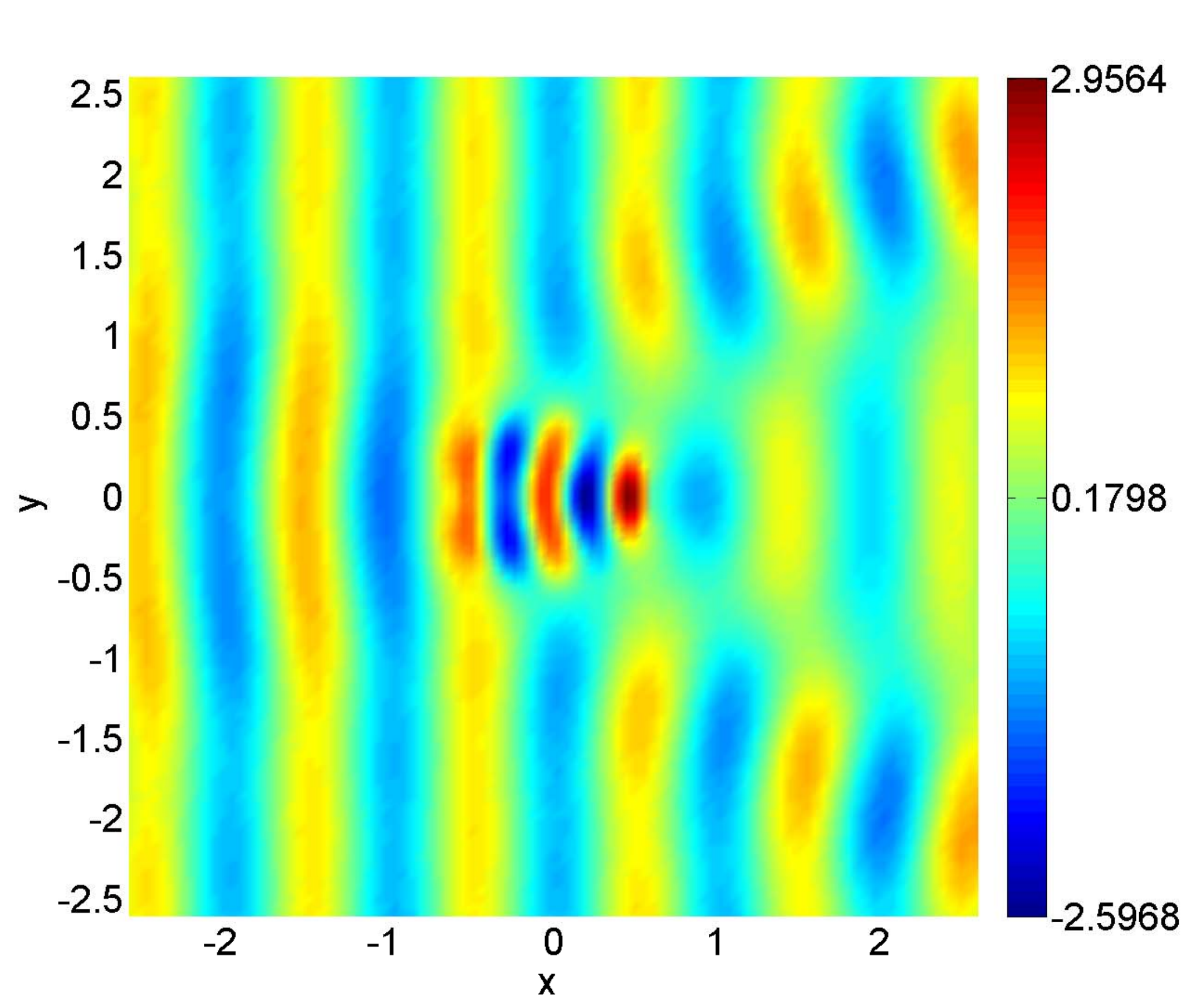}
}
\hspace{0.15cm}
\subfigure
{
\includegraphics[width=2.2in]{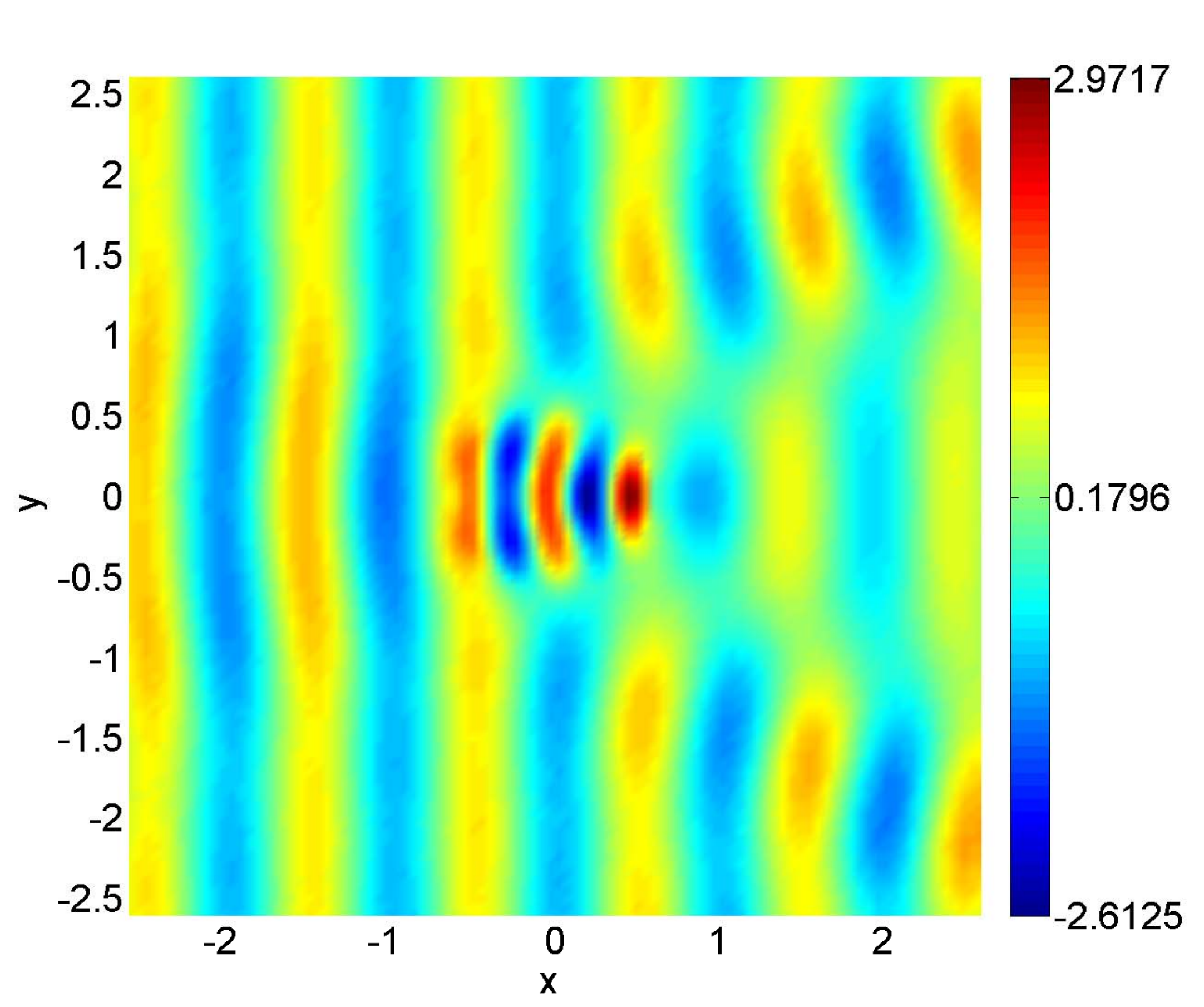}
}
\caption{Scattering of a plane wave by a dielectric cylinder: Comparison of the 2D distribution of the real part of $H_y$ between DGTD (left) and POD-GPR (right) of four test points: $\varepsilon_r=1.215$ (1st row), $\varepsilon_r=2.215$ (2nd row), $\varepsilon_r=3.215$ (3rd row) and $\varepsilon_r=4.215$ (4th row).}
\label{fig:8} 
\end{figure}

Firstly, the time evolution of $E_z$ and $H_y$ at a given point are compared in Fig.\ref{fig:5},
and in order to have a look on visual effects of electromagnetic fields, over the Fourier domain during the last oscillation period of the incident wave, we display in Fig.\ref{fig:6} the 1D x-wise appearance along $y=0$ of the real part of $E_z$ and $H_y$, plus their 2D distribution in Fig.\ref{fig:7} and Fig.\ref{fig:8}. As can be observed, the reduced-order solutions and the DGTD solutions are matching well with each other.

Secondly, a time performance comparison between the DGTD and the POD-GPR is given in Table \ref{tab:2}, 
\begin{table}[]
\caption{Scattering of a plane wave by a dielectric cylinder: Time performance comparison.}
\label{tab:2}
\centering
\begin{tabular}{cc}
\hline
Name& Time/s\\
\hline
Average DGTD solving time& $3.98\times10^{2}$\\
Average POD-GPR solving time& $1.19\times10^{1}$\\
GPR training time& $8.60\times10^{2}$\\
\hline
\end{tabular}
\end{table}
where we record the average consuming time of DGTD solver and ROM online output for above four test items, as well as the construction time of the GPR model. Although it takes a long time to build the new model, it is worthwhile especially when solving the electromagnetic field value under multiple parameters, since the online output time of ROM is greatly shortened compared with the DGTD solver, which demonstrates the appreciably improved efficiency of the approach employed herein.

\begin{figure}
\centering
\subfigure
{
\includegraphics[width=2.2in]{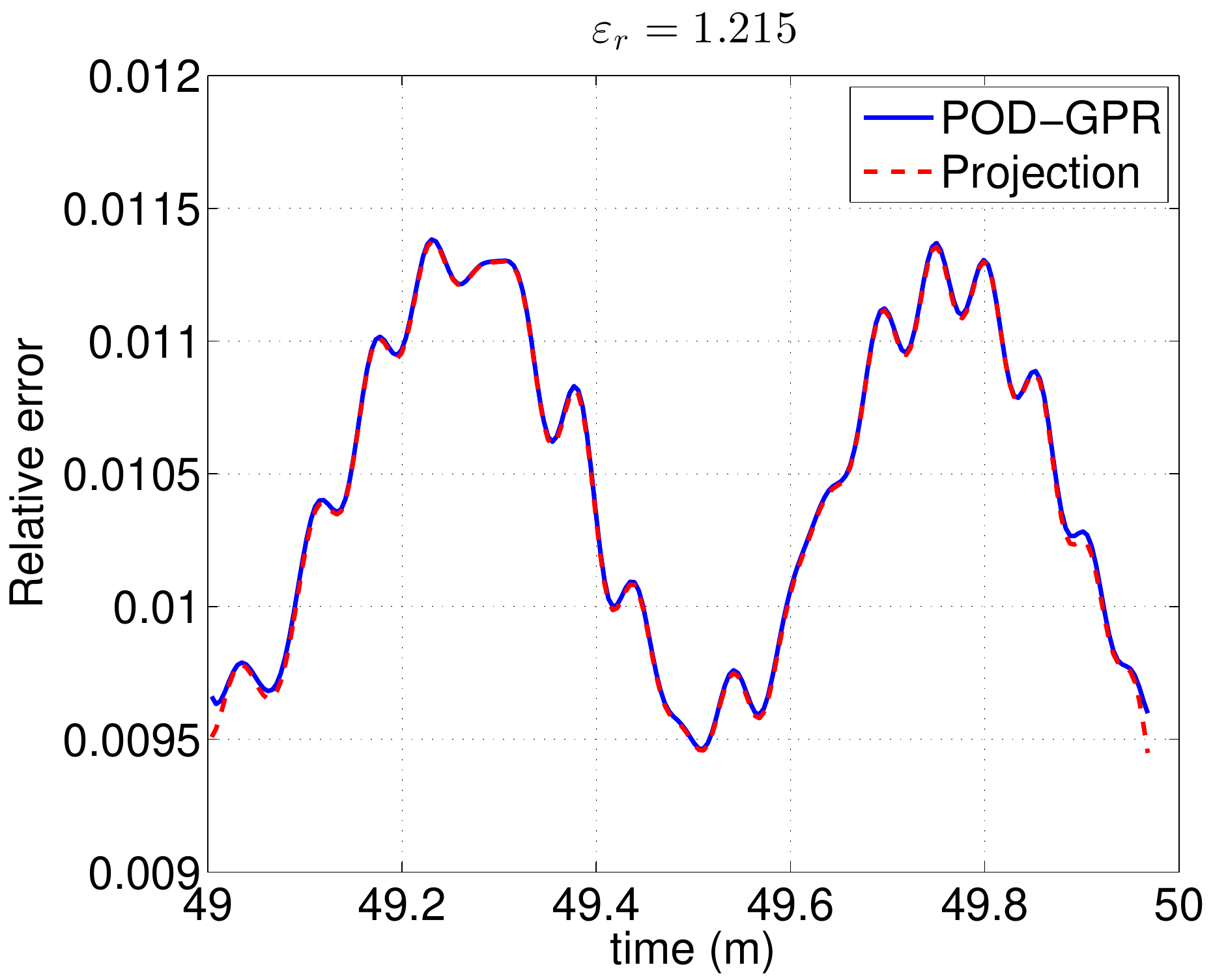}
}
\hspace{0.15cm}
\subfigure
{
\includegraphics[width=2.2in]{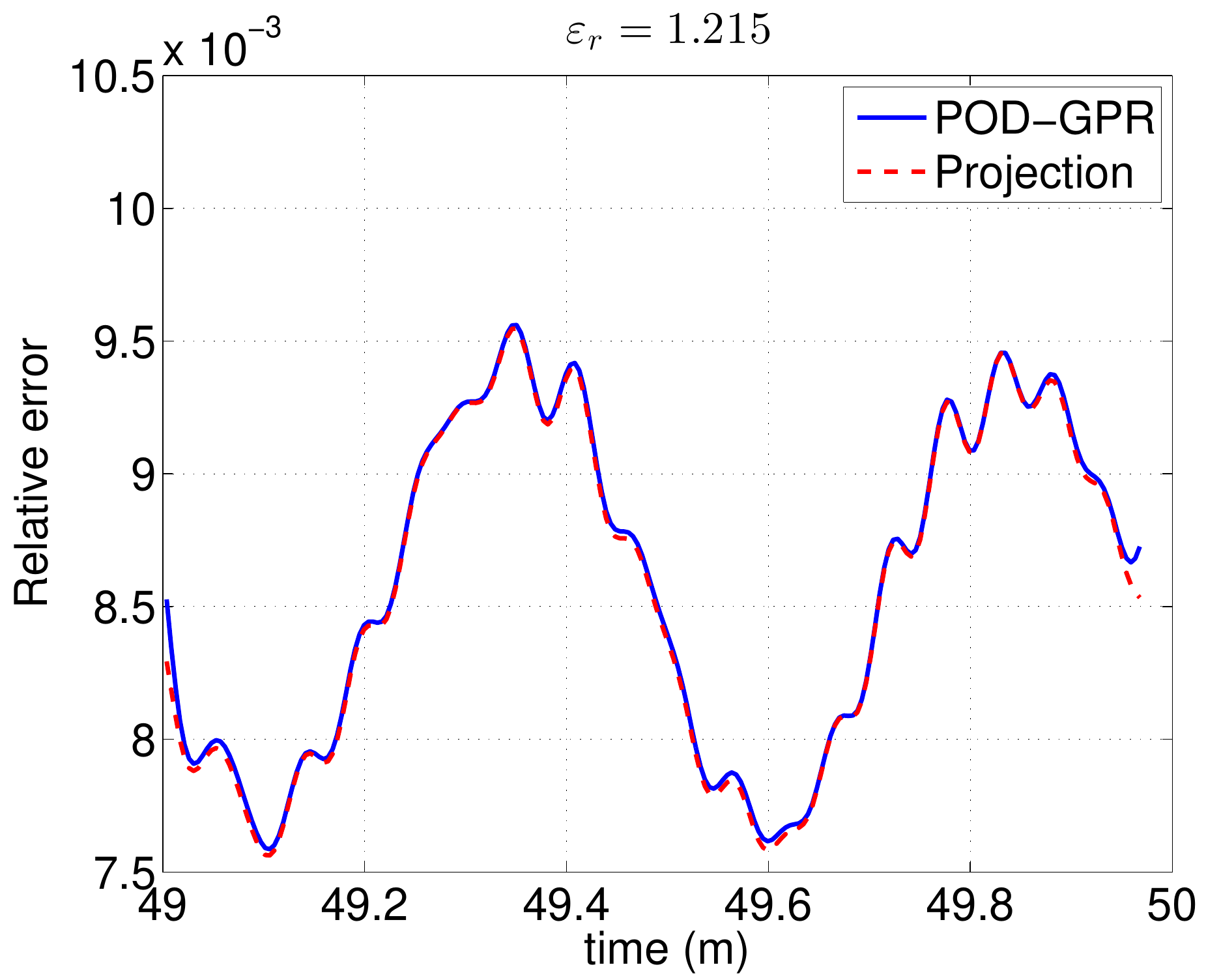}
}\\
\subfigure
{
\includegraphics[width=2.2in]{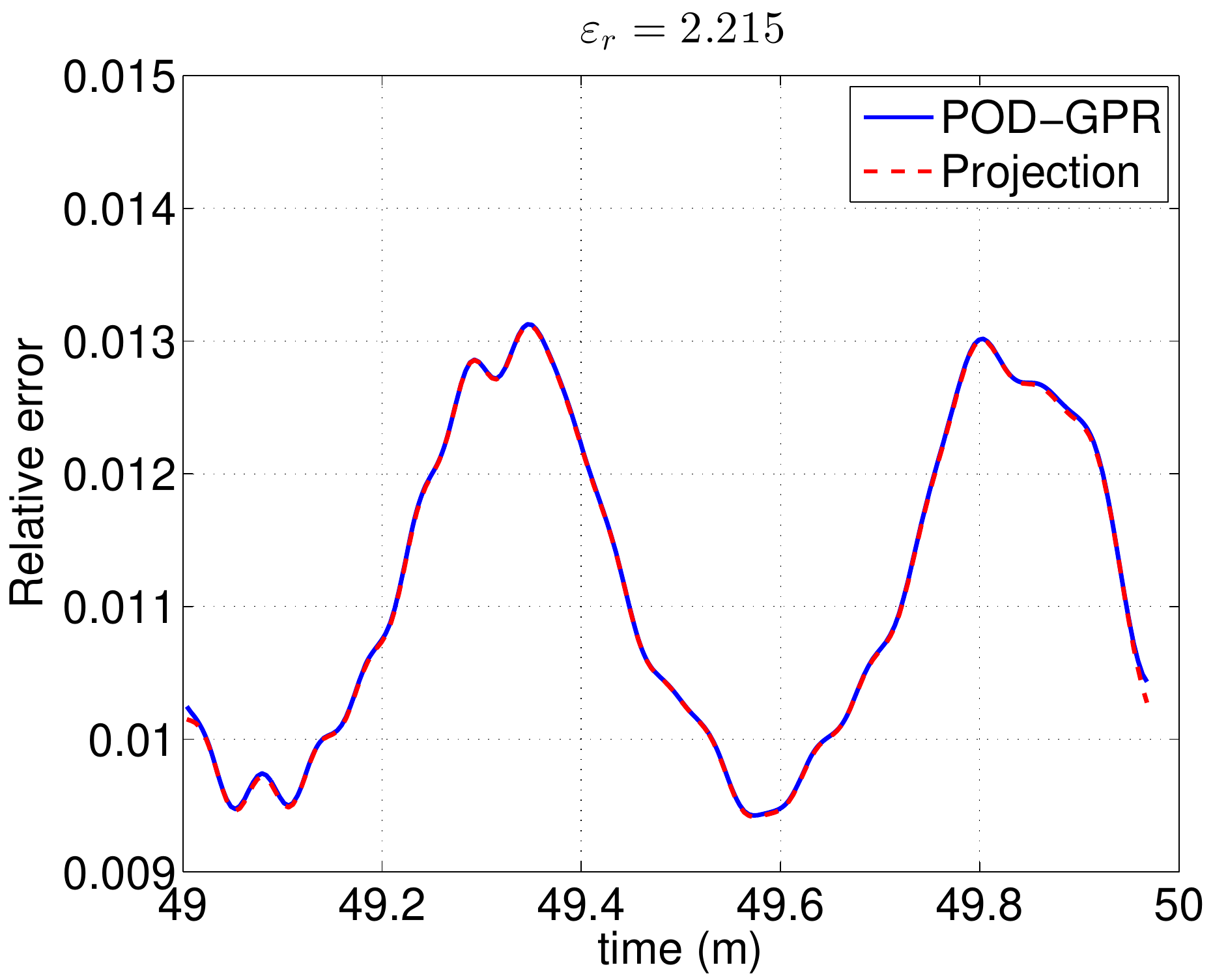}
}
\hspace{0.15cm}
\subfigure
{
\includegraphics[width=2.2in]{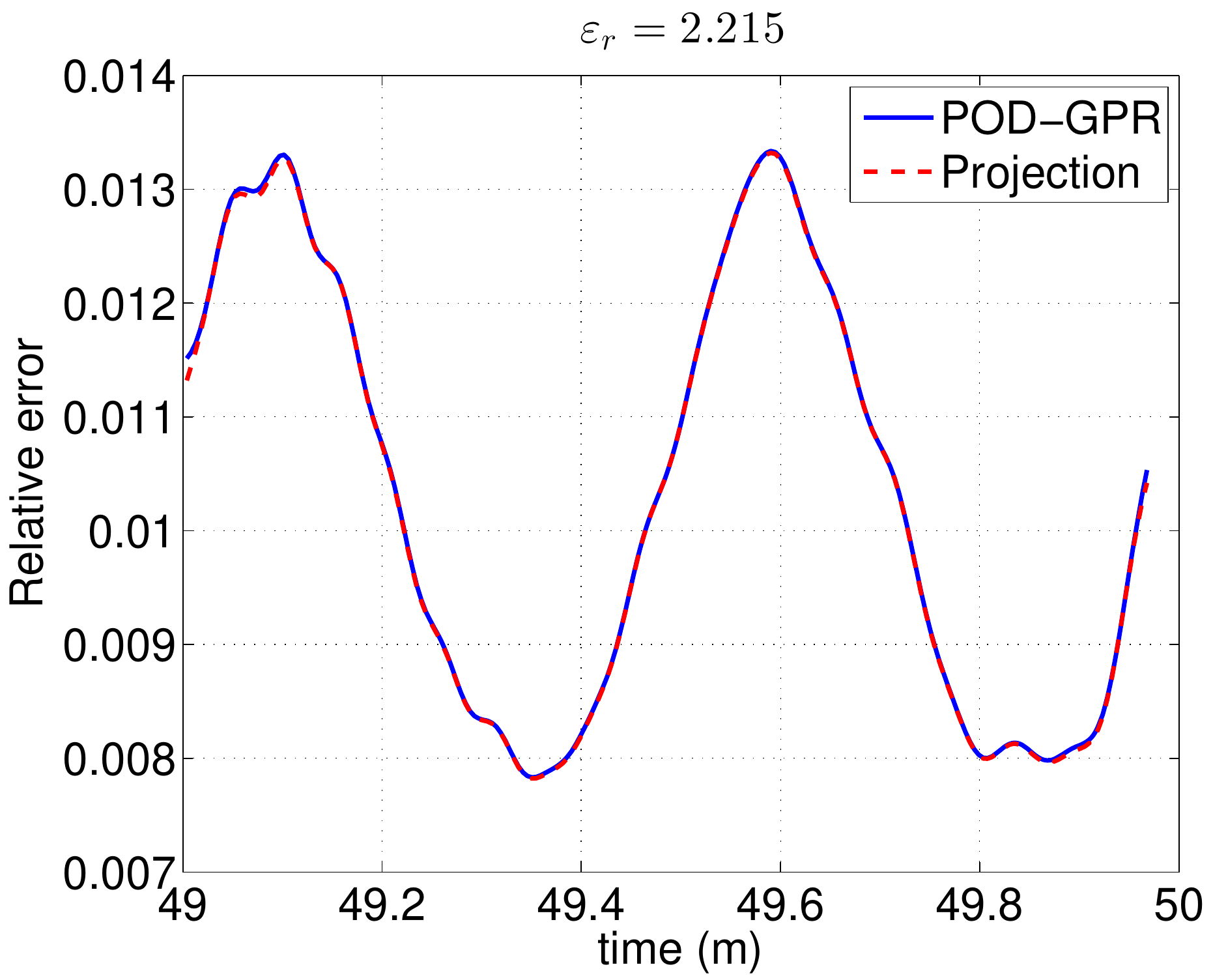}
}\\
\subfigure
{
\includegraphics[width=2.2in]{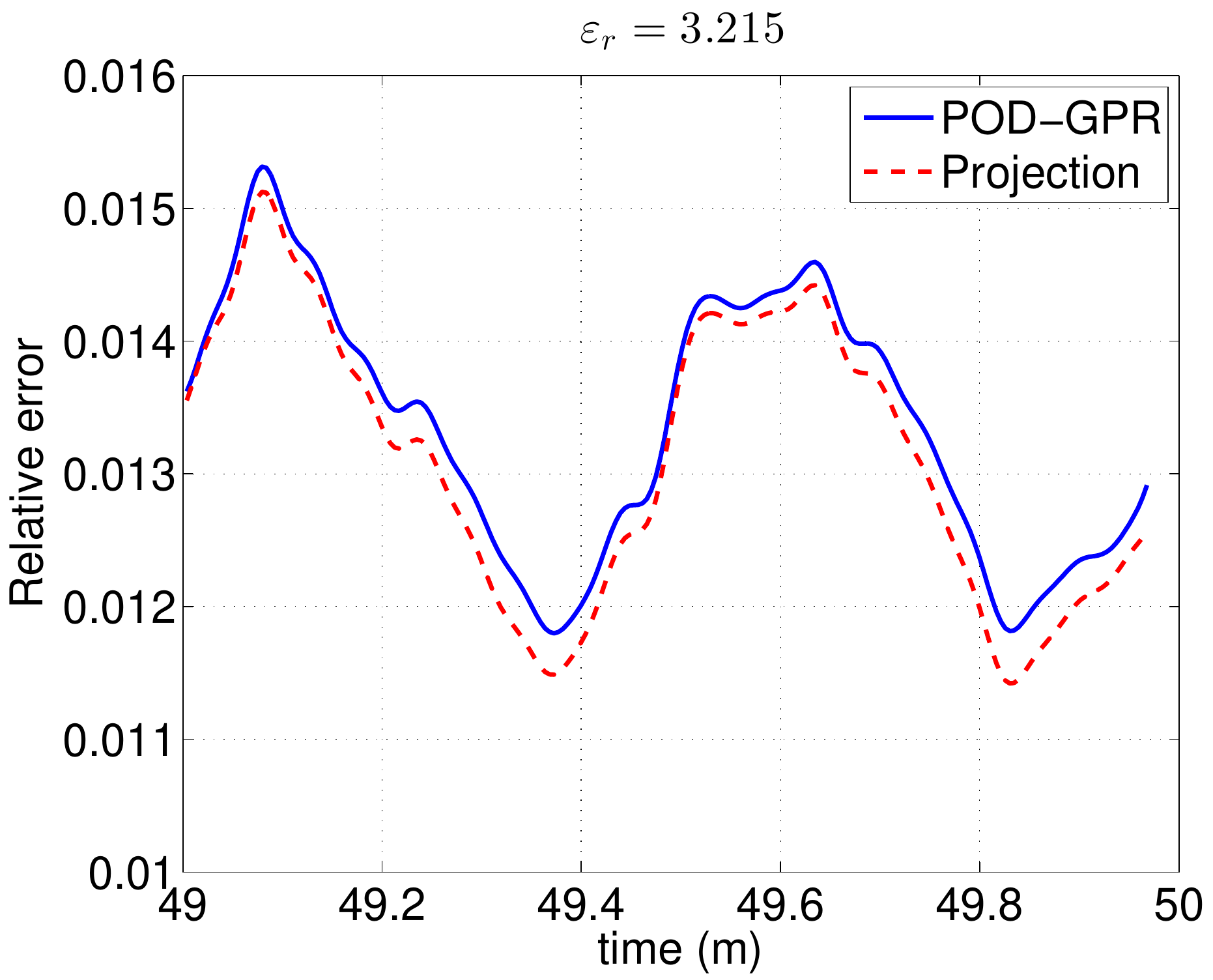}
}
\hspace{0.15cm}
\subfigure
{
\includegraphics[width=2.2in]{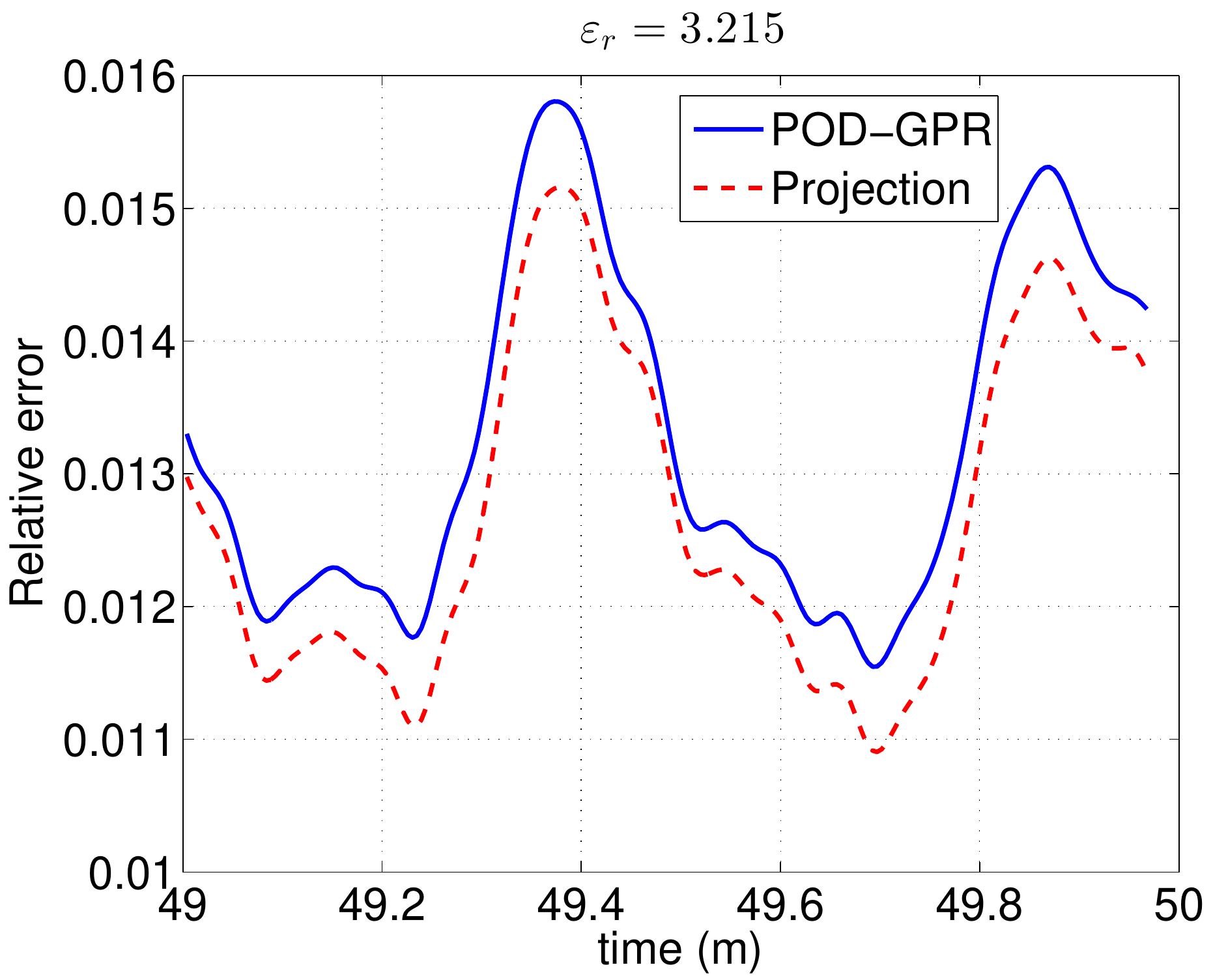}
}\\
\subfigure
{
\includegraphics[width=2.2in]{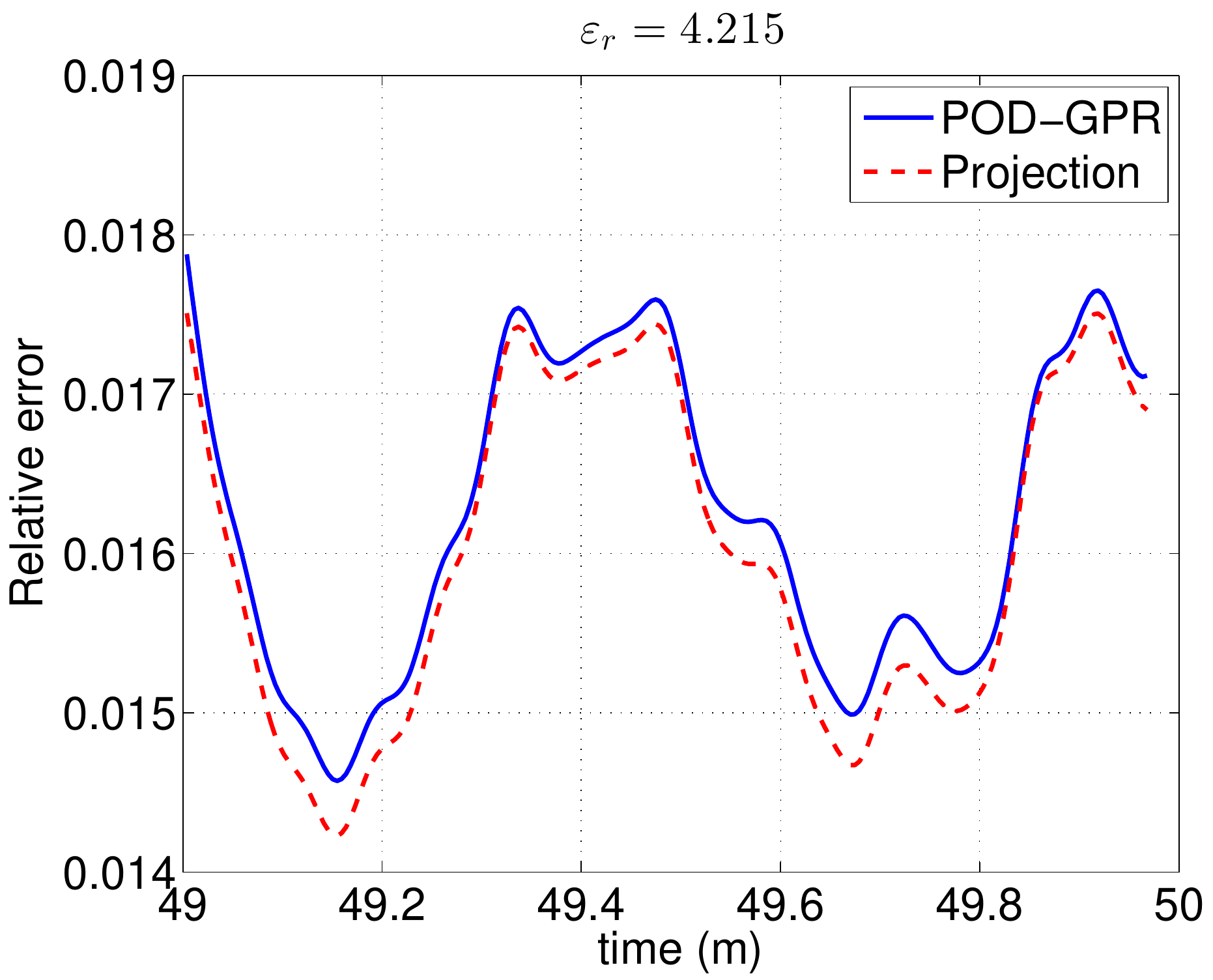}
}
\hspace{0.15cm}
\subfigure
{
\includegraphics[width=2.2in]{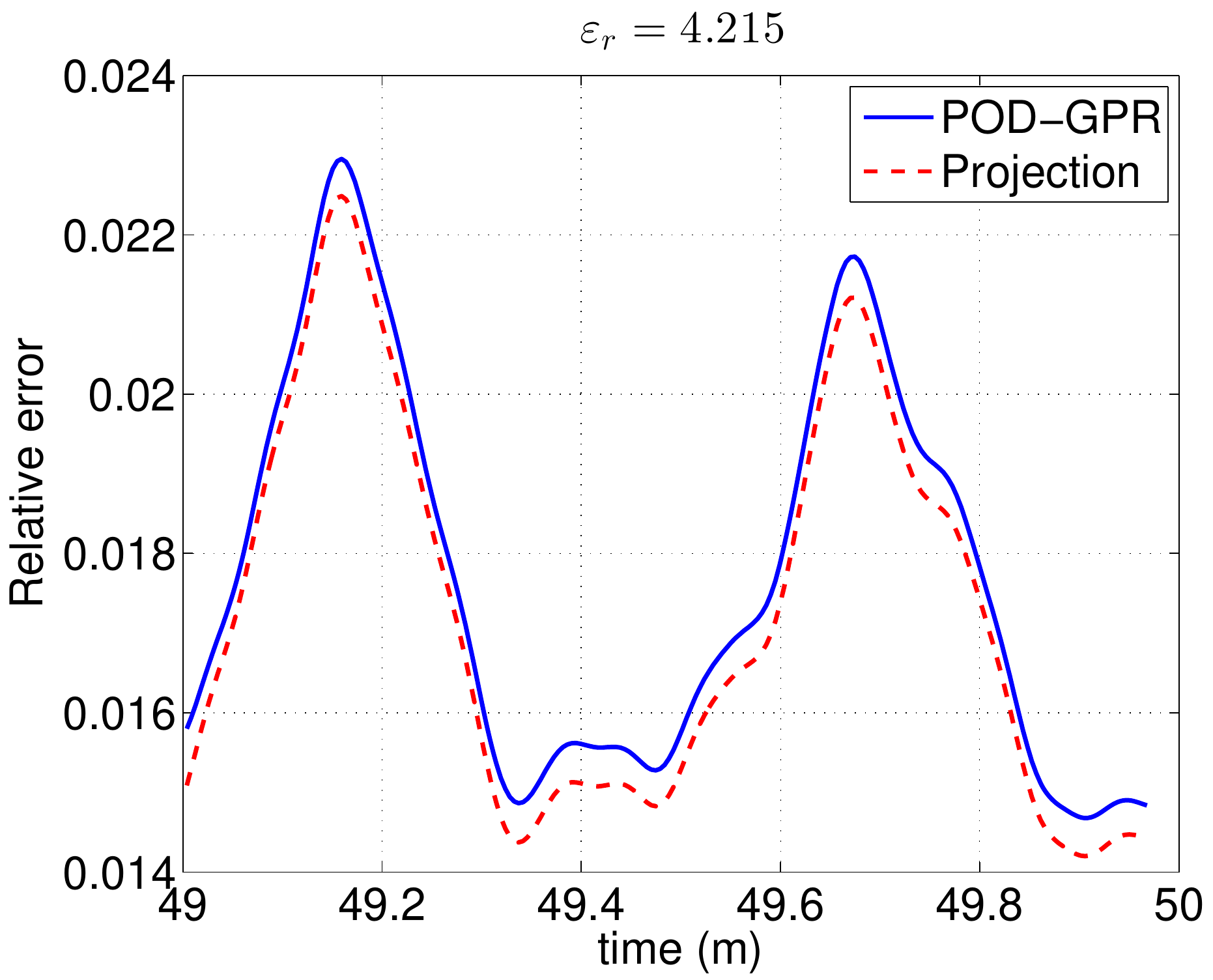}
}\\
\caption{Scattering of a plane wave by a dielectric cylinder: Comparison of the relative $L^2$ error between POD-GPR and DGTD for $E_z$ (left) and $H_y$ (right) of four test points.}
\label{fig:9} 
\end{figure}

Furthermore, the time evolution of the relative $L^2$ error between POD-GPR and DGTD is exhibited in Fig.\ref{fig:9}, 
which also contains the time evolution of the relative $L^2$ projection error generated by POD. Note that the curves of two errors are very close, meaning that the error caused by GPR is negligible and confirming its accuracy.
\subsection{Scattering of a plane wave by a multi-layer heterogeneous medium}
\begin{figure}[htbp]
\centering
\includegraphics[height=5.0cm,width=7.8cm]{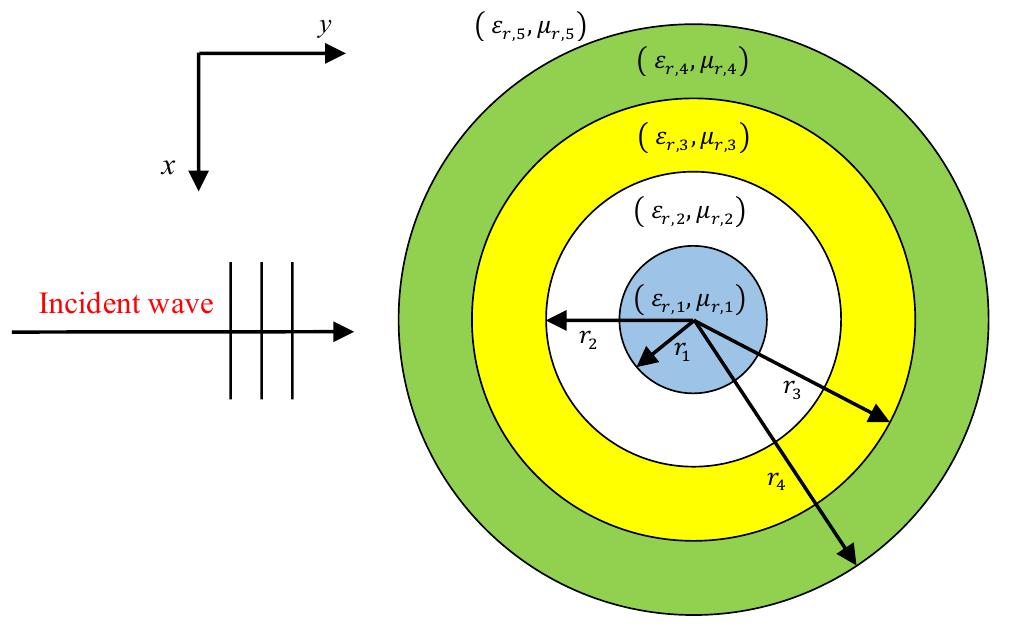}
\caption{Scattering of a plane wave by a multi-layer heterogeneous medium: Geometry of the multi-layer heterogeneous medium}
\label{fig:10}
\end{figure}

In this section, we are concerned with a more complex situation where a multi-heterogeneous medium, as is shown in Fig.\ref{fig:10}, 
is illuminated by an incident plane wave. The geometry of the domain is a square  $\Omega = [-3.2m,3.2m]\times[-3.2m,3.2m]$, with the ABC condition be enforced on its boundary. As with the last example, the external medium is vacuum, i.e., $\varepsilon_{r,5}=1$ and $\mu_{r,5}=1$. Besides, we only involve the nonmagnetic materials, i.e., $\mu_{r,i}=1,i=1,\cdots,4$. The size of each medium layer and their relative permittivity range are summarized in Table \ref{tab:3}, which also includes sampling methods for parameters.
\begin{table}[]
\caption{Scattering of a plane wave by a multi-layer heterogeneous medium: Physical and Sampling information of medium.}
\label{tab:3}
\centering
\begin{tabular}{ccccc}
\hline
Layer $i$ & 1 & 2 & 3 & 4\\
\hline
$r_i$ (m)& 0.15 & 0.3 & 0.45 & 0.6  \\
Range of $\varepsilon_{r,i}$ & [5.0,5.6]&[3.25,3.75]&[2.0,2.5]&[1.25,1.75]\\
Sampling method&3, uniform&3, uniform&3, uniform&3, uniform\\
\hline
\end{tabular}
\end{table}
Therefore, the parameter $\theta$ can be represented as a $4D$ vector $\theta=[\varepsilon_{r,1},\varepsilon_{r,2},\varepsilon_{r,3},\varepsilon_{r,4}]$, and $\mathcal{P}=[5.0,5.6]\times[3.25,3.75]\times[2.0,2.5]\times[1.25,1.75]$.

The computational mesh consists of 3256 nodes and 6206 elements, with 118 elements located inside the first layer, 308 in the second, 476 in the third and 604 in the fourth layer, resulting in $N_h=37236$ DOFs for the DGTD solver.

\begin{figure}[htbp]
\centering
\includegraphics[height=0.395\textheight, width=0.495\textheight]{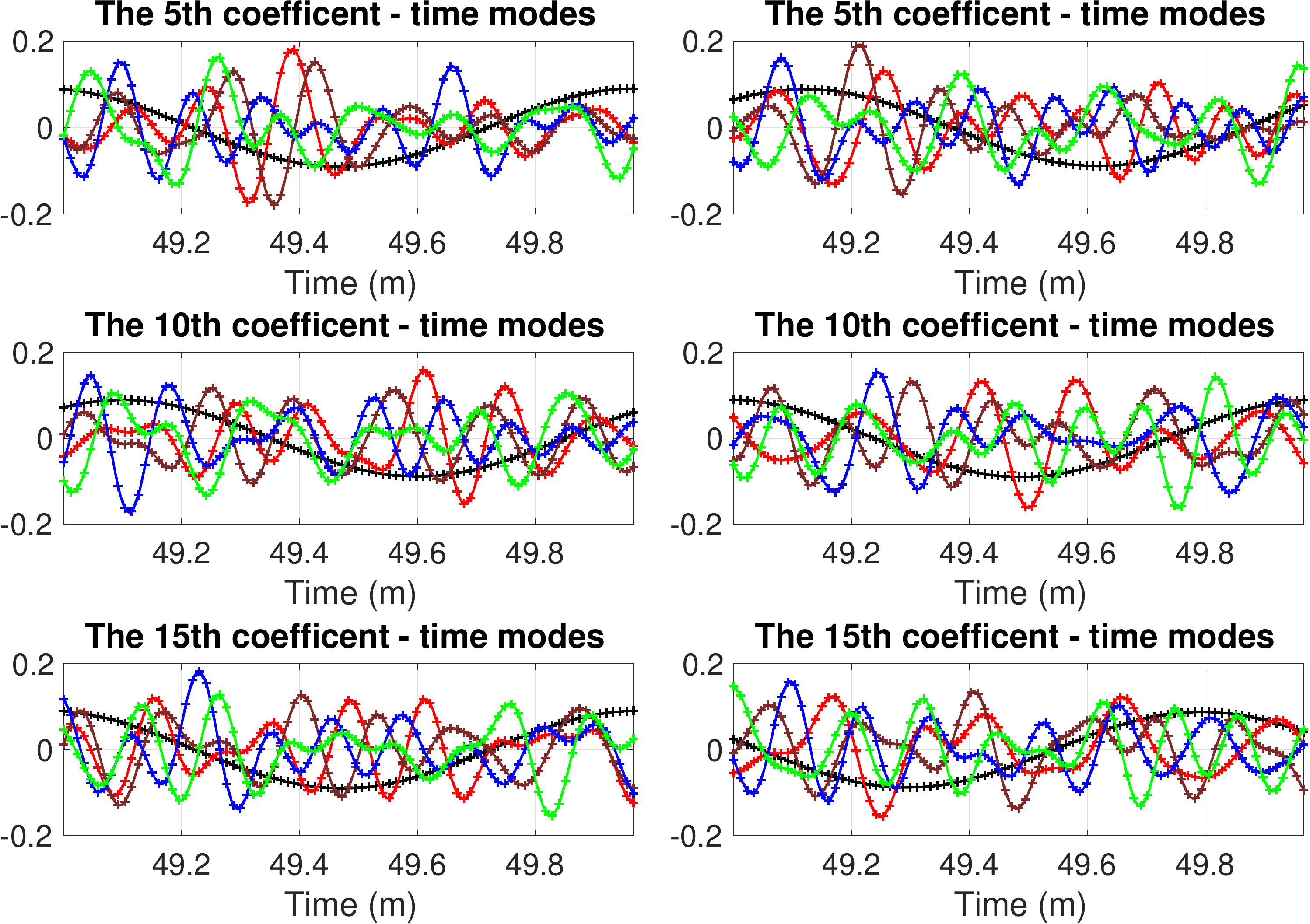}
\caption{Scattering of a plane wave by a multi-layer heterogeneous medium: Time-modes for the 5th, 10th and 15th projection coefficients for $E_z$ (left) and $H_y$ (right): the 2nd modes-black, the 4th modes-red, the 6th modes-brown, the 8th modes-blue, the 10th modes-green. }
\label{fig:11}
\end{figure}

In order to perform the offline preparation, we do some full-order simulations under $N_\theta=81$ parameter values points corresponding to the above sampling methods (i.e., $\mathcal{P}_h=[5.0:0.3:5.6]\times[3.25:0.25:3.75]\times[2.0:0.25:2.5]\times[1.25:0.25:1.75]$), with simulation time being 50 periods of the incident wave oscillation. As for single point of parameter, $N_l=253$ transient full-order solutions are extracted in the last oscillation period both as snapshots and as training data (i.e., $\mathcal{T}_{h}=\{49.0009,49.0042,49.0075,\cdots,49.9692\}$). Then, the reduced spaces $\mathcal{V}_{E_z,rb}$, $\mathcal{V}_{H_y,rb}$ and $\mathcal{V}_{H_x,rb}$ are spanned respectively by $d_{E_z}=15$, $d_{H_y}=15$ and $d_{H_x}=17$ basis functions, given by the two-step POD with $\epsilon_{\textbf{E},t}=\epsilon_{\textbf{H},t}=5\times e^{-4}$ and $\epsilon_{\textbf{E},\theta }=\epsilon_{\textbf{H},t}=1\times e^{-5}$. To apply GPR technique to approximate the map between $(t,\varepsilon_{r,1},\varepsilon_{r,2},\varepsilon_{r,3},\varepsilon_{r,4})$ and projection coefficients, we first perform SVD on training data matrix, and then the discrete time- and parameter-data are used to build GPR models, some of which can be looked up in Fig.\ref{fig:11}. In particular, Table \ref{tab:4} gives the grouping way of training.
\begin{table}[]
\caption{Scattering of a plane wave by a multi-layer heterogeneous medium: The SVD truncation tolerance.}
\label{tab:4}
\centering
\begin{tabular}{cc}
\hline
The projection coefficient item& $\delta_{\textbf{E},l}, \delta_{\textbf{H},l}$\\
\hline
$l\leq2$ & $6 \times 10^{-5}$\\
$2<l\leq5$ & $1 \times 10^{-4}$\\
$5<l\leq10$ & $5 \times 10^{-4}$\\
$10<l$ & $1 \times 10^{-3}$\\
\hline
\end{tabular}
\end{table}

\begin{figure}
\centering
\subfigure[]
{
\includegraphics[width=2.2in]{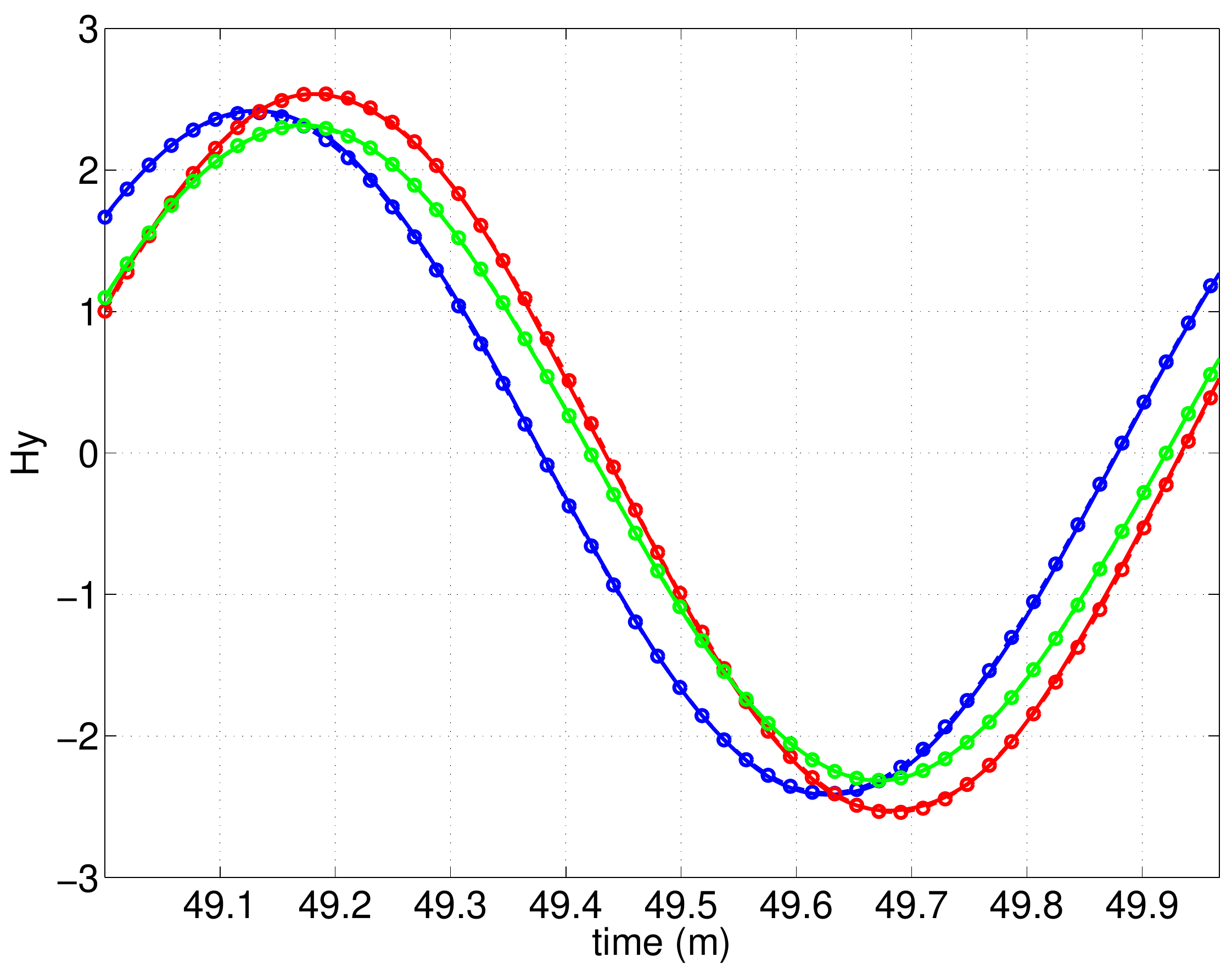}
}
\hspace{0.15cm}
\subfigure[]
{
\includegraphics[width=2.88in]{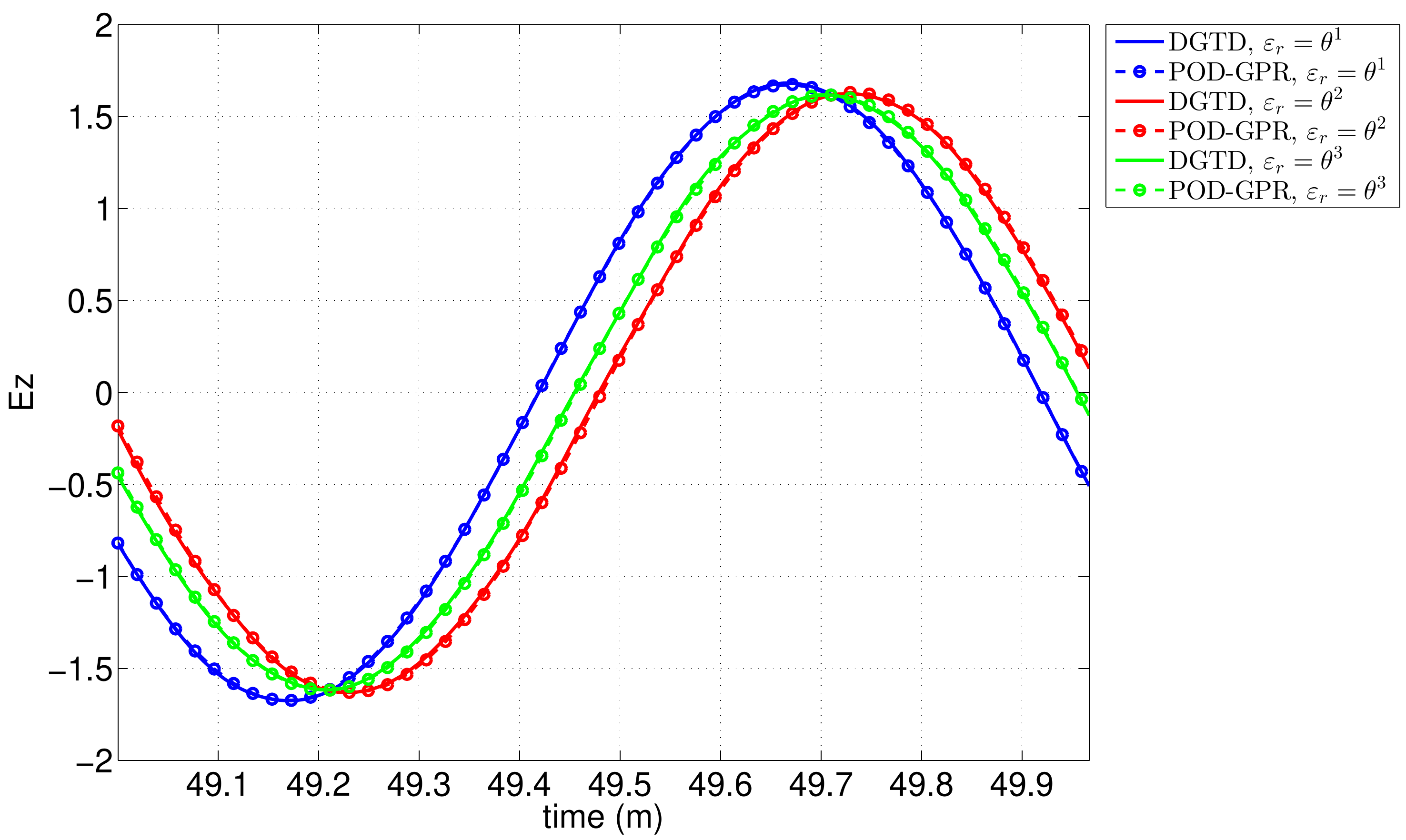}
}
\caption{Scattering of a plane wave by a multi-layer heterogeneous medium: Comparison of the time evolution of the field (a) $H_y$ and (b) $E_z$ at a given point}
\label{fig:12}
\end{figure}

\begin{figure}
\centering
\subfigure
{
\includegraphics[width=2.2in]{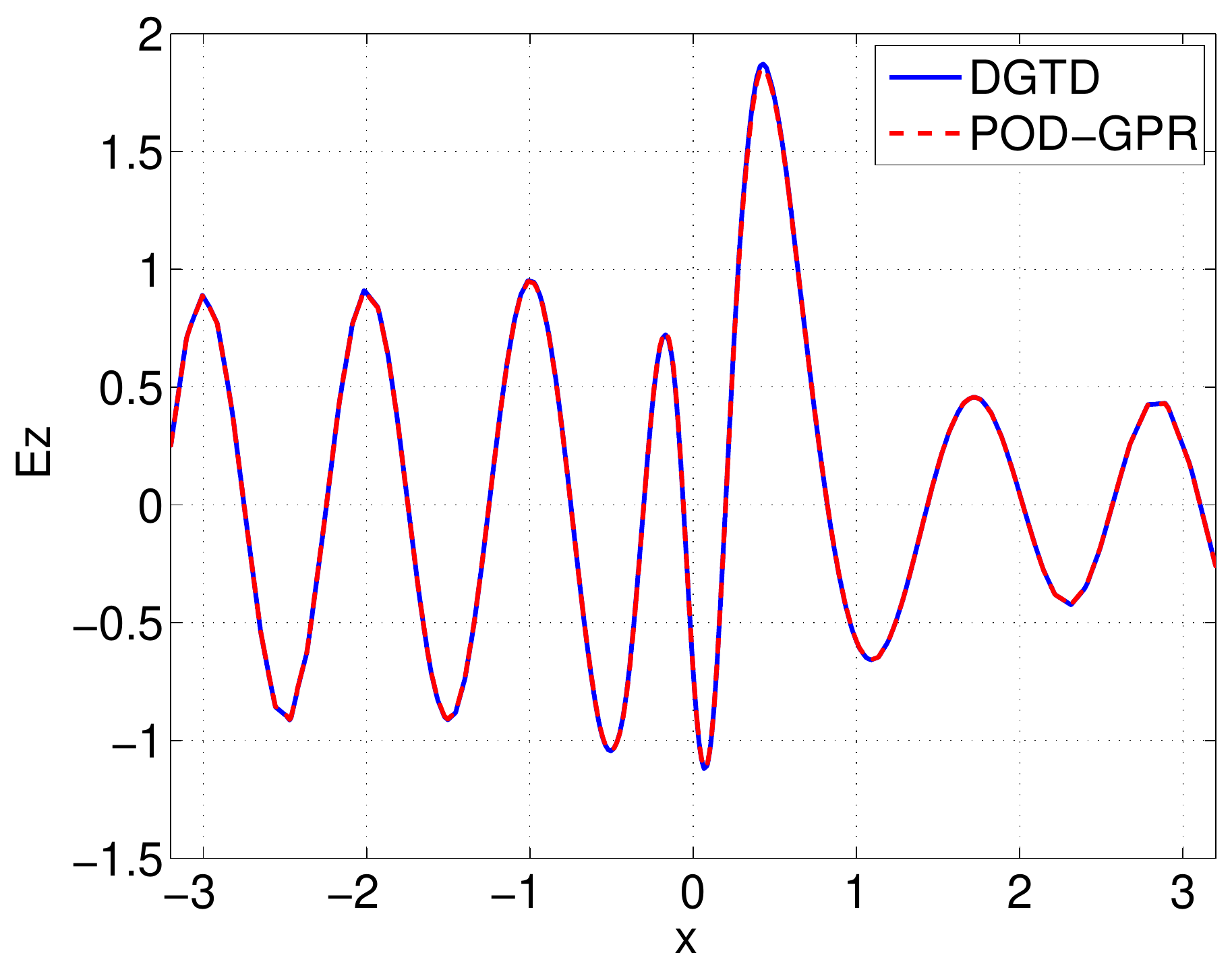}
}
\hspace{0.15cm}
\subfigure
{
\includegraphics[width=2.1in]{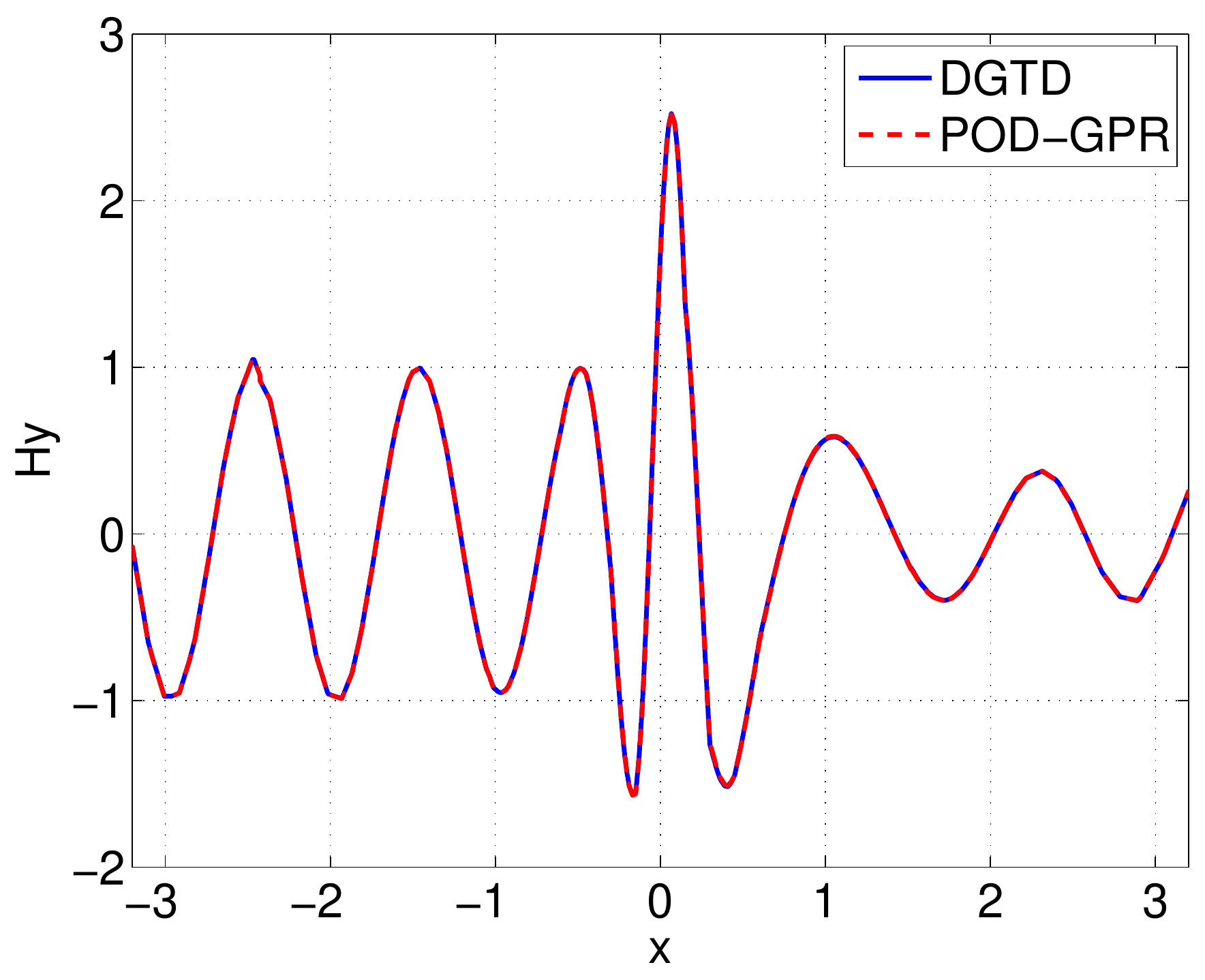}
}\\
\subfigure
{
\includegraphics[width=2.2in]{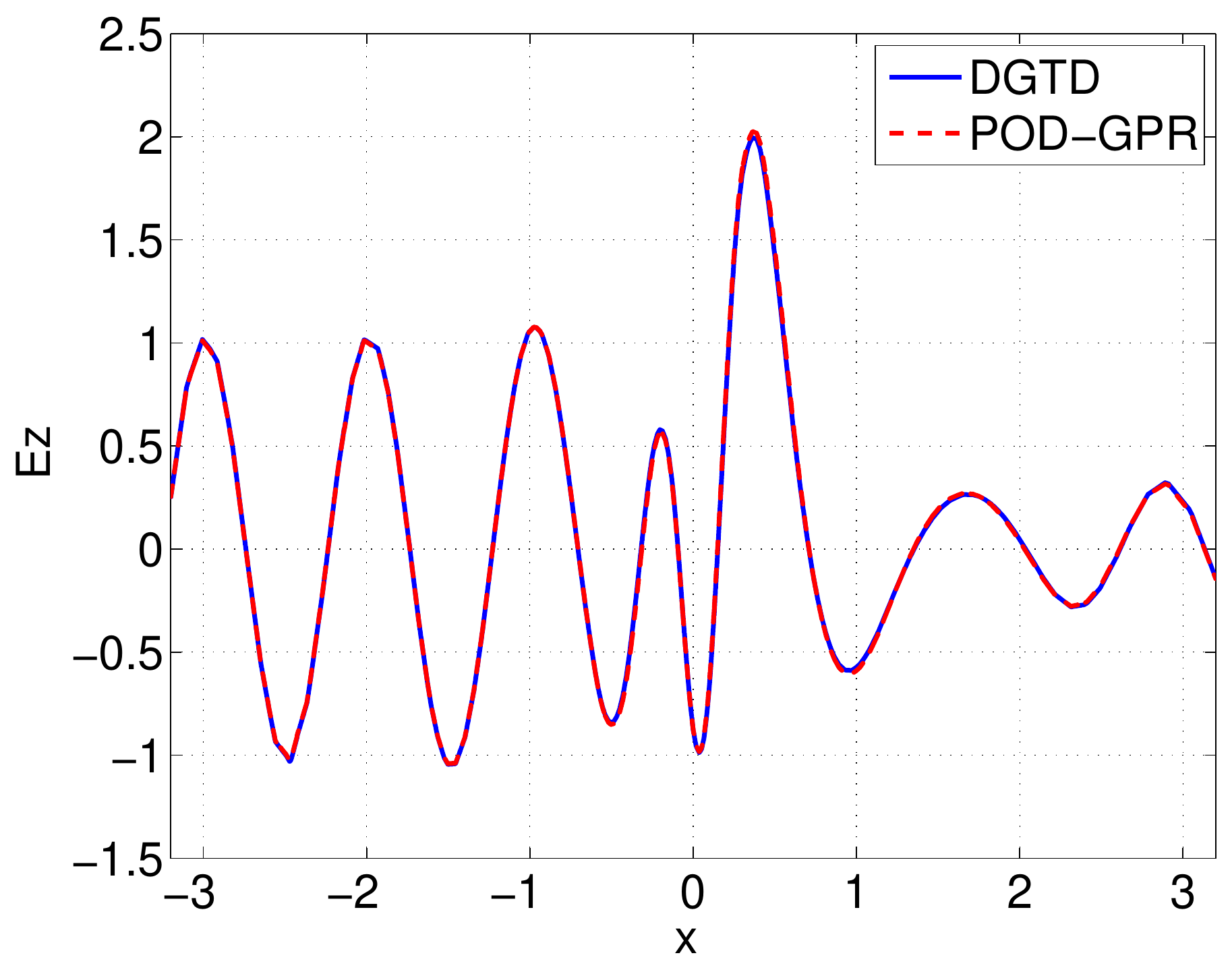}
}
\hspace{0.15cm}
\subfigure
{
\includegraphics[width=2.1in]{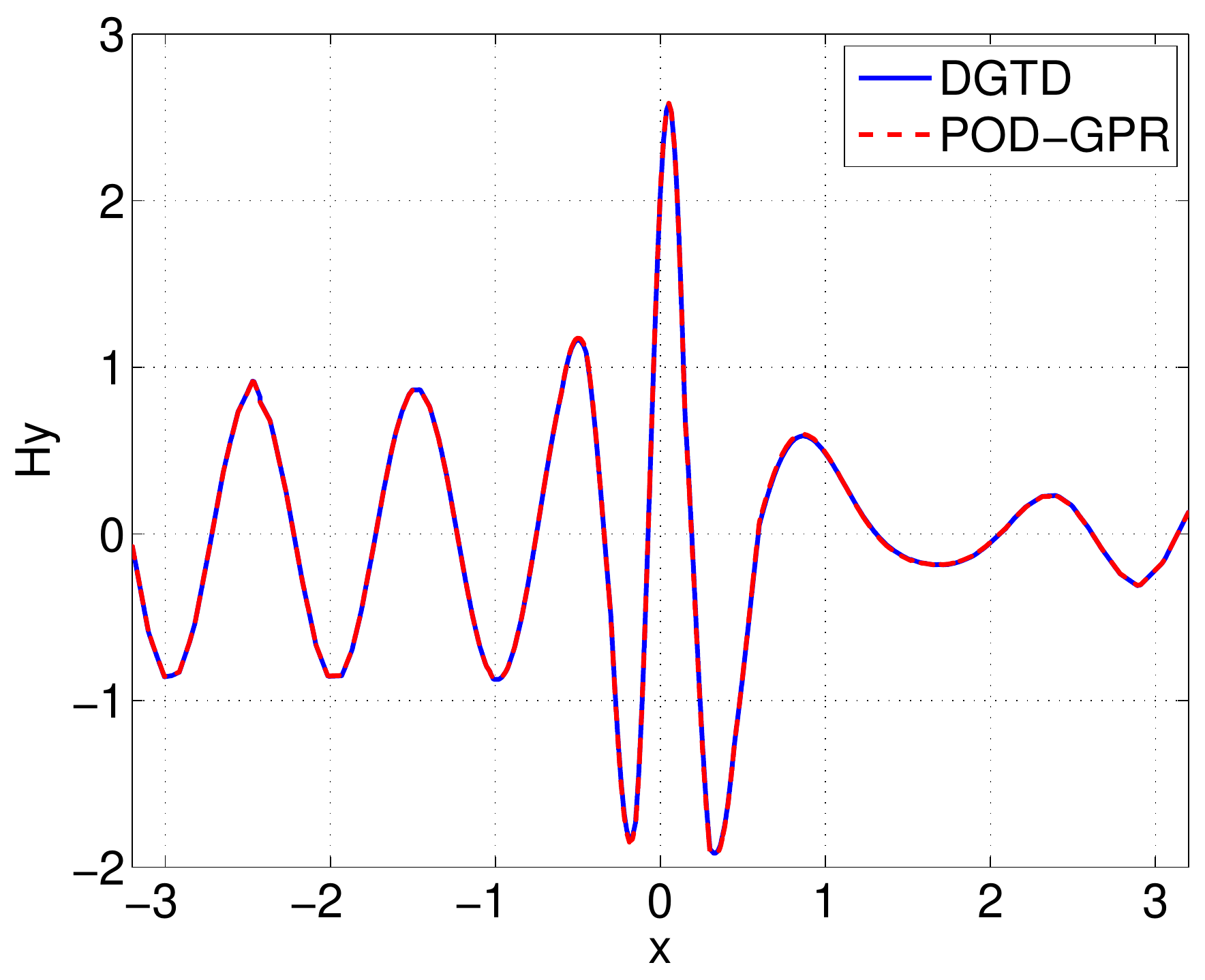}
}\\
\subfigure
{
\includegraphics[width=2.2in]{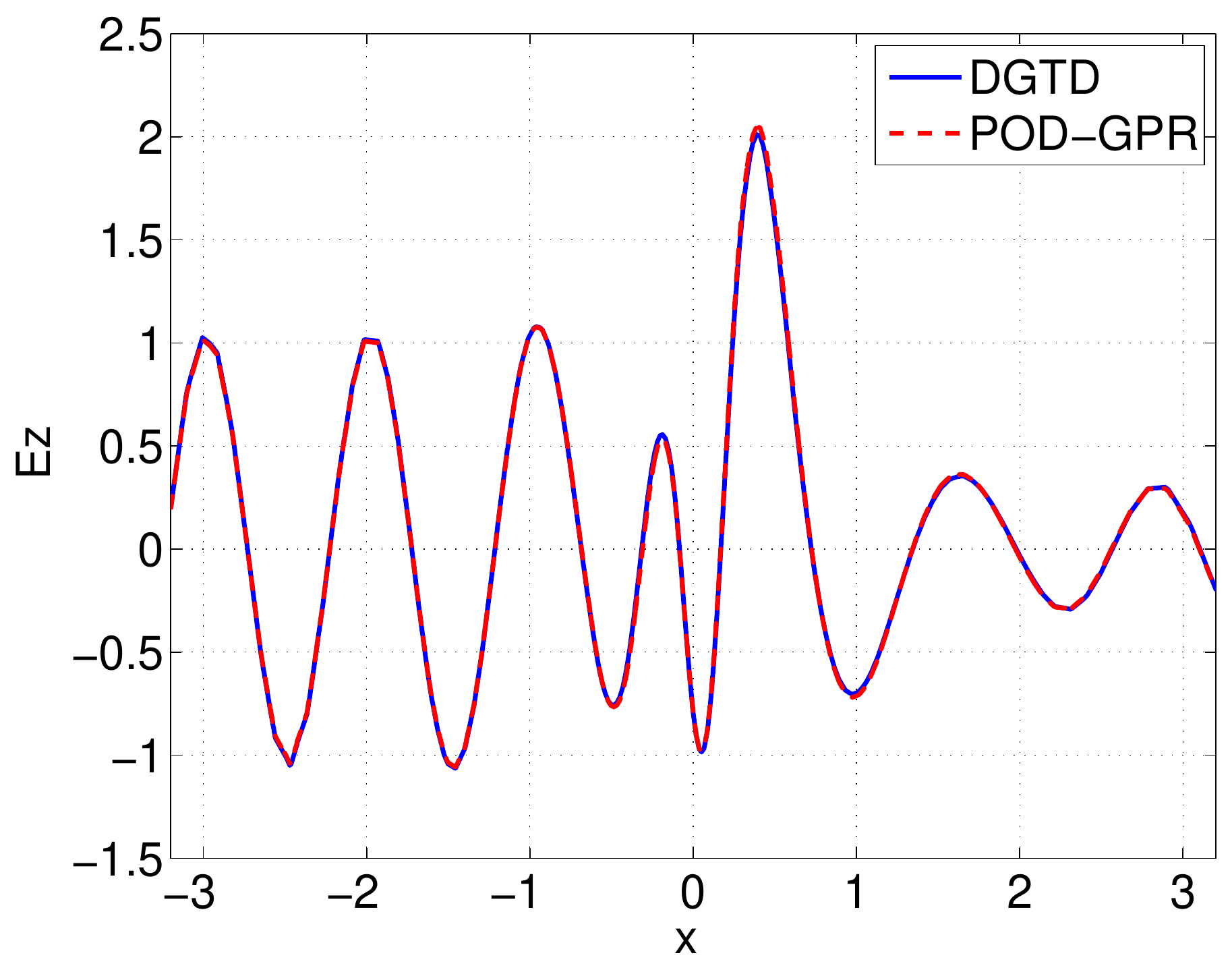}
}
\hspace{0.15cm}
\subfigure
{
\includegraphics[width=2.1in]{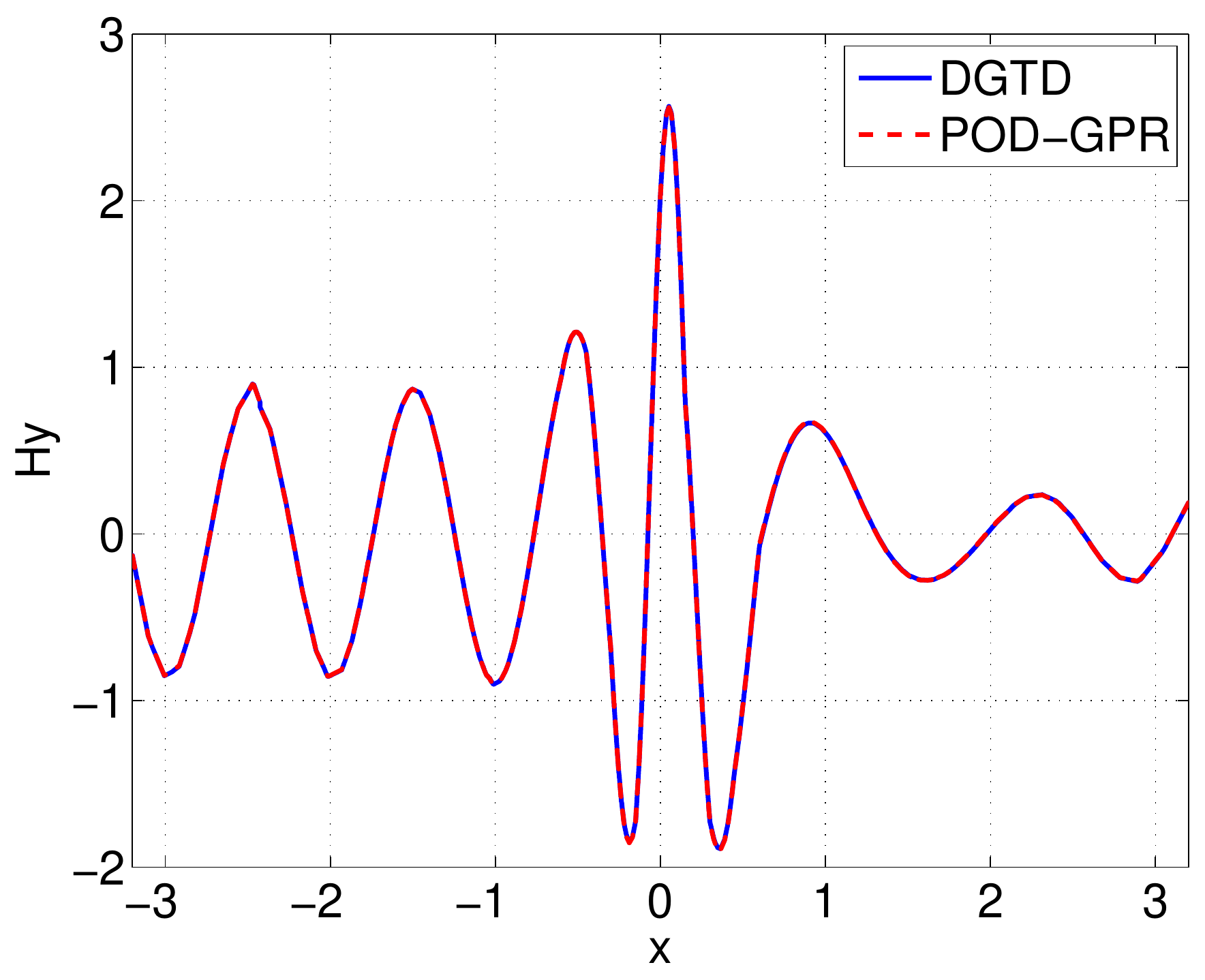}
}\\
\caption{Scattering of a plane wave by a multi-layer heterogeneous medium: Comparison of the 1D x-wise distribution along $y=0$ of the real part of $E_z$ (left) and $H_y$ (right) of three test points: $\theta^1$ (1st row), $\theta^2$ (2nd row) and $\theta^3$ (3rd row).}
\label{fig:13} 
\end{figure}

\begin{figure}
\centering
\subfigure
{
\includegraphics[width=2.2in]{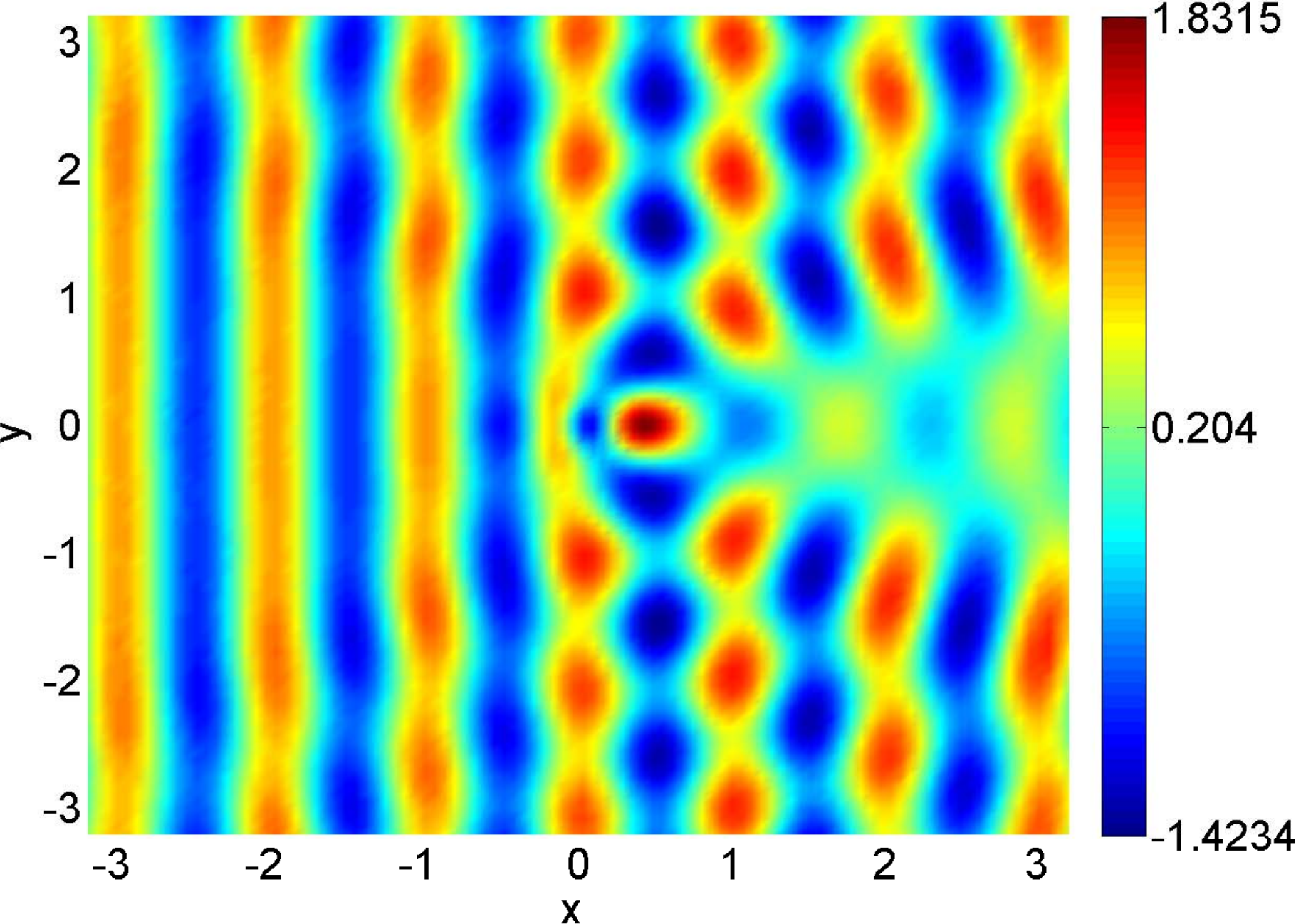}
}
\hspace{0.15cm}
\subfigure
{
\includegraphics[width=2.2in]{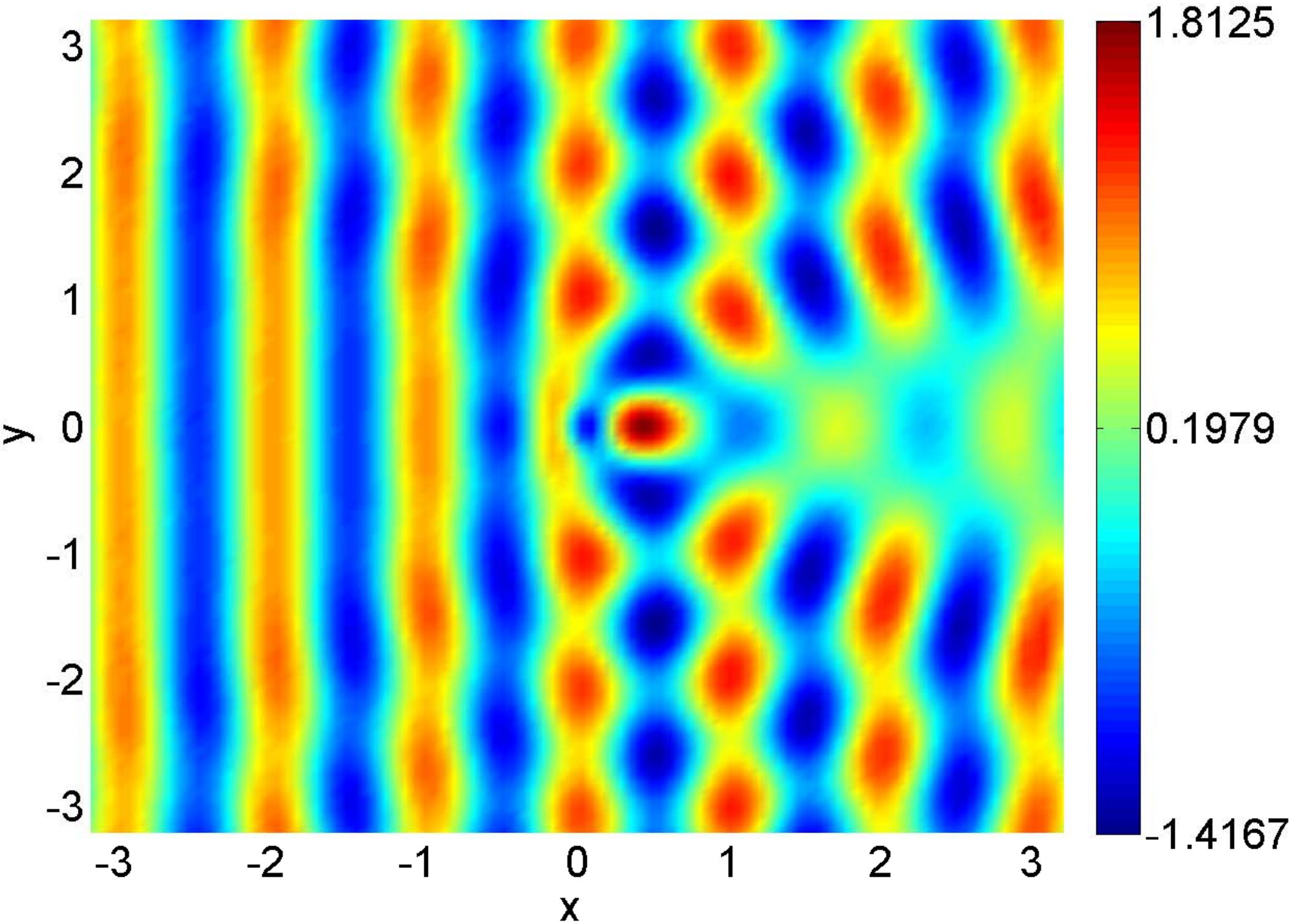}
}\\
\subfigure
{
\includegraphics[width=2.2in]{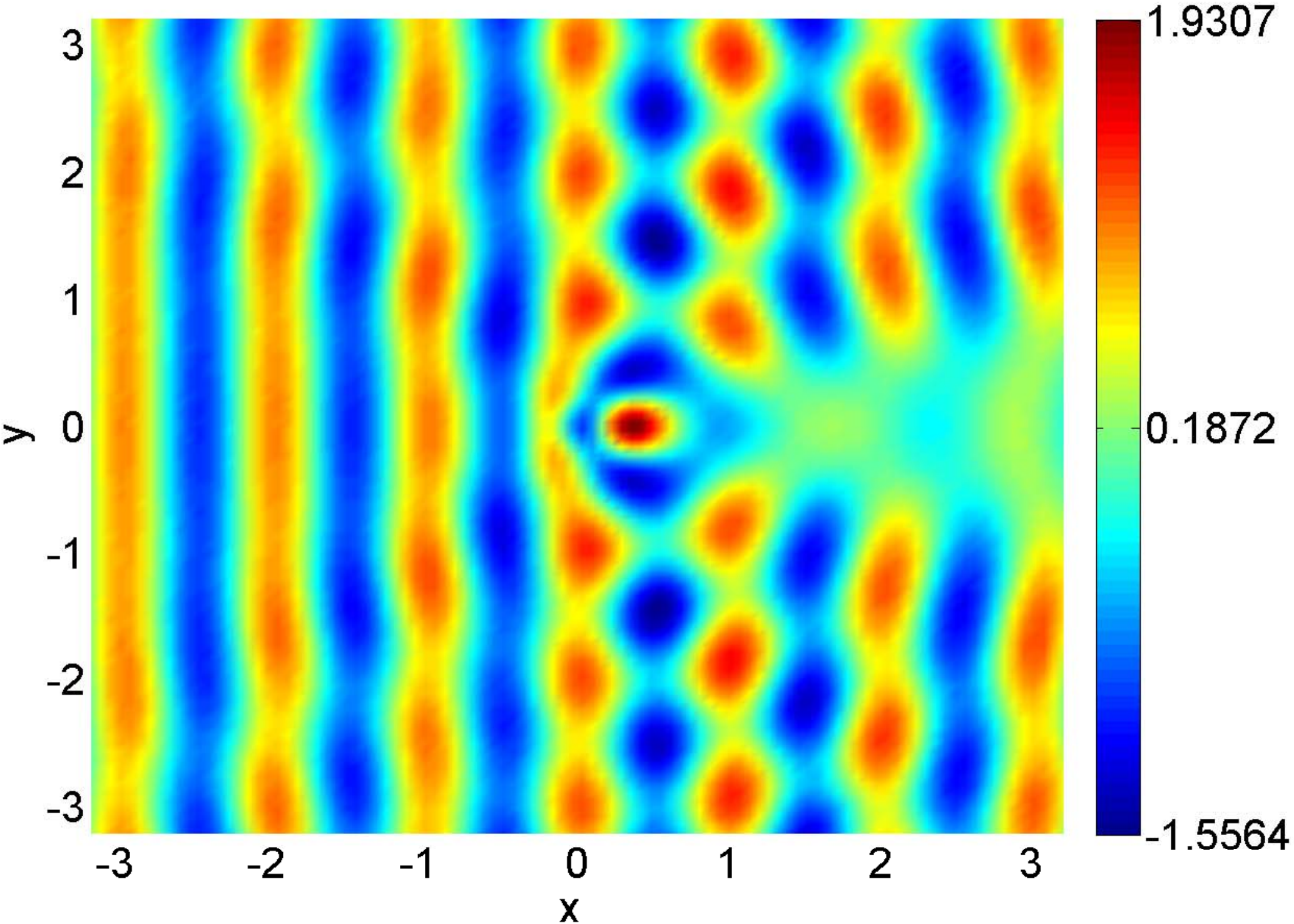}
}
\hspace{0.15cm}
\subfigure
{
\includegraphics[width=2.2in]{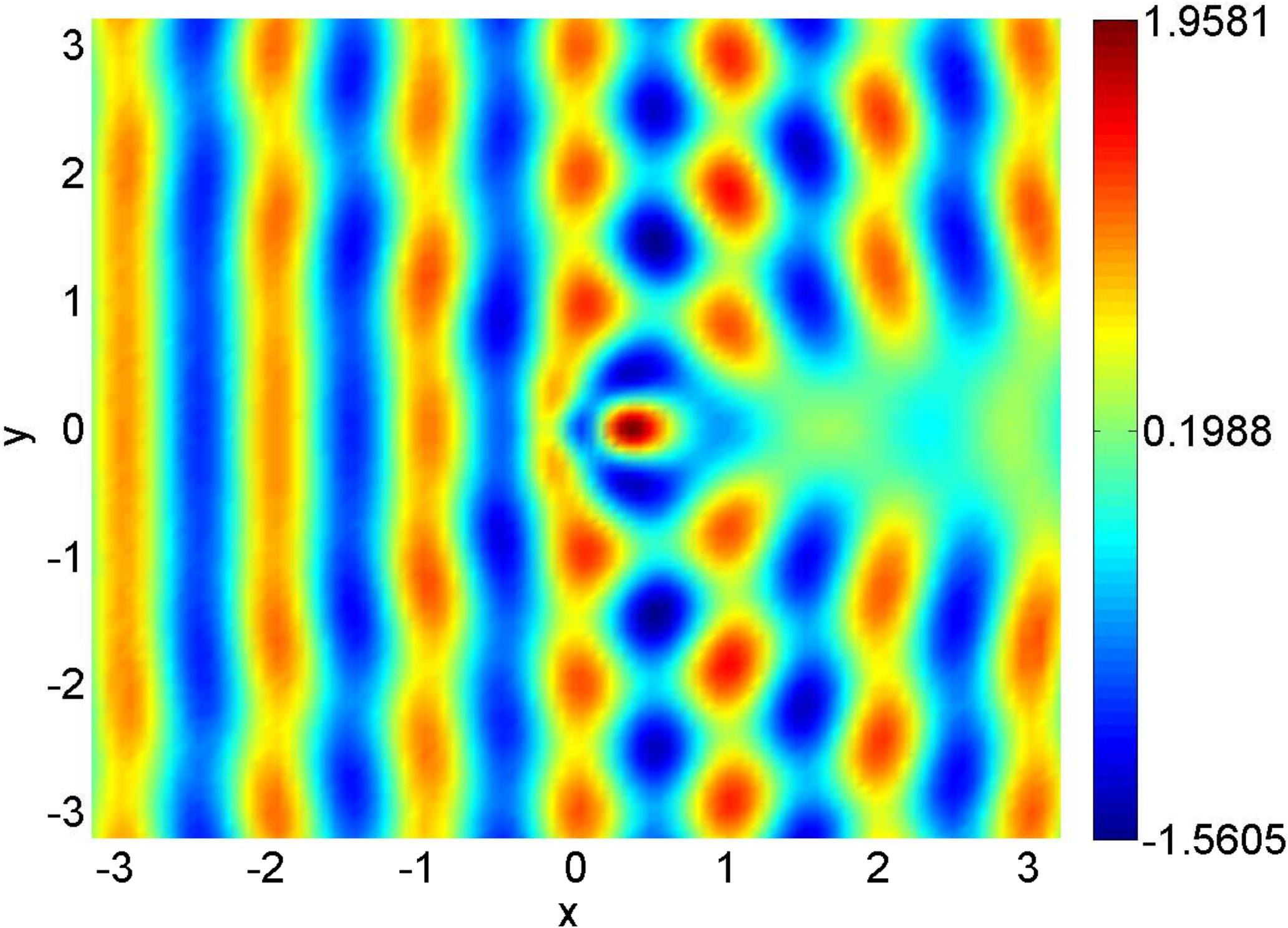}
}\\
\subfigure
{
\includegraphics[width=2.2in]{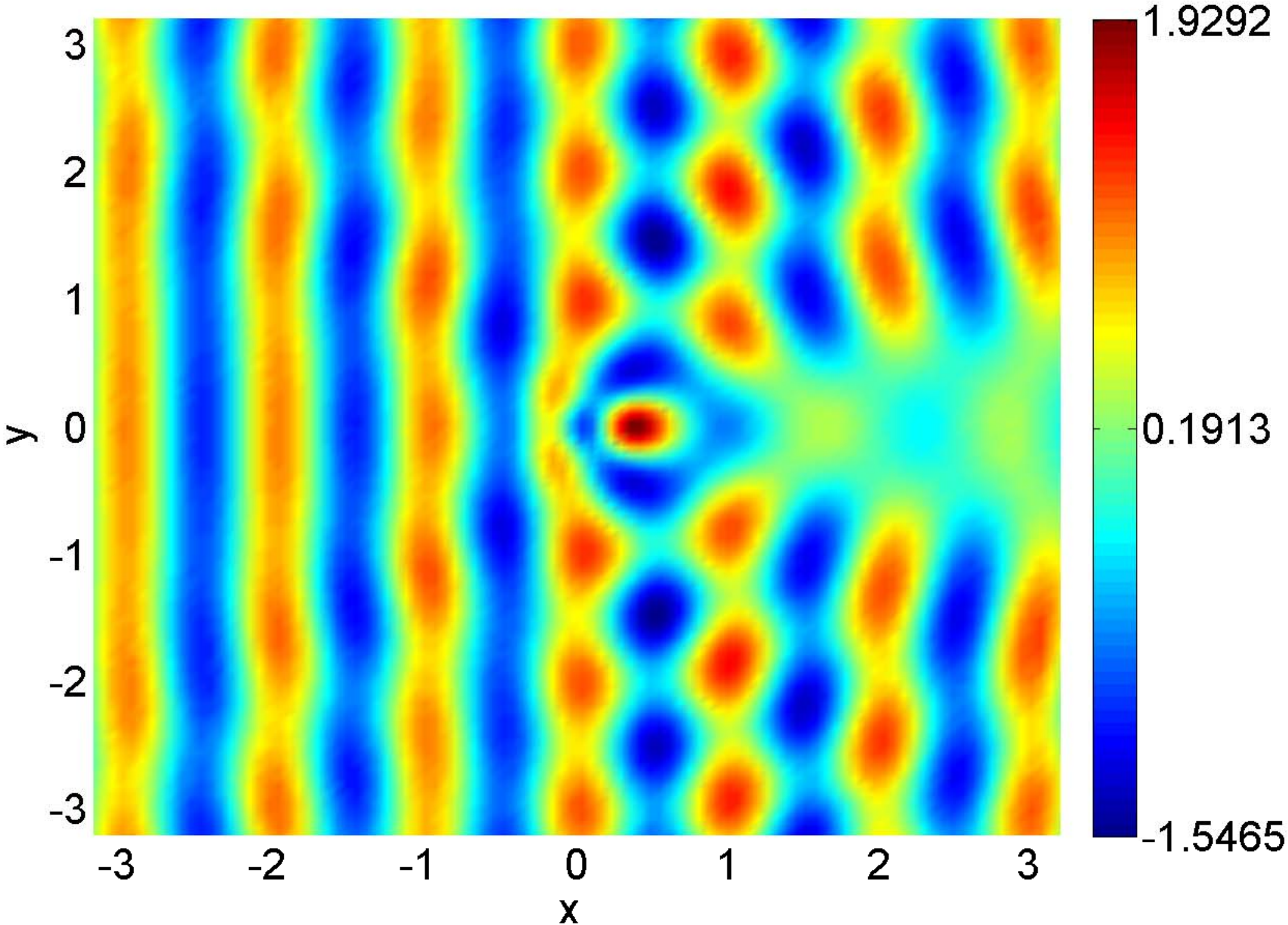}
}
\hspace{0.15cm}
\subfigure
{
\includegraphics[width=2.2in]{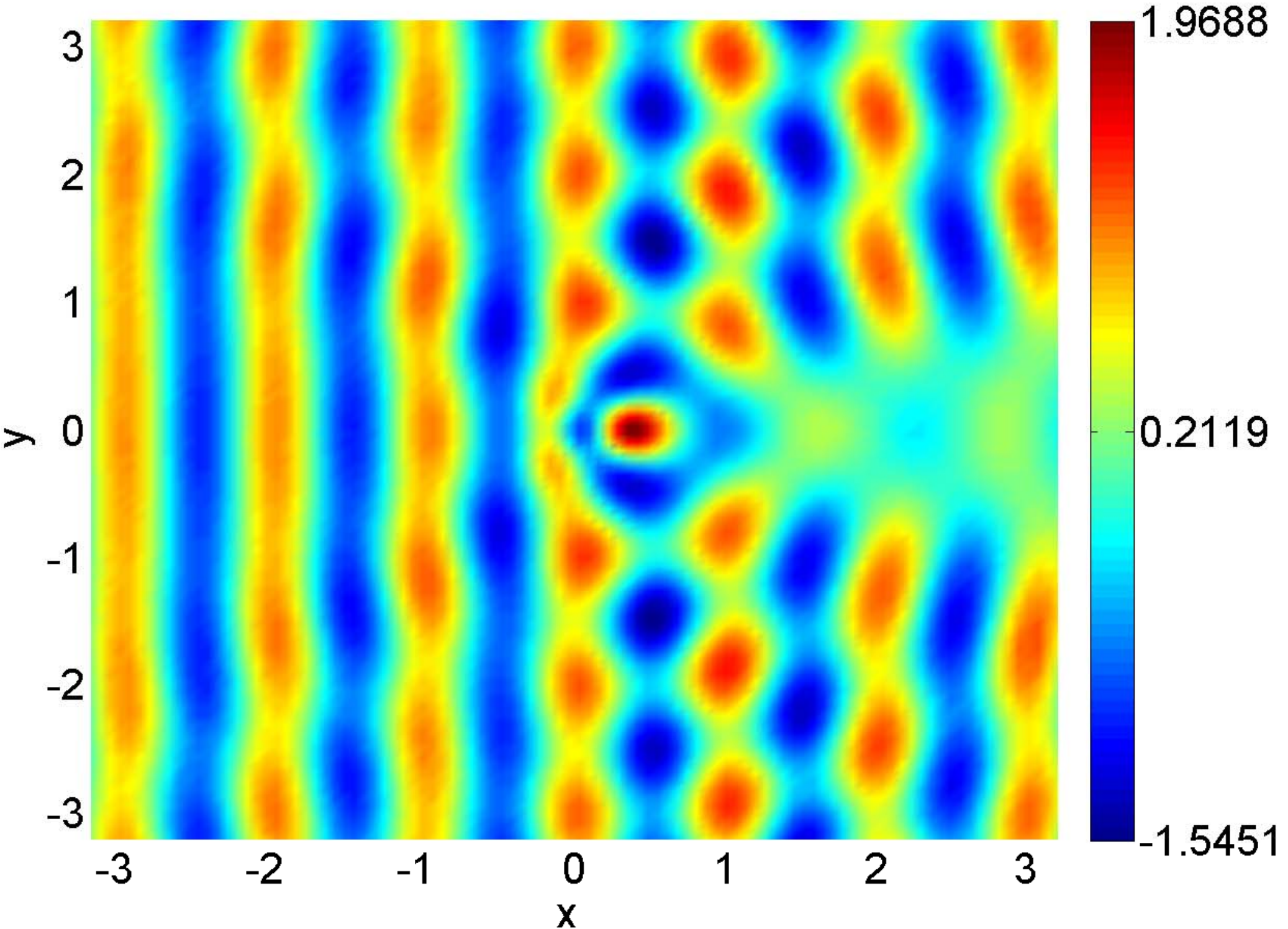}
}\\
\caption{Scattering of a plane wave by a multi-layer heterogeneous medium: Comparison of the 2D distribution of the real part of $E_z$ between DGTD (left) and POD-GPR (right) of of three test points: $\theta^1$ (1st row), $\theta^2$ (2nd row) and $\theta^3$ (3rd row).}
\label{fig:14} 
\end{figure}

\begin{figure}
\centering
\subfigure
{
\includegraphics[width=2.2in]{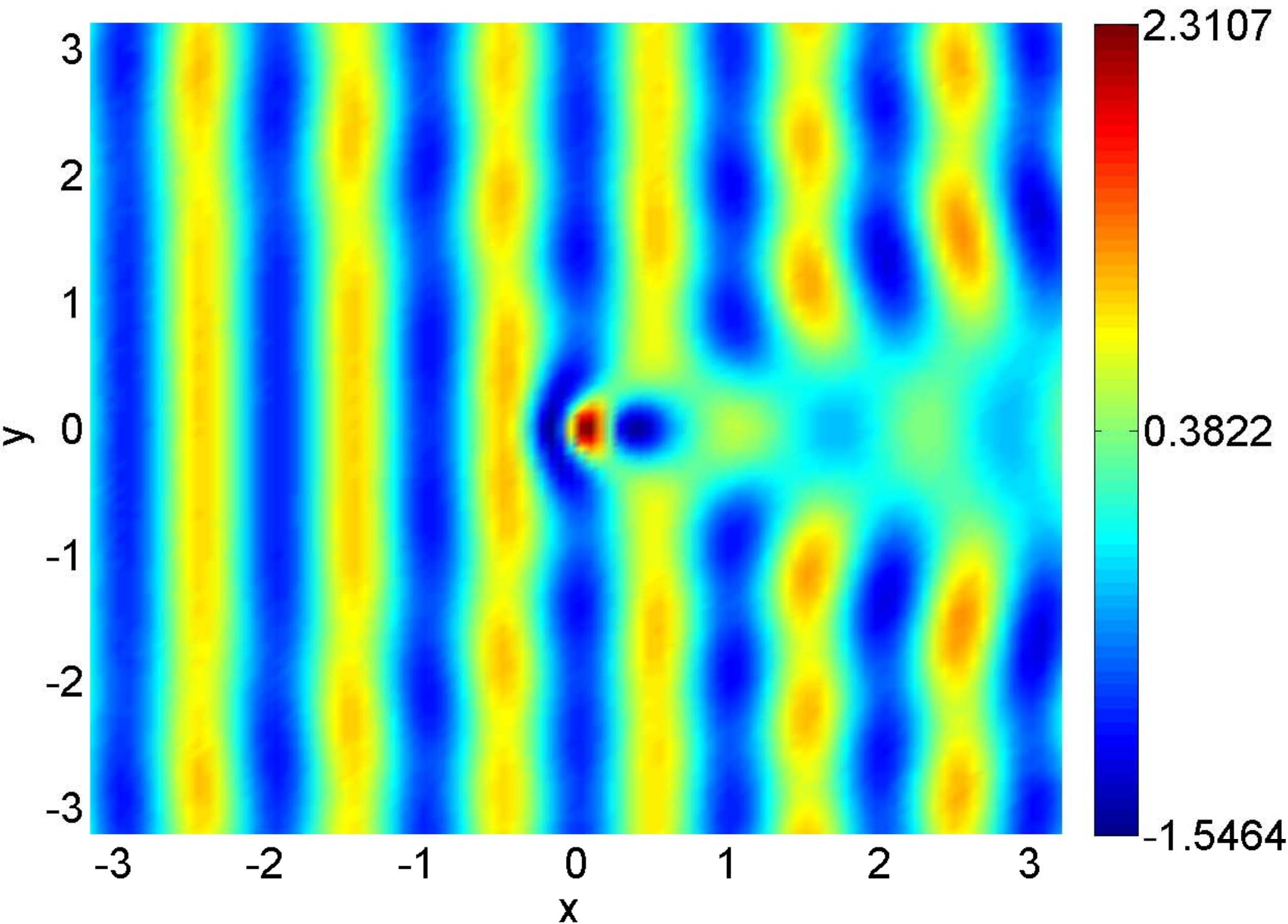}
}
\hspace{0.15cm}
\subfigure
{
\includegraphics[width=2.2in]{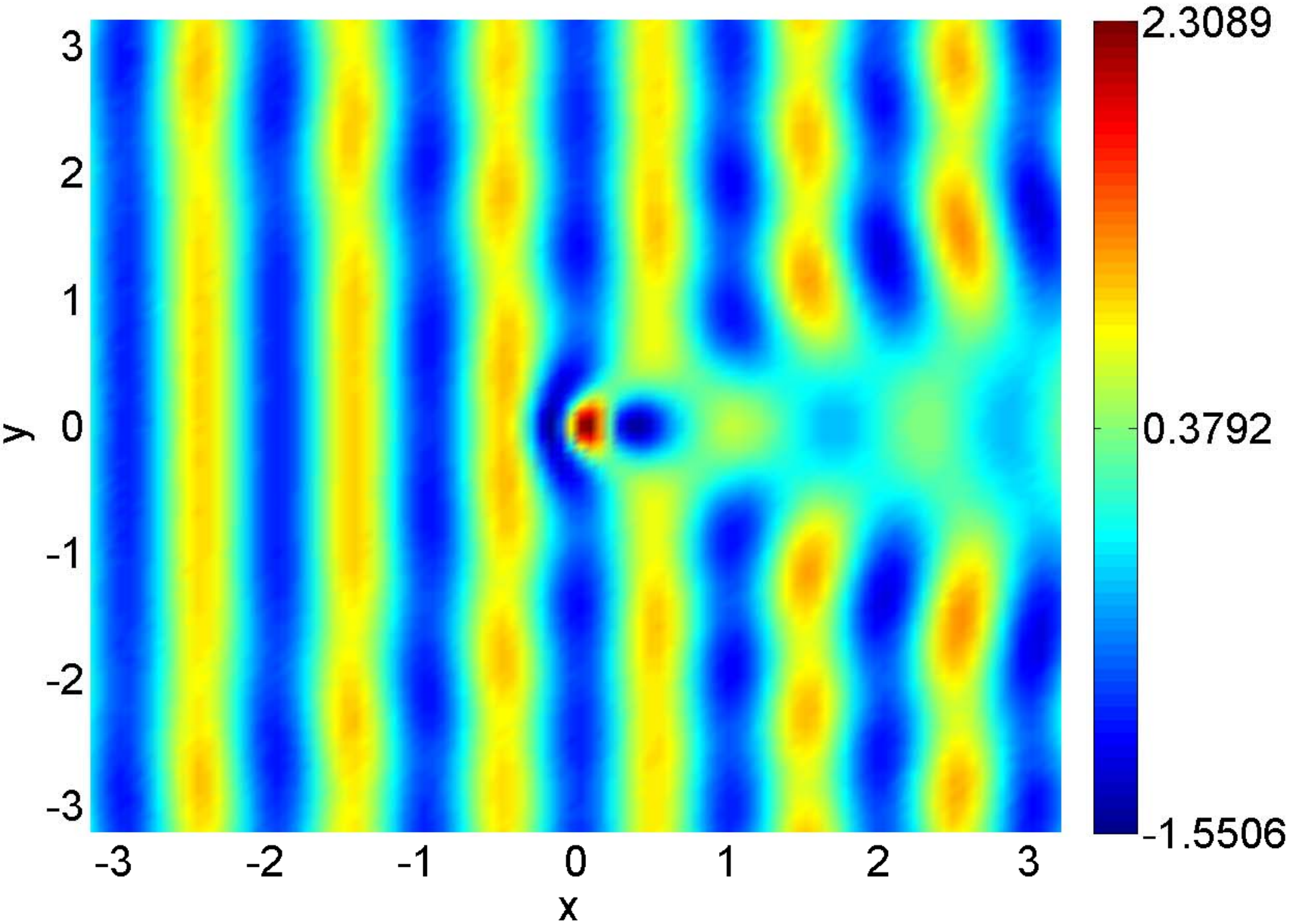}
}\\
\subfigure
{
\includegraphics[width=2.2in]{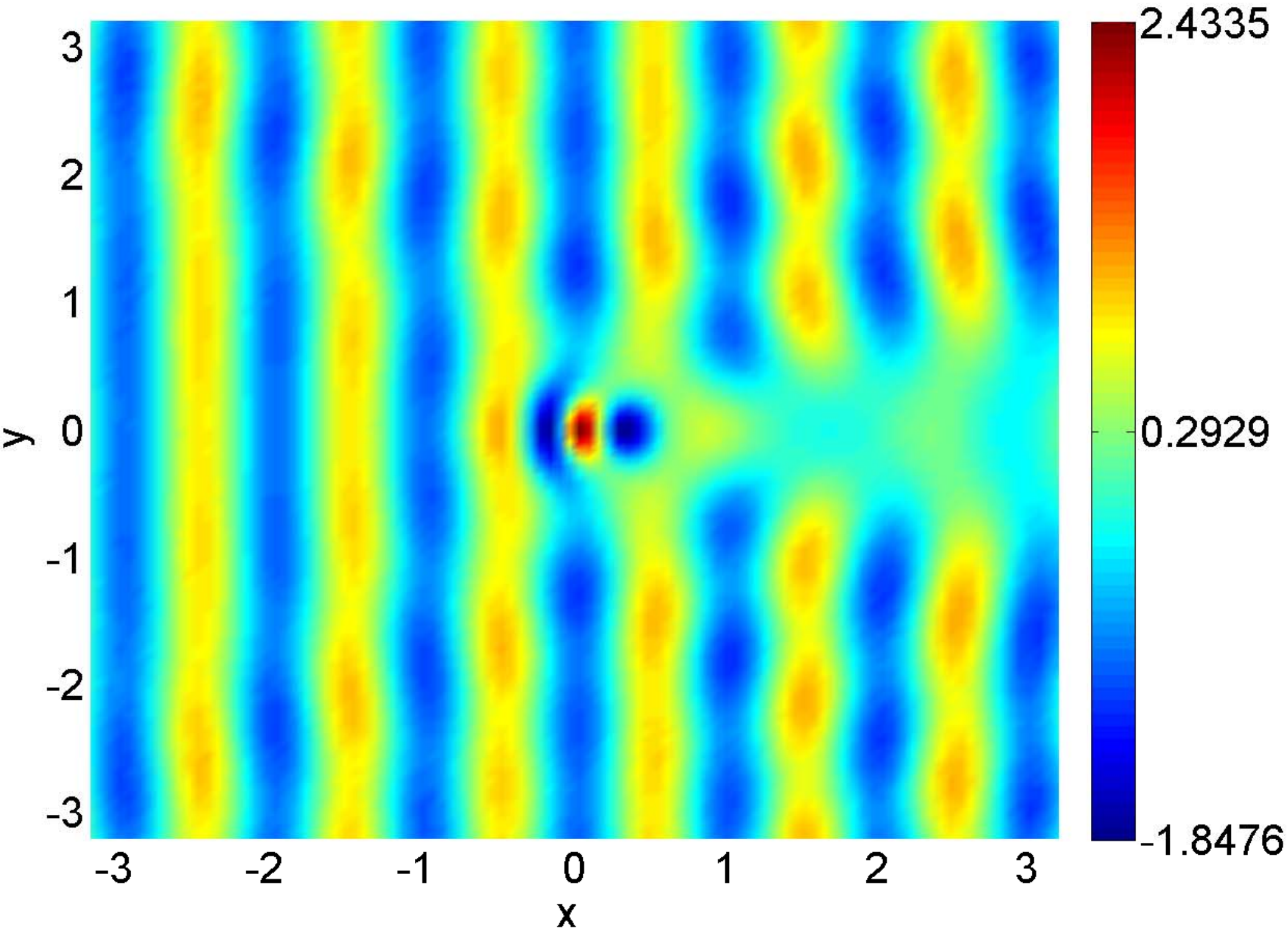}
}
\hspace{0.15cm}
\subfigure
{
\includegraphics[width=2.2in]{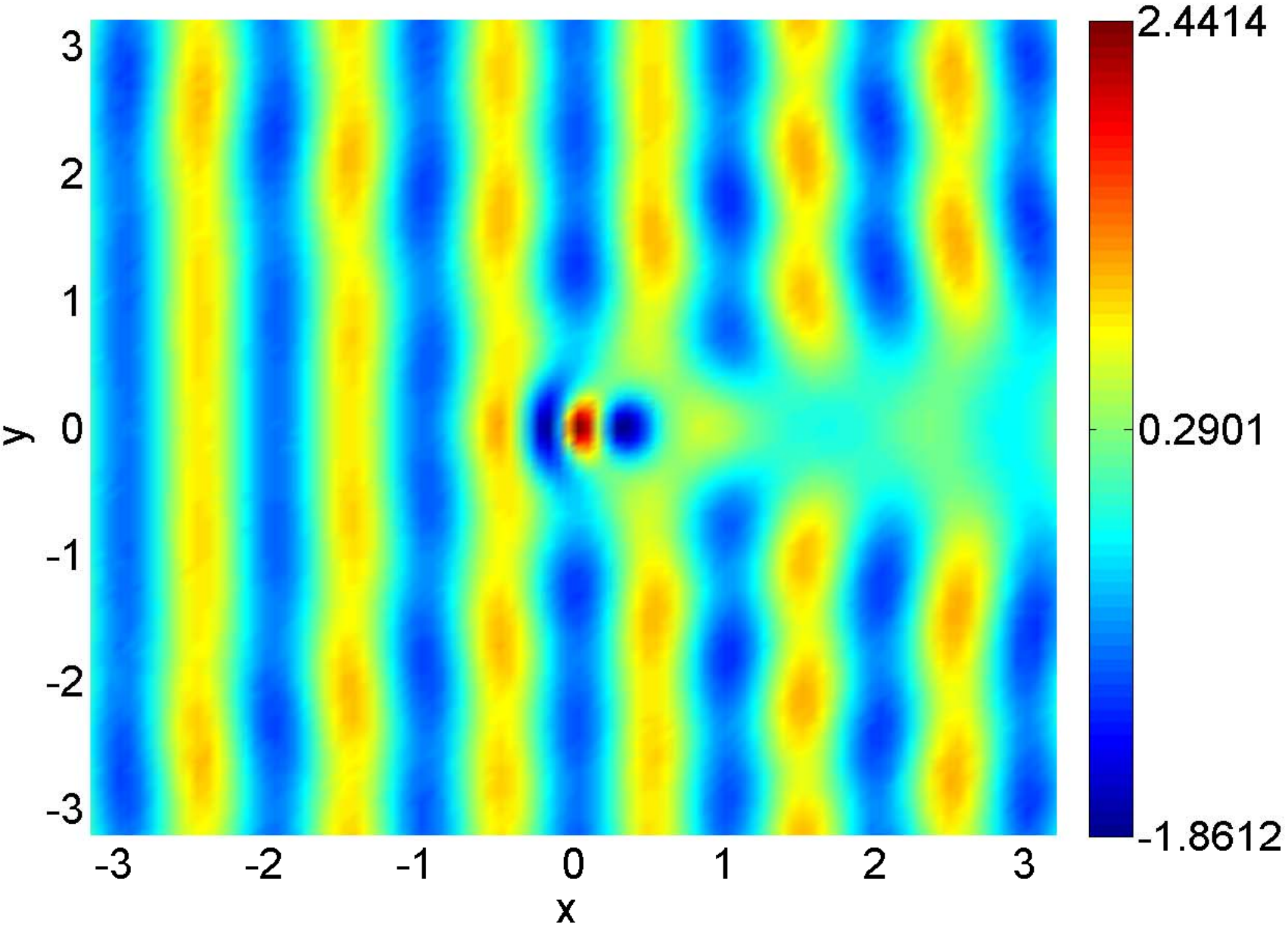}
}\\
\subfigure
{
\includegraphics[width=2.2in]{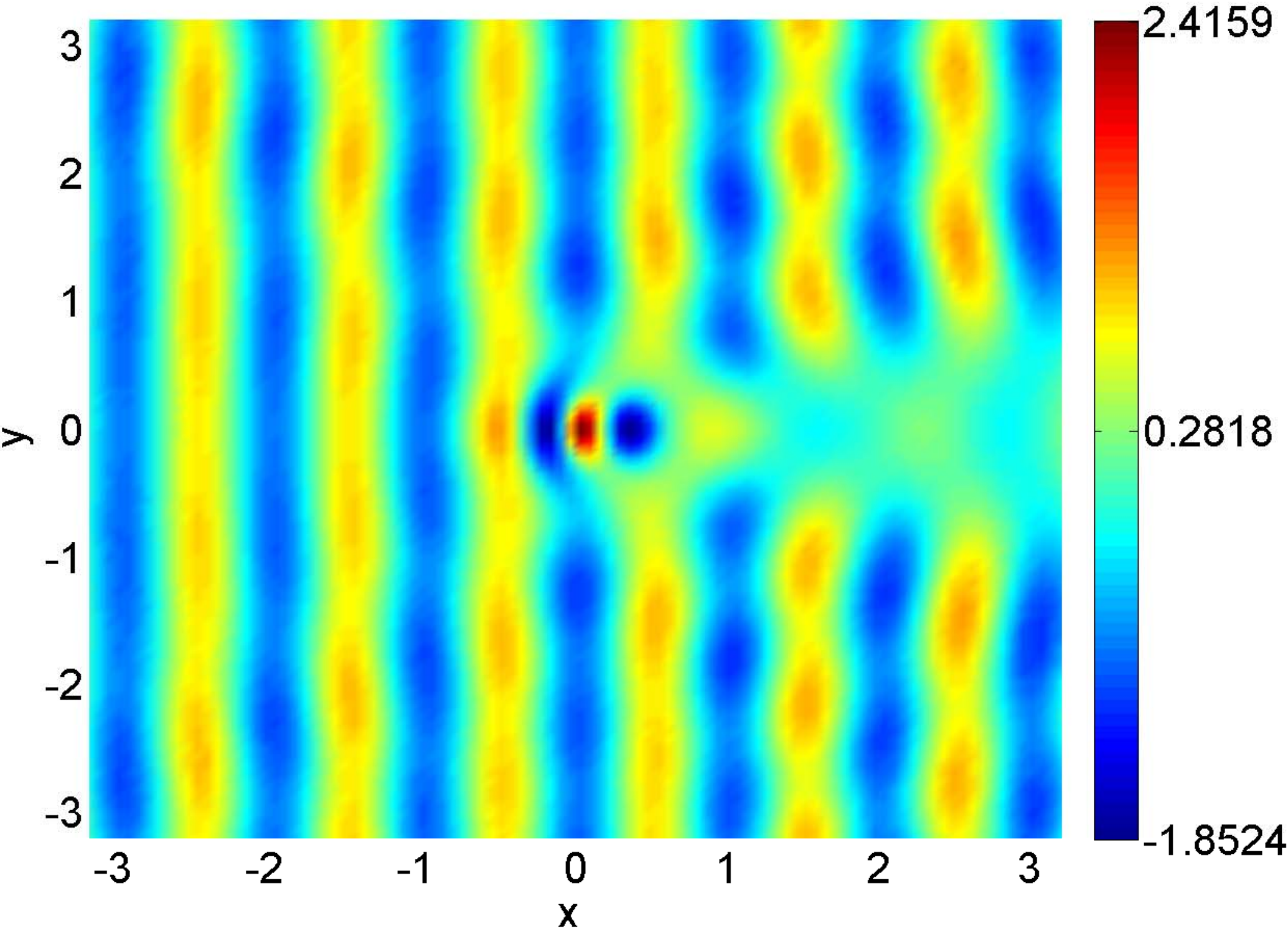}
}
\hspace{0.15cm}
\subfigure
{
\includegraphics[width=2.2in]{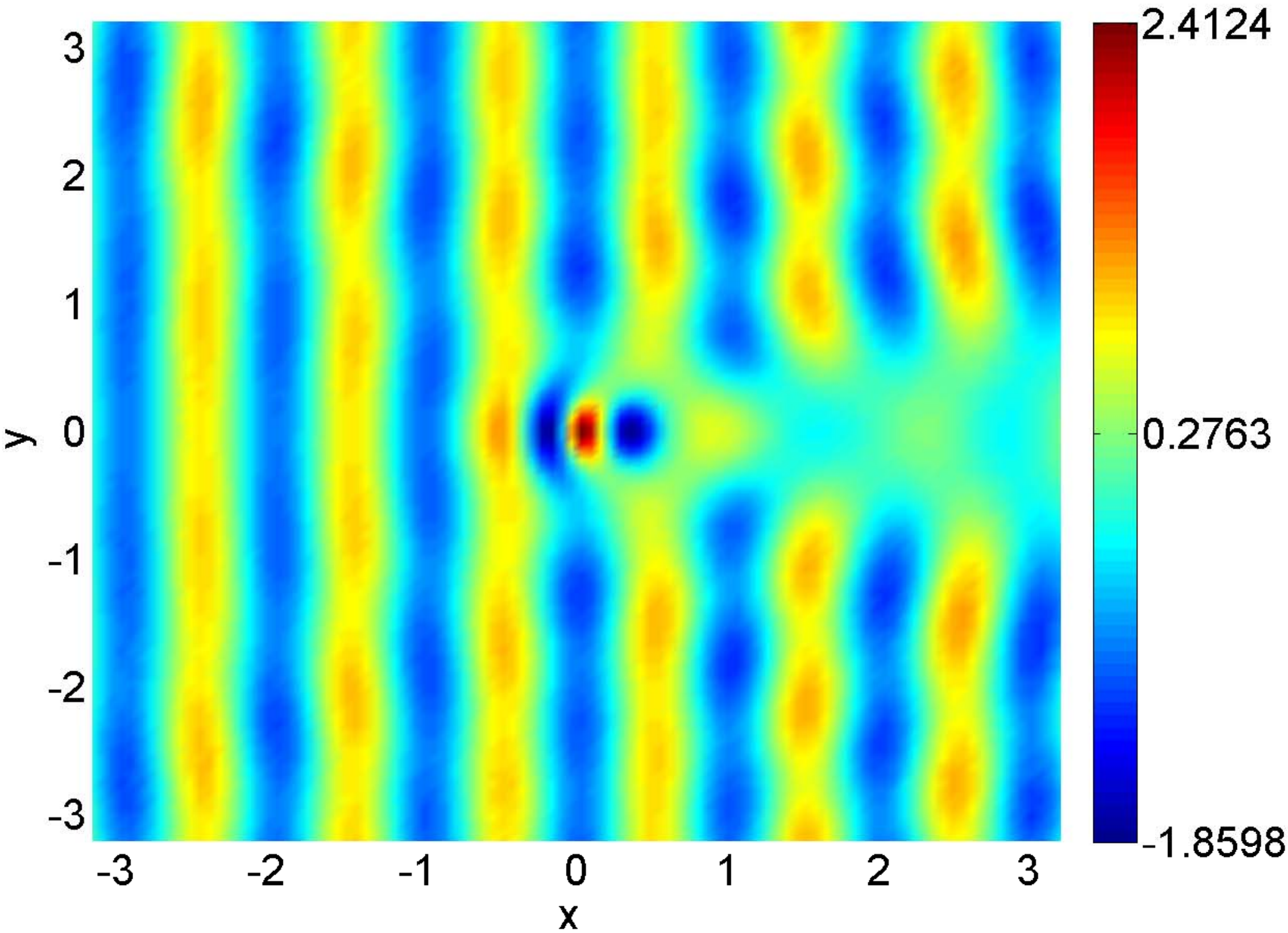}
}\\
\caption{Scattering of a plane wave by a multi-layer heterogeneous medium: Comparison of the 2D distribution of the real part of $H_y$ between DGTD (left) and POD-GPR (right) of of three test points: $\theta^1$ (1st row), $\theta^2$ (2nd row) and $\theta^3$ (3rd row).}
\label{fig:15}
\end{figure}

Following up with the offline phase, online tests are implemented for three non-trained $\theta$, $\theta^1=[5.15,3.375,2.125,1.375]$, $\theta^2=[5.45,3.625,2.375,1.625]$, $\theta^3=[5.215,3.325,2.455,1.655]$, whose POD basis expansion coefficients are sought by GPR models. An evidence of the efficacy of the POD-GPR is given in Fig.\ref{fig:12}, 
which reports the comparison between the time evolution of full-order solutions at a given point and their reduced ones. Moreover, to more intuitively see the simulated electromagnetic field, in Fig.\ref{fig:13} displays the1D x-wise distributions of the real part of $E_z$ and $H_y$ in the Fourier domain during the last period of wave oscillation, and Fig.\ref{fig:14} and Fig.\ref{fig:15} are their 2D contour lines distributions, demonstrating a fine matching between DGTD solutions and reduced-order solutions.

\begin{figure}
\centering
\subfigure
{
\includegraphics[width=2.2in]{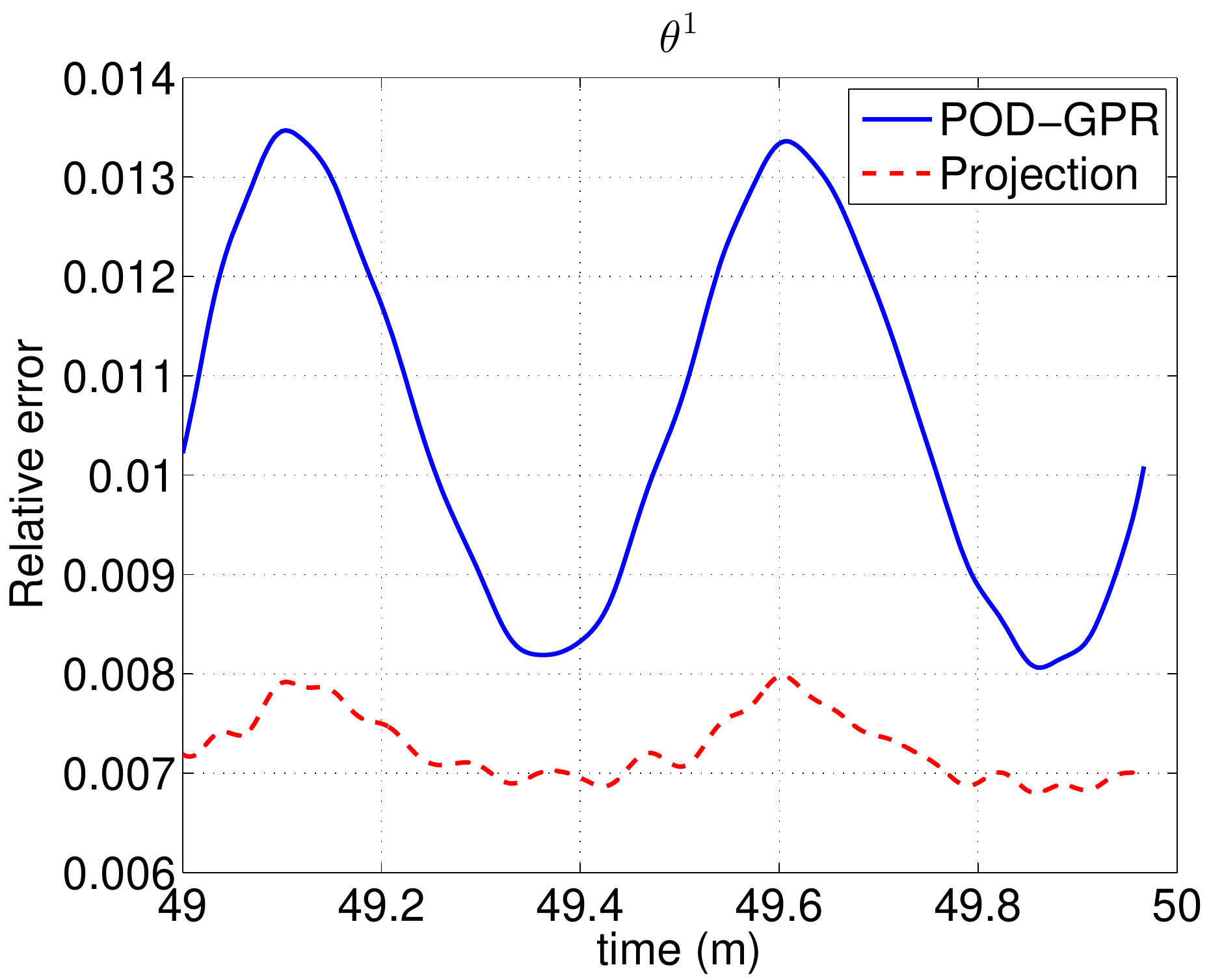}
}
\hspace{0.15cm}
\subfigure
{
\includegraphics[width=2.1in]{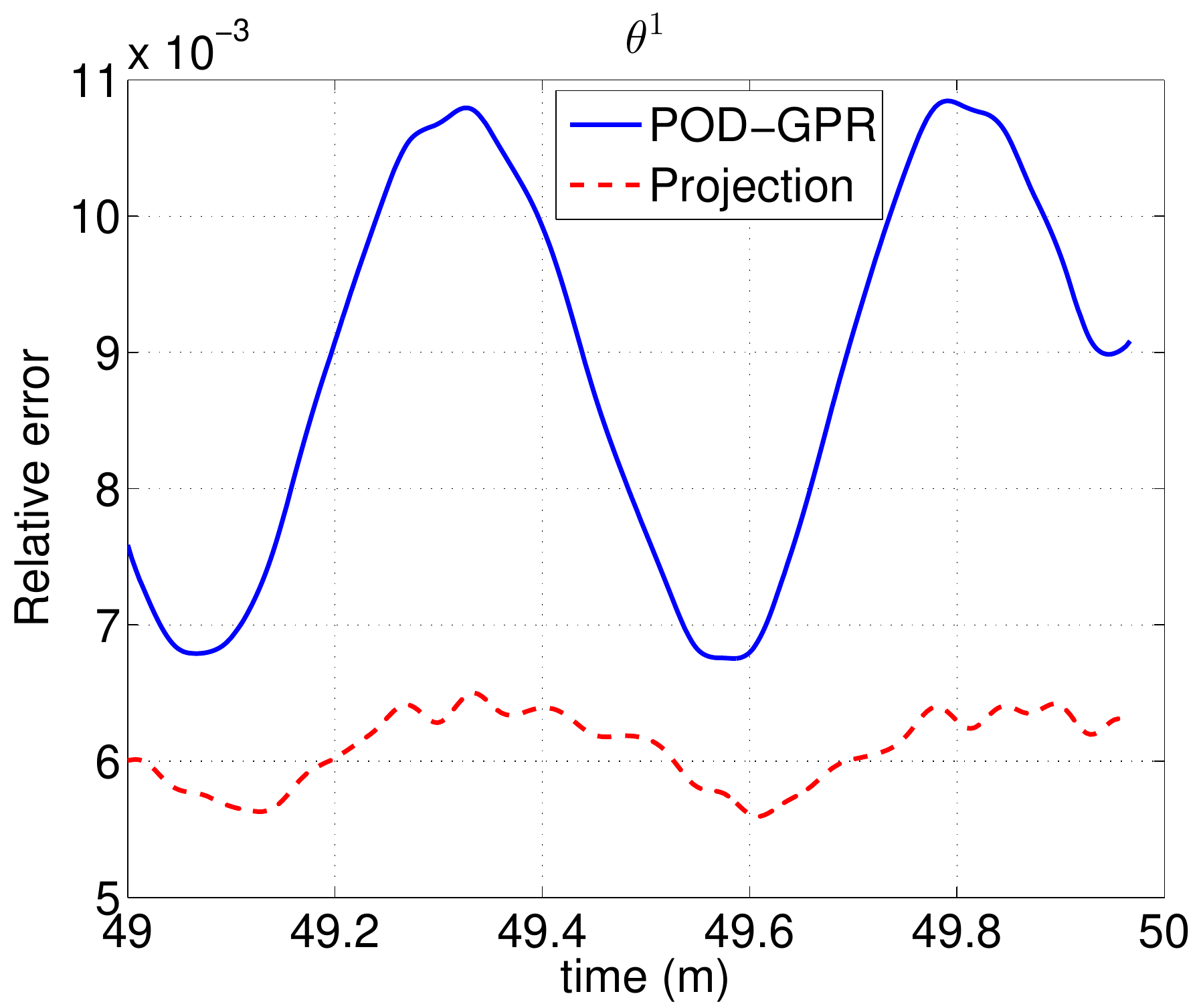}
}\\
\subfigure
{
\includegraphics[width=2.2in]{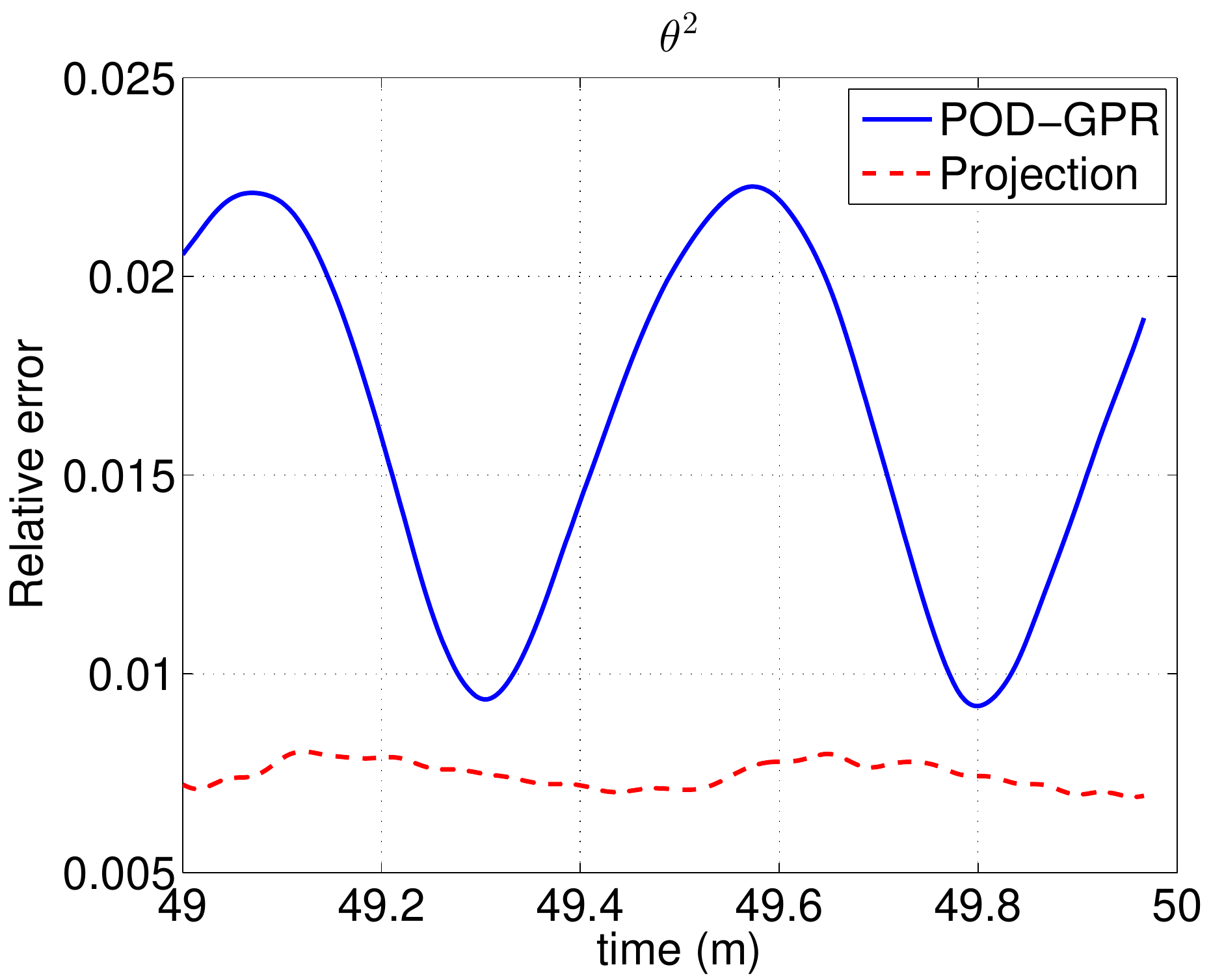}
}
\hspace{0.15cm}
\subfigure
{
\includegraphics[width=2.2in]{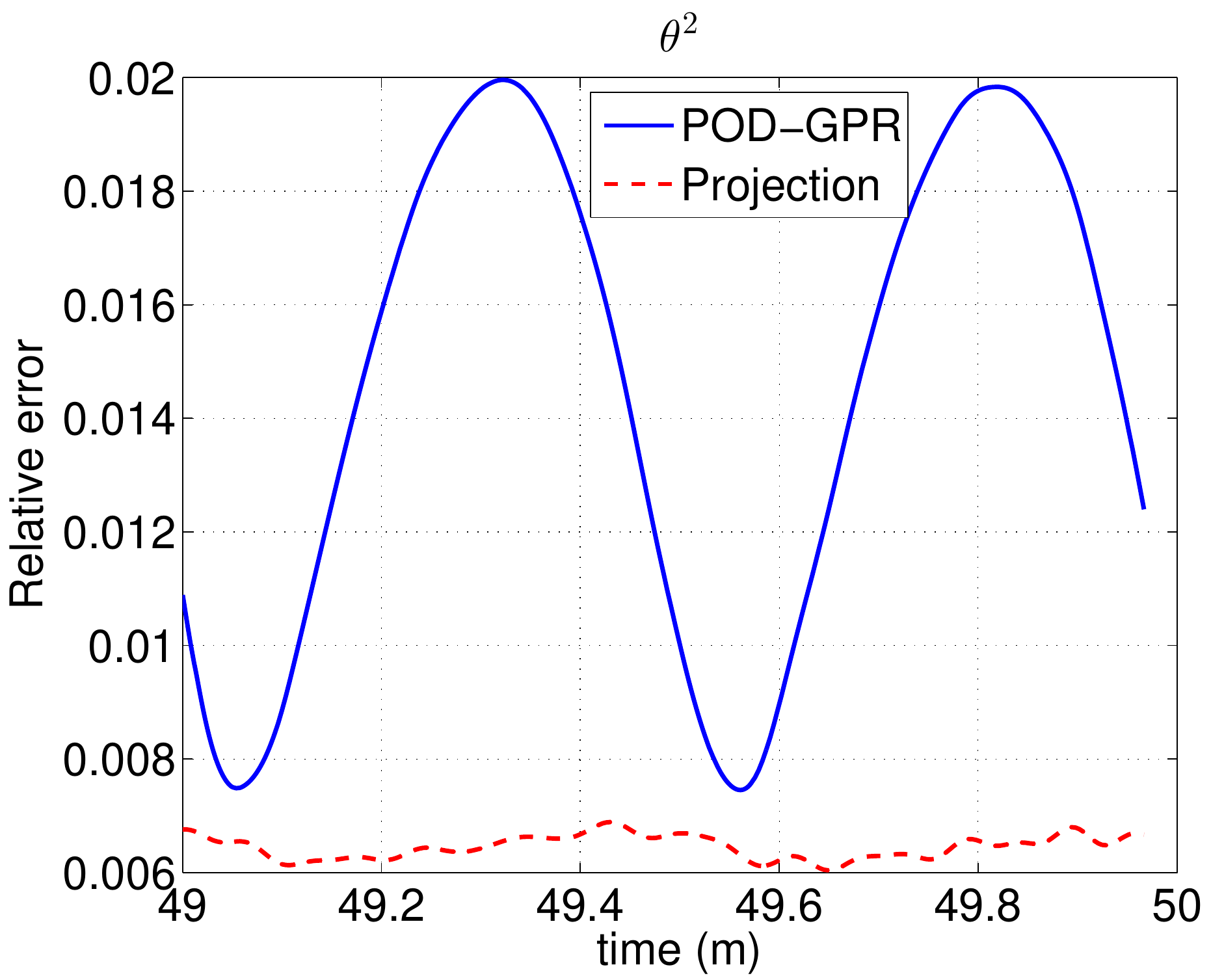}
}\\
\subfigure
{
\includegraphics[width=2.3in]{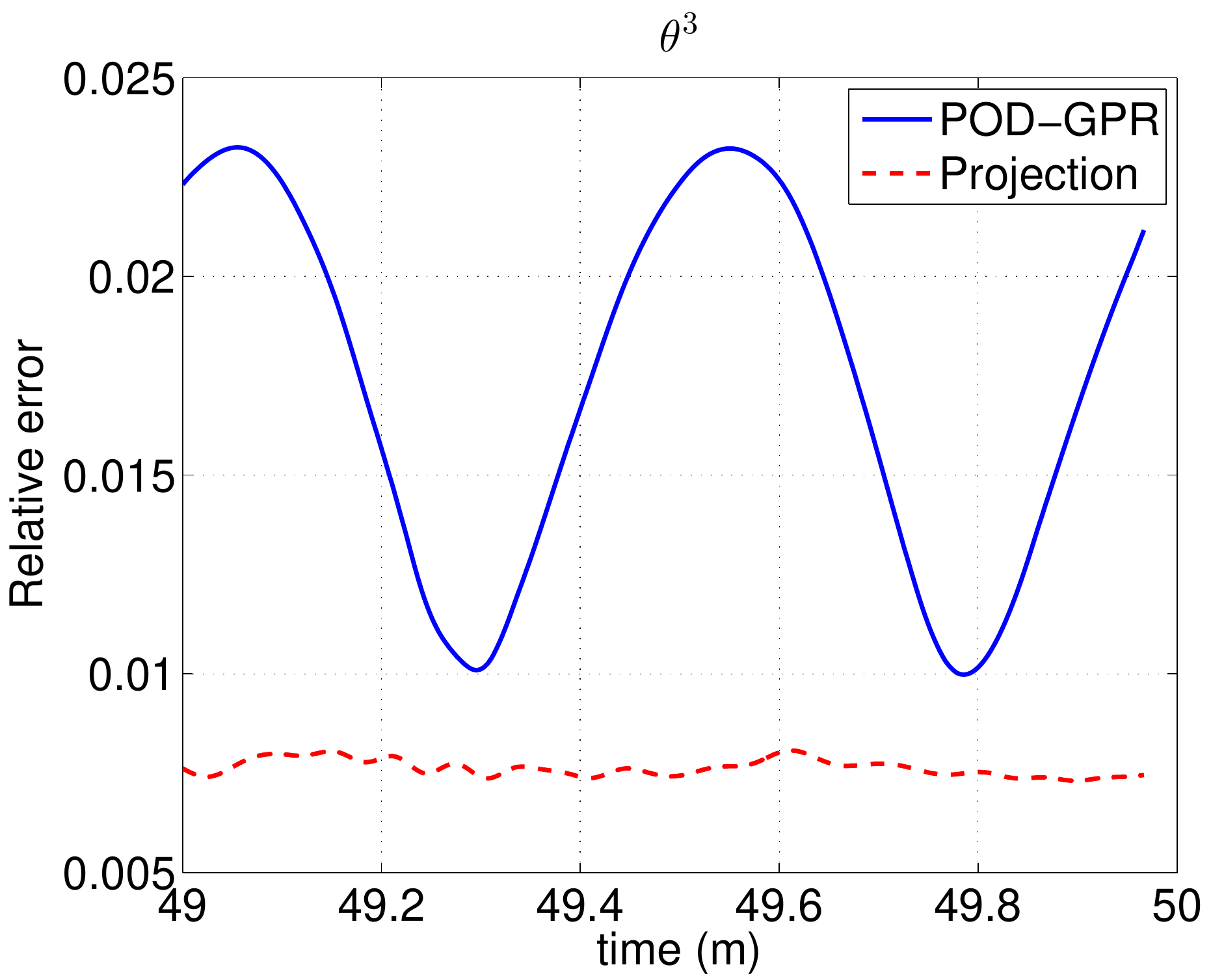}
}
\hspace{0.15cm}
\subfigure
{
\includegraphics[width=2.2in]{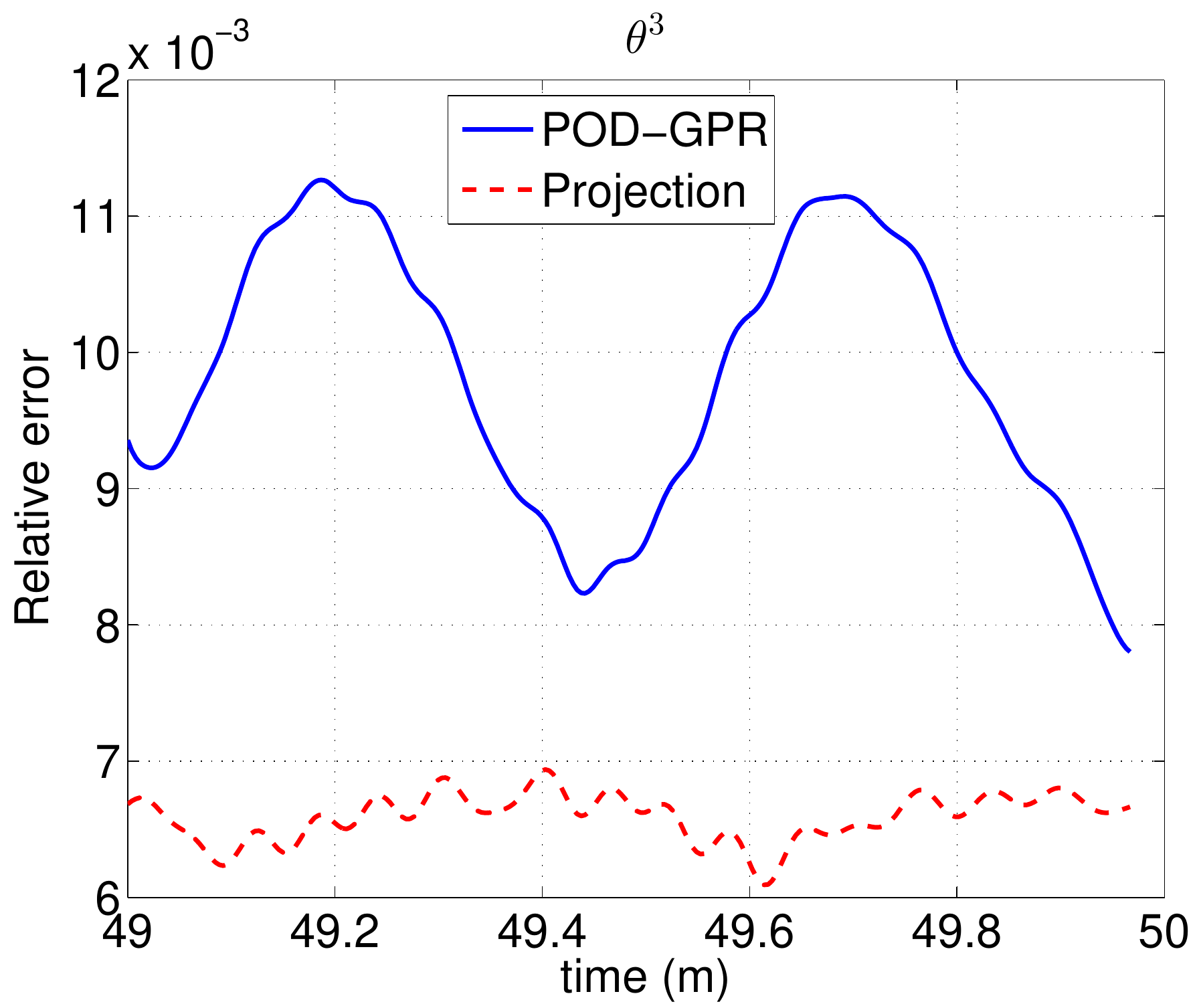}
}\\

\caption{Scattering of a plane wave by a multi-layer heterogeneous medium: Comparison of the relative $L^2$ error between POD-GPR and DGTD for $E_z$ (left) and $H_y$ (right) of four test points.}
\label{fig:16} 
\end{figure}

As can be observed in Table \ref{tab:5},
\begin{table}[]
\caption{Scattering of a plane wave by a multi-layer heterogeneous medium: Time performance comparison.}
\label{tab:5}
\centering
\begin{tabular}{cc}
\hline
Name & Time/s\\
\hline
Average DGTD solving time& $9.06\times10^{2}$\\
Average POD-GPR solving time& $4.37\times10^{0}$\\
GPR training time& $1.35\times10^{2}$\\
\hline
\end{tabular}
\end{table}
for these three test points, the average time consuming of POD-GPR is also much shorter than that of DGTD solver, implying an excellent time performance of our method. Finally, the time evolution of relative $L^2$ error that stem from POD-GPR is provided in Fig.\ref{fig:16}, 
resembling the one caused by projection in POD, by which the reliability of the POD-GPR approach is further confirmed. 
\section{Conclusion}
In this work, a non-intrusive POD-GPR model order reduction method combined with DGTD solver is introduced for the numerical simulation of parametric time-domain Maxwell's equations. The method apply DGTD to prepare a collection of full-order snapshots and extract a reduced basis through a two-step POD procedure.  The map between time/parameter value and the projection coefficients is approximated by GPR, which is represented as the combination of some time- and parameter-type GPs via SVDs. The offline stage completes the extraction of reduced basis and the training of GPRs, thus the online phase only do the output of GPR model, ensuring a full decoupling between offline and online phase and providing a fast and efficient tool for the parametric time-domain Maxwell's equations. Numerical results demonstrate the robustness and the high performance of this POD-GPR method. In the future, we will consider some more realistic 3D electromagnetic simulations and other type parameters such as frequencies, incident directions and some geometric parameters.

\end{document}